\documentclass[11pt,a4paper]{article}

\usepackage{color}
\usepackage{graphics,epsfig}
\usepackage{amsfonts,amssymb,amsmath,amsthm}

\newtheorem{lem}{Lemma}[section]
\newtheorem{cor}[lem]{Corollary}
\newtheorem{thm}[lem]{Theorem}
\newtheorem{prop}[lem]{Proposition}
\theoremstyle{definition}
\newtheorem{clm}[lem]{Claim}

\theoremstyle{remark}
\newtheorem{rem}[lem]{Remark}

\numberwithin{equation}{section}

\newcommand{\ep}{\varepsilon}
\newcommand{\ue}{u^\ep}
\newcommand{\ve}{v^\ep}
\newcommand{\fe}{f^\ep}
\newcommand{\fee}{f_2^{\;\!\ep}}
\newcommand{\g}{g^\ep}
\newcommand{\F}{f(u)-\ep \g(x,t,u)}
\newcommand{\tblow}{T^{max}}
\newcommand{\tmaxi}{t_{max}}
\newcommand{\de}{\delta}

\newcommand{\U}{U_0}
\newcommand{\Vep}{U_1 ^\ep}
\newcommand{\UU}{U_{0z}}
\newcommand{\UUU}{U_{0zz}}
\newcommand{\VV}{U_{1z}}
\newcommand{\VVep}{U_{1z}^\ep}
\newcommand{\VVV}{U_{1zz}}
\newcommand{\VVVep}{U_{1zz}^\ep}
\newcommand{\hm}{\alpha_-(\delta)}
\newcommand{\h}{a(\delta)}
\newcommand{\hp}{\alpha_+(\delta)}
\newcommand{\m}{\mu(\de)}
\newcommand{\mm}{\mu(\de)}
\newcommand{\n}{\nabla }
\newcommand{\emut}{e^{\mu t/\ep ^2}}
\newcommand{\emutt}{e^{\mu (\ep \mathcal G) t/\ep ^2}}
\newcommand{\dudn}{\displaystyle{\frac{\partial u}{\partial\nu}}}
\newcommand{\dom}{\partial\Omega}
\newcommand{\om}{\Omega}
\newcommand{\ombar}{\overline{\Omega}}
\newcommand{\edeux}{\displaystyle{\frac{1}{\ep^2}}}
\newcommand{\R}{\mathbb{R}}

\newcommand{\vsp}{\vspace{8pt}}

\newcommand{\di}{\displaystyle}

\newcommand{\f}{f_{\delta}}
\newcommand{\am}{\alpha _ -}
\newcommand{\ap}{\alpha _ +}
\newcommand{\apm}{\alpha _ \pm}
\newcommand{\EB}{e^{-\beta t/\ep ^2}}
\newcommand{\cha}{\gamma }
\newcommand{\chaep}{\gamma ^\ep}
\newcommand{\barg}{\bar{g}}
\newcommand{\tildeg}{\tilde{g}}
\newcommand{\gzero}{\Gamma _0}
\newcommand{\RD}{(RD^{\;\!\ep})}
\newcommand{\RDz}{(RD^{\;\! 0})}
\newcommand{\Pe}{(P^{\;\!\ep})}
\newcommand{\Pz}{(P^{\;\!0})}

\newcommand{\Pzggbar}{(P^{\;\!0}_{\barg})}
\newcommand{\gbragbar}{[ \barg]}
\newcommand{\gbragtilde}{[ \tildeg]}
\newcommand{\artificial}{d^{new}}
\newcommand{\phiunep}{\Phi ^\ep _1(v)}
\newcommand{\phiunzero}{\Phi ^0 _1(v)}
\newcommand{\phiep}{\Phi ^\ep(v)}
\newcommand{\phizero}{\Phi ^0(v)}

\title{The singular limit of the Allen-Cahn equation and
the FitzHugh-Nagumo system}
\author{ }
\date{}

\begin{document}
\maketitle \vspace{-15 mm}
\begin{center}

{\large\bf Matthieu Alfaro, Danielle Hilhorst,}\\[1ex]
Laboratoire de Math\'ematiques, Analyse Num\'erique et EDP, \\
Universit\'e de Paris Sud, 91405 Orsay Cedex, France, \\[2ex]

{\large\bf Hiroshi Matano}\\[1ex]
Graduate School of Mathematical Sciences, University of Tokyo,\\
3-8-1 Komaba, Tokyo 153-8914, Japan.\\[2ex]
\end{center}

\vspace{15pt}

\begin{abstract}
We consider an Allen-Cahn type equation of the form $u_t=\Delta
u+\ep^{-2}f^\ep(x,t,u)$, where $\ep$ is a small parameter and
$f^\ep(x,t,u)=f(u)-\ep g^\ep(x,t,u)$ a bistable nonlinearity
associated with a double-well potential whose well-depths can be
slightly unbalanced. Given a rather general initial data $u_0$
that is independent of $\ep$, we perform a rigorous analysis of
both the generation and the motion of interface. More precisely we
show that the solution develops a steep transition layer within
the time scale of order $\ep^2|\ln\ep|$, and that the layer obeys
the law of motion that coincides with the formal asymptotic limit
within an error margin of order $\ep$. This is an optimal estimate
that has not been known before for solutions with general initial
data, even in the case where $\g \equiv 0$.

Next we consider systems of reaction-diffusion equations of the
form
\[
\begin{cases}
 \,u_t=\Delta u+\ep^{-2}\, f^\ep(u,v)\\
 \,v_t=D \Delta v +h(u,v),
 \end{cases}\hspace{30pt}
\]
which include the FitzHugh-Nagumo system as a special case. Given
a rather general initial data $(u_0,v_0)$, we show that the
component $u$ develops a steep transition layer and that all the
above-mentioned results remain true for the $u$-component of these
systems.\\

\end{abstract}

%%%%%%%%%%%%%%%%%%%%%%%%%%%%%
\section{Introduction}\label{s:intro}

\subsection{Perturbed Allen-Cahn equation}

In some classes of nonlinear diffusion equations, solutions often
develop internal transition layers --- or ``interfaces"--- that
separate the spatial domain into different phase regions. This
happens, in particular, when the diffusion coefficient is very
small or the reaction term is very large. The motion of such
interfaces is often driven by their curvature. A typical example
is the Allen-Cahn equation $u_t=\Delta u+\ep^{-2}f(u)$, where
$\ep>0$ is a small parameter and $f(u)$ is a bistable
nonlinearity, whose meaning will be explained below. A usual
strategy for studying such phenomena is to first derive the
``sharp interface limit" as $\ep\to 0$ by a formal analysis, then
to check if this limit gives good approximation of the behavior of
actual layers.

In this paper we study a perturbed Allen-Cahn type equation of the
form
\[
 \Pe \quad\begin{cases}
 u_t=\Delta u+\edeux (\F) &\text{in }\om \times (0,+\infty)\\
 \dudn = 0 &\text{on }\partial \om \times (0,+\infty)\vspace{3pt}\\
 u(x,0)=u_0(x) &\text{in }\om,
 \end{cases}
\]
and study the behavior of layers near the sharp interface limit as
$\ep\to 0$. Here $\om$ is a smooth bounded domain in $\R^N$
($N\geq 2$) and $\nu$ is the Euclidian unit normal vector exterior
to $\partial \om$. The nonlinearity is given by $f(u):=-W'(u)$,
where $W(u)$ is a double-well potential with equal well-depth,
taking its global minimum value at $u=\apm$. More precisely we
assume that $f$ is $C^2$ and has exactly three zeros $\am<a<\ap$
such that
\begin{equation}\label{der-f}
f'(\apm)<0, \quad f'(a)>0\quad\ \hbox{(bistable nonlinearity)},
\end{equation}
and that
\begin{equation}\label{int-f}
\int _ {\am} ^ {\ap} f(u)\,du=0.
\end{equation}
The condition \eqref{der-f} implies that the potential $W(u)$
attains its local minima at $u=\am, \ap$, and \eqref{int-f}
implies that $W(\am)=W(\ap)$. In other words, the two stable zeros
of $f$, namely $\am$ and $\ap$, have ``balanced" stability. A
typical example is the cubic nonlinearity $f(u)= u(1-u^2)$.

The term $\ep\g$ represents a small perturbation, where
$\g(x,t,u)$ is a function defined on $\overline{\om}\times
[0,+\infty)\times\R$.  This has the role of breaking the balance
of the two stable zeros slightly. In the special case where
$\g\equiv 0$, problem $\Pe$ reduces to the usual Allen-Cahn
equation. As we will explain later, our main results are new even
for this special case.

We assume that $\g$ is $C^2$ in $x$ and $C^1$ in $t,\,u$, and
that, for any $T>0$ there exist $\vartheta \in (0,1)$ and $C>0$
such that, for all $(x,t,u) \in \overline{\om}\times[0,T]\times
\R$,
\begin{equation}\label{g-est1}
 |\Delta_{x} \g (x,t,u)| \leq C \ep^{-1} \quad \text{ and }
 \quad | \g_t(x,t,u)| \leq C \ep^{-1},
\end{equation}
\begin{equation}\label{g-est2}
 |\g_u (x,t,u)| \leq C,
\end{equation}
\begin{equation}\label{g-est3}
 \Vert \g(\cdot,\cdot,u)
 \Vert_{C^{1+\vartheta,\frac{1+\vartheta}{2}}(\ombar\times[0,T])}
 \leq C.
\end{equation}
Moreover, we assume that there exists a function $g(x,t,u)$ and a
constant, which we denote again by $C$, such that
\begin{equation}\label{g-est4}
 |\g (x,t,u)-g(x,t,u)| \leq C \ep,
\end{equation}
for all small $\ep>0$.  Note that the estimate \eqref{g-est3} and
the pointwise convergence $g^\ep\to g\;(\hbox{as }\ep\to 0)$ imply
that $g$ satisfies the same estimate as \eqref{g-est3}. For
technical reasons we also assume that
\begin{equation}\label{g-neumann}
\frac{\partial \g}{\partial \nu} =0 \quad \text { on }\dom\times
[0,T]\times \mathbb R,
\end{equation}
which, in turn, implies the same boundary condition for $g$. Apart
from these bounds and regularity requirements, we do not make any
specific assumptions on the perturbation term $\g$.

\vskip 6pt
\begin{rem}\label{rm:g0}
Since we will consider only bounded solutions in this paper, it is
sufficient to assume \eqref{g-est1}--\eqref{g-est3} to hold in
some bounded interval $-M\leq u\leq M$. Note that if $\g$ does not
depend on $\ep$, then the assumptions
\eqref{g-est1}--\eqref{g-est3} are automatically satisfied on any
interval $-M\leq u\leq M$. \hfill $\square$
\end{rem}

\vskip 6pt
\begin{rem}\label{rm:g1}
The reason why we do not assume more smoothness on $g$ is that we
will later apply our results to systems of equations including
FitzHugh-Nagumo system, in which $\g$ loses $C^{2,1}$-smoothness
as $\ep\to 0$. \hfill $\square$
\end{rem}

\vskip 6pt As for the initial data $u_0(x)$, we assume $u_0 \in
C^2(\ombar)$. Throughout the present paper the constant $C_0$ will
stand for the following quantity:
\begin{equation}\label{int1}
C_0:=\|u_0\|_{C^0(\ombar)}+\|\n u_0\|_{C^0(\ombar)}+\| \Delta
u_0\|_{C^0(\ombar)}.
\end{equation}
Furthermore we define the ``initial interface" $\Gamma _0$ by
\begin{equation}\label{initial-interface}
\Gamma_0:=\{x\in\om, \; u_0(x)=a\},
\end{equation}
and suppose that $\Gamma_0$ is a $C^{3+\vartheta}$ hypersurface
without boundary such that, $n$ being the outward unit normal
vector to $\Gamma_0$,
\begin{equation}\label{dalltint}
\Gamma_0 \subset\subset \Omega \quad \mbox { and } \quad \n
u_0(x)\cdot n(x) \neq 0\quad\text{if $x\in\Gamma_0,$}
\end{equation}
\begin{equation}\label{initial-data}
u_0>a \quad \text { in }\om ^+ _0,\quad u_0<a \quad \text { in
}\om ^- _0 ,
\end{equation}
where $\om ^- _0$ denotes the region enclosed by $\Gamma _0$ and
$\om ^+ _0$ the region enclosed between $\partial \om$ and
$\Gamma_0$.

It is standard that problem $\Pe$ has a unique smooth solution,
which we denote by $\ue$. As $\ep \rightarrow 0$, a formal
asymptotic analysis shows the following: in the very early stage,
the diffusion term $\Delta u$ is negligible compared with the
reaction term $\ep ^{-2}(f(u)-\ep \g(x,t,u))$ so that, in the
rescaled time scale $\tau=t /\ep ^2$, the equation is well
approximated by the ordinary differential equation $u_\tau
=f(u)+O(\ep)$. Hence, in view of the profile of $f$, the value of
$u^{\ep}$ quickly becomes close to either $\ap$ or $\am$ in most
part of $\Omega$, creating a steep interface (transition layer)
between the regions $\{\ue\approx \am\}$ and $\{\ue\approx \ap\}$.
Once such an interface develops, the diffusion term becomes large
near the interface, and comes to balance with the reaction term.
As a result, the interface ceases rapid development and starts to
propagate in a much slower time scale.

To study such interfacial behavior, it is useful to consider a
formal asymptotic limit of $\Pe$ as $\ep\rightarrow 0$. Then the
limit solution $\tilde u (x,t)$ will be a step function taking the
value $\ap$ on one side of the interface, and $\am$ on the other
side. This sharp interface, which we will denote by $\Gamma_t$,
obeys a certain law of motion, which is expressed as follows (see
Section \ref{s:formal} for details):
\[
 \Pz\quad\begin{cases}
 \, V_{n}=-(N-1)\kappa + c_0 (G(x,t,\ap)-G(x,t,\am))
 \quad \text { on } \Gamma_t \vspace{3pt}\\
 \, \Gamma_t\big|_{t=0}=\Gamma_0,
\end{cases}
\]
where $V_n$ is the normal velocity of $\Gamma _t$ in the exterior
direction, $\kappa$ the mean curvature at each point of
$\Gamma_t$,
\begin{equation}\label{c0}
c_0=\Big [\sqrt 2\int_{\am}^{\ap} (W(s) - W(\am))^{1/2} ds \Big ]
^{-1},
\end{equation}
$$
W(s) =-\int_a^s f(r) dr,\quad \quad  G(x,t,s)=\int _a^s
g(x,t,r)dr.
$$
It is well known that problem $\Pz$ possesses locally in time a
unique smooth solution. Let $0\leq t <\tblow$, $\tblow \in
(0,+\infty]$, be the maximal time interval for the existence of
the solution of $\Pz$ and denote this solution by $\Gamma=\bigcup
_{0\leq t < \tblow} (\Gamma_t\times\{t\})$. Hereafter, we fix $T$
such that $0<T<\tblow$ and work on $[0,T]$. More precisely, so as
$g(\cdot,\cdot,u)$, the function $G(\cdot,\cdot,u)$ is of class
$C^{1+\vartheta,\frac{1+\vartheta}2}$, which implies, by the
standard theory of parabolic equations, that $\Gamma $ is of class
$C^{3+\vartheta,\frac{3+\vartheta}2}$. For more details, we refer
to \cite{CR}, Lemma 2.1.

Next we set
\[
Q_T:=\om \times (0,T),
\]
and, for each $t\in [0,T]$, we denote by $\om^-_t$ the region
enclosed by the hypersurface $\Gamma_t$, and by $\om^+_t$ the
region enclosed between $\partial \om$ and $\Gamma_t$. We define a
step function $\tilde u(x,t)$ by
\begin{equation}\label{u}
\tilde u(x,t)=\begin{cases}
\, \ap &\text{in } \om^+_t\\
\, \am &\text{in } \om^-_t
\end{cases} \quad\text{for } t\in[0,T],
\end{equation}
which represents the formal asymptotic limit of $\ue$ (or the {\it
sharp interface limit}) as $\ep\to 0$.

The aim of the present paper is to make a detailed study of the
limiting behavior of the solution $u^\ep$ of problem $\Pe$ as
$\ep\to 0$. Our first main result, Theorem \ref{width}, describes
the profile of the solution after a very short initial period. It
asserts that: given a virtually arbitrary initial data $u_0$, the
solution $\ue$ quickly becomes close to $\apm$, except in a small
neighborhood of the initial interface $\Gamma _0$, creating a
steep transition layer around $\Gamma _0$ ({\it generation of
interface}). The time needed to develop such a transition layer,
which we will denote by $t ^\ep$, is of order $\ep^2|\ln\ep|$. The
theorem then states that the solution $\ue$ remains close to the
step function $\tilde u$ on the time interval $[t^\ep,T]$ ({\it
motion of interface}); in other words, the motion of the
transition layer is well approximated by the limit interface
equation $\Pz$.

\begin{thm}[Generation and motion of interface]\label{width}
Let $\eta$ be an arbitrary constant satisfying $0< \eta <
\min(a-\am,\ap -a)$ and set
$$
\mu=f'(a).
$$
Then there exist positive constants $\ep _0 $ and $C$ such that,
for all $\,\ep \in (0,\ep _0)$ and for all $\,t^\ep \leq t \leq
T$, where $t^\ep:=\mu ^{-1} \ep ^2 |\ln \ep|$, we have
\begin{equation}\label{resultat}
\ue(x,t) \in
\begin{cases}
\,[\am-\eta,\ap+\eta]\quad\text{if}\quad
x\in\mathcal N_{C\ep}(\Gamma_t)\\
\,[\am-\eta,\am+\eta]\quad\text{if}\quad
x\in\om_t^-\setminus\mathcal N_{C\ep}(\Gamma
_t)\\
\,[\ap-\eta,\ap+\eta]\quad\text{if}\quad x\in\om
_t^+\setminus\mathcal N_{C\ep}(\Gamma _t),
\end{cases}
\end{equation}
where $\mathcal N _r(\Gamma _t):=\{x\in \om, dist(x,\Gamma
_t)<r\}$ denotes the $r$-neighborhood of $\Gamma _t$.
\end{thm}

\begin{cor}[Convergence]\label{total}
As $\ep\to 0$, $\ue$ converges to $\tilde u$ everywhere in
$\bigcup _{0<t\leq T}(\om^\pm_t\times\{ t\})$.
\end{cor}

The next theorem is concerned with the relation between the actual
interface $\Gamma^\ep _t:=\{x \in \om,\; u^\ep(x,t)=a\}$ and the
formal asymptotic limit $\Gamma _t$, which is given as the
solution of $\Pz$.

\begin{thm}[Error estimate]\label{error}
There exists $C>0$ such that
\begin{equation}\label{thi-2}
\Gamma _t ^\ep \subset \mathcal N _{C \ep} (\Gamma _t)\quad \text
{ for }\ 0 \leq t \leq T.
\end{equation}
\end{thm}

\begin{cor}[Convergence of interface]\label{total-2}
There exists $C>0$ such that
\begin{equation}\label{thi-3}
d_\mathcal H (\Gamma ^\ep _t,\Gamma _t)\leq C\ep \quad \text { for
}\ 0\leq t \leq T,
\end{equation}
where $d_\mathcal H (A,B):=\max \{ \sup_{a \in A}
d(a,B),\,\sup_{b\in B} d(b,A)\}$ denotes the Hausdorff distance
between two compact sets $A$ and $B$.  Consequently, $\Gamma^\ep_t
\to \Gamma_t$ as $\ep\to 0$ uniformly in $0\leq t\leq T$, in the
sense of Hausdorff distance.
\end{cor}

Note that the estimates \eqref{thi-2} and \eqref{thi-3} follow
from Theorem \ref{width} in the range $t^\ep \leq t \leq T$, but
the range $0\leq t \leq t^\ep$ has to be treated by a separate
argument since the behavior of the solution in this time range is
quite different from that of the later stage.

The estimate \eqref{resultat} in our Theorem \ref{width} implies
that, once a transition layer is formed, its thickness remains
within order $\ep$ for the rest of time.  Here, by ``thickness of
interface" we mean the smallest $r>0$ satisfying
\[
\{\,x\in\Omega,\;\ue(x,t)\not\in[\am-\eta,\,\am+\eta] \cup
[\ap-\eta,\,\ap+\eta]\,\} \subset {\mathcal N}_r (\Gamma^\ep_t).
\]
Naturally this quantity depends on $\eta$, but the estimates
\eqref{resultat} and \eqref{thi-3} assert that it is bounded by
$2C\ep$ (with the constant $C$ depending on $\eta$) regardless of
the choice of $\eta>0$.

\vskip 6pt
\begin{rem}[Optimality of the thickness estimate]
\label{rm:thickness} The above $O(\ep)$ estimate is optimal, {\it
i.e.} the interface cannot be thinner than this order.  In fact,
rescaling time and space as $\tau:=t/\ep^2,$ $y:=x/\ep$, we get
\[
u_{\tau}=\Delta_y \vspace{5pt}u+f(u)-\ep\,\g.
\]
Thus, by the uniform boundedness of $u$ and by standard parabolic
estimates, we have $|\nabla_y\vspace{4pt} u|\leq M$ for some
constant $M>0$, which implies
\[
|\nabla_x \vspace{4pt}u(x,t)|\leq \frac{M}{\ep}.
\]
{}From this bound it is clear that the thickness of interface
cannot be smaller than $M^{-1}(\ap-\am)\,\ep$, hence, by
\eqref{resultat}, it has to be exactly of order $\ep$.
Intuitively, this $O(\ep)$ estimate follows also from the formal
asymptotic expansion \eqref{inner}, but the validity of such an
expansion is far from obvious for solutions with arbitrary initial
data. \hfill $\square$
\end{rem}
\vskip 6pt

Our $O(\ep)$ estimate is new, even in the special case where
$\g\equiv 0$, provided that $N\geq 2$.  Previously, the best
thickness estimate in the literature was of order $\ep |\ln \ep|$
(see \cite{C1}), except that X. Chen has recently obtained an
$O(\ep)$ estimate for the case $N=1$ by a different argument
(private communication). We also refer to the forthcoming papers
\cite{HMS} and \cite{HKMN}, in which the same $O(\ep)$ estimate is
established for different but related problems.  The paper
\cite{HMS} is concerned with a ``balanced type" Allen-Cahn
equation with large spatial inhomogeneity, namely an equation of
the form $u_t = \nabla (k(x) \nabla u) + \ep^{-2} h(x)f(u)$, and
\cite{HKMN} is concerned with a Lotka-Volterra
competition-diffusion system with large spatial inhomogeneity
whose nonlinearity is of the balanced bistable type.

\vskip 6pt
\begin{rem}[Optimality of the generation time]
\label{rm:time} The estimate \eqref{resultat} also implies that
the generation of interface takes place within the time span of
$t^\ep$. This estimate is optimal.  In other words, a
well-developed interface cannot appear much earlier; see
Proposition \ref{pr:optimal-time} for details. \hfill $\square$
\end{rem}

\vskip 6pt The singular limit of the Allen-Cahn equation was first
studied in the pioneering work of Allen and Cahn \cite{AC} and,
slightly later, in Kawasaki and Ohta \cite{KO} from the point of
view of physicists. They derived the interface equation by formal
asymptotic analysis, thereby revealing that the interface moves by
the mean curvature. These early observations triggered a flow of
mathematical studies aiming at rigorous justification of the above
limiting procedure; see, for example, \cite{BS}, \cite{C1,C2} and
\cite{MS1, MS2} for results in the framework of classical
solutions, and \cite{ESS}, \cite{BSS,BaS} and \cite{I} for the
case where $\Gamma_t$ is a viscosity solution of the interface
equation.

As for problem $\Pe$, whose nonlinearity is slightly unbalanced,
the limit interface equation involves a pressure term as well as
the curvature term as indicated in $\Pz$.  This fact has been long
known on a formal level; see e.g. \cite{RSK}. Ei, Iida and
Yanagida \cite{EIY} proved rigorously that the motion of the
layers of $\Pe$ is well approximated by the limit interface
equation $\Pz$, on the condition that the initial data has already
a well developed transition layer. In other words, they studied
the motion of interface, but not the generation of interface.

\subsection{Singular limit of reaction-diffusion systems}

Our results can be extended to reaction-diffusion systems of the
form
\[
 \RD \quad
 \begin{cases}
 \,u_t=\Delta u+\edeux\, \fe(u,v) \quad&\textrm{in}\ \
 \Omega\times(0,+\infty)\vspace{2pt}\\
 \,v_t=D \Delta v +h(u,v) \quad&\textrm{in} \ \
 \Omega\times(0,+\infty)\vspace{4pt}\\
 \,\dudn =\displaystyle{\frac{\partial v}{\partial\nu}}
 = 0 \quad &\textrm{on}\ \
 \partial\Omega\times(0,+\infty) \vspace{4pt}\\
 \,u(x,0)=u_0(x)\quad &\textrm{in}\ \ \Omega\vspace{3pt}\\
 \,v(x,0)=v_0(x)\quad &\textrm{in}\ \ \Omega,
 \end{cases}\hspace{30pt}
\]
where $D$ is a positive constant, and $f^\ep,\;h$ are $C^2$
functions such that
\begin{itemize}
\item[{\bf (F)}] there exist $C^2$ functions
$f_1(u,v),\;\fee(u,v)$ such that
\begin{equation}\label{f(u,v)}
 \fe(u,v)=f(u)+\ep f_1(u,v)+\ep^2 \fee(u,v),
\end{equation}
where $f(u)$ is a bistable nonlinearity satisfying \eqref{der-f},
\eqref{int-f}, and $\fee$, along with its derivatives in $u,v$,
remain bounded as $\ep\to 0$; \item[{\bf (H)}] for any constant
$L,M>0$ there exists a constant $M_1\geq M$ such that
\begin{equation}\label{h(u,v)}
 h(u,-M_1)\geq 0\geq h(u,M_1) \qquad\hbox{for}\ \ |u|\leq L.
\end{equation}
\end{itemize}

The conditions (F) and (H) imply that the ODE system
\[
\dot u =\frac{1}{\ep^2}\,\fe(u,v),\quad\ \dot v= h(u,v)
\]
has a family of invariant rectangles of the form $\{|u|\leq
L,\,|v|\leq M\}$, provided that $\ep$ is sufficiently small. The
maximum principle and standard parabolic estimates then guarantee
that every solution of $\RD$ exists globally for $t\geq 0$ and
remains bounded as $t\to\infty$ (see Proposition \ref{pr:global}).
Apart from \eqref{h(u,v)}, we do not make any specific assumptions
on the function $h$.

Problem $\RD$ represents a large class of important
reaction-diffusion systems including the FitzHugh-Nagumo system
\begin{equation}\label{FN}
 \begin{cases}
 \,u_t=\Delta u+\edeux (f(u)-\ep v) \vspace{4pt}\\
 \,v_t=D \Delta v +\alpha u - \beta v ,
 \end{cases}\hspace{50pt}
\end{equation}
which is a simplified model for nervous transmission, and the
following type of prey-predator system:
\begin{equation}\label{PP}
 \begin{cases}
 \,u_t=\Delta u+\edeux \big((1-u)(u-1/2)-\ep v\big)u\vspace{4pt}\\
 \,v_t=D \Delta v +(\alpha u - \beta v)v.
 \end{cases}\hspace{20pt}
\end{equation}

\vskip 6pt
\begin{rem}\label{rm:positive}
In some equations such as the prey-predator system \eqref{PP},
only nonnegative solutions are to be considered.  In such a case,
we replace the condition \eqref{h(u,v)} by
\[
 h(u,0)\geq 0\geq h(u,M_1) \qquad
 \hbox{for}\ \ 0\leq u\leq L,
\]
and assume $f^\ep(0,v)\geq 0$. The rest of the argument remains
the same. \hfill$\square$
\end{rem}

\vskip 6pt Now the same formal analysis as is used to derive $\Pz$
shows that the singular limit of $\RD$ as $\ep\to 0$ is the
following moving boundary problem:
\[
\RDz\quad
\begin{cases}
\,V_{n}=-(N-1)\kappa - c_0\,F_1(\tilde{v}(x,t))
\quad &\textrm{on}\ \ \Gamma_t \vspace{4pt}\\
\,\tilde{v}_t=D \Delta \tilde{v} + h(\tilde{u},\tilde{v})
\quad &\textrm{in}\ \ \Omega\times(0,+\infty) \vspace{4pt}\\
\,\Gamma_t\big|_{t=0}=\Gamma_0 \vspace{4pt}\\
\,\displaystyle{\frac{\partial \tilde{v}}{\partial \nu}} = 0
\quad&\textrm{on}\ \ \partial \Omega \times (0,+\infty)
\vspace{4pt}\\
\,\tilde v(x,0)=v_0(x) \quad&\textrm{in}\ \ \Omega,
\end{cases}
\]
where $\tilde{u}$ is the step function defined in \eqref{u} and
\[
F_1(v)=\int_{\am}^{\ap}f_1(r,v)\,dr.
\]
This is a system consisting of an equation of surface motion and a
partial differential equation. Since $\tilde{u}$ is determined
straightforwardly from $\Gamma_t$, in what follows, by a solution
of $\RDz$ we mean the pair $(\Gamma,\tilde{v}):=(\Gamma _t,
\tilde{v}(x,t))$. In the case of the FitzHugh-Nagumo system
\eqref{FN}, $\RDz$ reduces to
\[
\begin{cases}
\,V_{n}=-(N-1)\kappa + c_0\,(\ap-\am)\tilde{v}(x,t)\vsp\\
\,\tilde{v}_t=D \Delta \tilde{v} + \alpha \tilde u-\beta \tilde v,
\end{cases}
\]
while in the prey-predator system \eqref{PP}, $\RDz$ reduces to
\[
\begin{cases}
\,V_{n}=-(N-1)\kappa + c_0\,\tilde{v}(x,t)/2\vsp\\
\,\tilde{v}_t=D \Delta \tilde{v} + (\alpha \tilde u-\beta \tilde
v)\tilde u.
\end{cases}
\]
Note that the positive sign in front of the term $c_0\tilde v
(x,t)$ in the interface equation implies an inhibitory effect on
$\tilde u$, since the velocity $V_n$ is measured in the exterior
normal direction, toward which $\tilde u$ decreases.

\begin{lem}[Local existence]\label{existence}
Assume that $v_0 \in C^2(\overline{\Omega})$ and that $\Gamma_0$
is a $C^2$ hypersurface which is the boundary of a domain $D_0
\subset \subset \Omega$. Then there exists $\tblow \in
(0,+\infty]$ such that the limit free boundary problem $\RDz$ has
a unique solution $(\Gamma,\tilde v)$ in the interval
$[0,\tblow)$.
\end{lem}

This existence result was established in \cite{CXY}. The
uniqueness can be obtained by using the estimates in \cite{C2}.

Hereafter, we fix $T$ such that $0<T<\tblow$ and work on $[0,T]$.
Our main results for the system $\RD$ are the following:

\begin{thm}[Thickness of interface]\label{width-N}
Let \eqref{f(u,v)} and \eqref{h(u,v)} hold (or let the assumptions
in Remark \ref{rm:positive} hold). Assume also that $u_0$
satisfies \eqref{dalltint} and \eqref{initial-data}. Then the same
conclusion as in Theorem \ref{width} holds for $\RD$.
\end{thm}

\begin{cor}[Convergence]\label{total-N}
Under the assumptions of Theorem \ref{width-N}, the same
conclusion as in Corollary \ref{total} holds for $\RD$.
\end{cor}

\begin{thm}[Error estimate]\label{error-N}
Let the assumptions of Theorem \ref{width-N} hold. Then the same
conclusion as in Theorem \ref{error} holds for $\RD$. Moreover,
there exists a constant $C>0$ such that
$$
\Vert\ve - \tilde v\Vert_{ L^\infty (\om\times(0,T))} \leq C\ep.
$$
\end{thm}

\begin{cor}[Convergence of interface]\label{total-2-N}
Under the assumptions of Theorem \ref{width-N}, the same
conclusion as in Corollary \ref{total-2} holds for $\RD$.
\end{cor}

The organization of this paper is as follows. In Section
\ref{s:formal}, we derive the interface equation $\Pz$ from $\Pe$
by formal asymptotic expansions which involve the so-called signed
distance function.

In Sections \ref{s:generation-0} and \ref{s:generation-g}, we
present basic estimates concerning the generation of interface for
$\Pe$. For the clarity of underlying ideas, we first consider the
special case where $\g\equiv 0$ in Section \ref{s:generation-0},
and deal with the general case in Section \ref{s:generation-g}.

In Section \ref{s:motion} we prove a preliminary result on the
motion of interface (Lemma \ref{fix}), which implies that if the
initial data has already a well-developed transition layer, then
the layer remains to exist for $0\leq t\leq T$ and its motion is
well approximated by the interface equation $\Pz$.

Our approach in Sections \ref{s:generation-0} to \ref{s:motion} is
based on the sub- and super-solution method, but we use two
completely different sets of sub- and super-solutions. More
precisely, the sub- and super-solutions for the motion of
interface are constructed by using the first two terms of the
formal asymptotic expansion \eqref{inner}, while those for the
generation of interface are constructed by modifying the solution
of the equation in the absence of diffusion: $u_t=\ep ^{-2}f(u)$.

In Section \ref{s:proof}, we prove our main results for $\Pe$:
Theorems \ref{width}, \ref{error} and their respective
corollaries.

In the final section, we study the reaction-diffusion system $\RD$
and prove Theorems \ref{width-N}, \ref{error-N} and their
corollaries. These results are obtained by applying a slightly
modified version of the results for $\Pe$. The strategy is to
regard $\fe(u,v)$ as a perturbation of $f(u)$. Indeed, the
equation for $u$ in $\RD$ is identical to $\Pe$ if we set
$\g=-f_1-\ep\fee$. However, what makes the analysis difficult is
the fact that $\g$ is no longer a given function but a quantity
that depends on the unknown function $\ve$. In particular, the
existence of the limit $\g\to g\;(\ep\to 0)$ is not a priori
guaranteed, and the estimate \eqref{g-est4} is far from obvious.
As it turns out, the standard $L^p$ or Schauder estimates for
$\ve$ would not yield \eqref{g-est4}, because of the fact that
$\ue$ converges to a discontinuous function as $\ep\to 0$. In
order to overcome this difficulty, we derive a fine estimate of
$\ve$ that is based on estimates of the heat kernel and the fact
that $\ue$ remains uniformly smooth outside of an $O(\ep)$
neighborhood of the smooth hypersurface $\Gamma_t$.

\section{Formal derivation of the interface motion equation}
\label{s:formal}

In this section we derive the equation of interface motion
corresponding to problem $\Pe$ by using a formal asymptotic
expansion. The resulting interface equation can be regarded as the
singular limit of $\Pe$ as $\ep \to 0$. Our argument is basically
along the same lines with the formal derivation given by Nakamura,
Matano, Hilhorst and Sch\"atzle \cite{NMHS}, who studied a similar
but slightly different type of spatially inhomogeneous equations
by formal analysis.  Let us also mention some earlier papers
\cite{ABC}, \cite{Fi} and \cite{RSK} involving the method of
matched asymptotic expansions for problems that are related to
ours.

As in \cite{NMHS}, the first two terms of the asymptotic expansion
determine the interface equation.  Though our analysis in this
section is for the most part formal, the observations we make here
will help the rigorous analysis in later sections.\\

Let $u^\ep$ be the solution of $\Pe$. We recall that $\Gamma _t
^\ep:=\{x \in \om, u^\ep(x,t)=a\}$ is the interface at time $t$
and call $\Gamma ^\ep:=\bigcup _{t \geq 0} (\Gamma _t ^\ep \times
\{t\})$ the interface. Let $\Gamma=\bigcup _{0\leq t \leq
T}(\Gamma_t\times\{t\})$ be the unique solution of the limit
geometric motion problem $\Pz$ and let $\widetilde d$ be the
signed distance function to $\Gamma$ defined by:
\begin{equation}\label{eq:dist}
\widetilde d (x,t)=
\begin{cases}
&\hspace{-10pt}\mbox{dist}(x,\Gamma _t)\quad\text{for }x\in\om_t^+ \\
-&\hspace{-10pt} \mbox{dist}(x,\Gamma _t) \quad \text{for }
x\in\om _t^- ,
\end{cases}
\end{equation}
where $\mbox{dist}(x,\Gamma _t)$ is the distance from $x$ to the
hypersurface  $\Gamma _t$ in $\om$. We remark that $\widetilde
d=0$ on $\Gamma$ and that $|\nabla \widetilde d|=1$ in a
neighborhood of $\Gamma.$ We then define
\begin{equation*}
Q^+_T = \bigcup_{\,0<t\leq T}(\Omega^+_t \times\{t\}),\qquad Q^-_T
= \bigcup_{\,0<t\leq T}(\Omega^-_t \times\{t\}).
\end{equation*}
We also assume that the solution $\ue$ is of the form
\begin{equation} \label{outer}
\ue(x,t)= \apm  + \ep u_1^\pm (x,t) + \cdots
\quad\ \hbox{in}\ \ Q^{\pm}_T
\end{equation}
away from the interface $\Gamma$ (the outer expansion), and
\begin{equation}\label{inner}
\ue(x,t)=\U(x,t,\xi)+\ep U_1(x,t,\xi)+\cdots
\end{equation}
near $\Gamma$ (the inner expansion), where $U_j(x,t,z)$,
$j=0,1,\cdots$, are defined for $x\in \overline \Omega$, $t\geq
0$, $z\in \R$ and $\xi:=\widetilde d(x,t)/\ep$. The stretched
space variable $\xi$ gives exactly the right spatial scaling to
describe the rapid transition between the regions $\{\ue \approx
\am\}$ and $\{\ue \approx \ap\}$. We normalize $U_0$ in such a way
that
$$
\U(x,t,0)=a
$$
(normalization conditions). To make the inner and outer expansions
consistent, we require that
\begin{equation}\label{match}
\begin{array}{ll}
\U(x,t,+\infty)= \ap, \quad \U(x,t,-\infty)= \am.
\end{array}
\end{equation}
As we will see below this will determine $U_0$ uniquely, which will
then determine $U_1$.

In what follows we will substitute the inner expansion
\eqref{inner} into the parabolic equation of problem $\Pe$ and
collect the $\ep^{-2}$ and $\ep^{-1}$ terms. To that purpose we
compute the needed terms and get
\[
\begin{array}{l}
\di\ue_t= U_{0t}+\UU \frac{\widetilde d_t}{\ep}
+\ep U_{1t}+\VV \widetilde d _t+\cdots \vsp \\
\di \nabla \ue= \nabla  \U+ \UU \frac{\nabla\widetilde d}{\ep}
+\ep \nabla  U_1 + \VV\nabla \widetilde d+\cdots \vsp \\
\di \Delta  \ue=\Delta\U+2 \frac{\nabla \widetilde d}{\ep}\cdot
\nabla \UU+\UU \frac{\Delta \widetilde d}{\ep}+
\UUU \frac{|\nabla \widetilde d|^2}{\ep^2}+\ep \Delta U_1 \vsp \\
\di\hspace{50pt}+2\nabla \widetilde d \cdot \nabla \VV + \VV
\Delta
\widetilde d +\VVV \frac{|\nabla\widetilde d|^2}{\ep}+\cdots\vsp \\
\di f( \ue)=f(\U)+\ep f'(\U)U_1+O(\ep^2) \vsp \\
\di \g(x,t,\ue)=g(x,t,\ue)+O(\ep) \qquad \big(\ \longleftarrow \text{ in view of \eqref{g-est4} }\big)\vsp\\
\di\hspace{56pt}=g(x,t,U_0)+O(\ep),\end{array}
\]
where the functions $U_i$ $(i = 0, 1)$, as well as their
derivatives, are taken at the point $(x,t,\widetilde d
(x,t)/\ep)$. Note also that $\nabla$ and $\Delta$ stand for
$\nabla _x$ and $\Delta _x$, respectively. Collecting the $\ep
^{-2}$ terms yields
$$
\begin{array}{ll}
\UUU + f(\U) &=0.
\end{array}
$$
In view of the normalization and matching conditions, we can now
assert that $\U(x,t,z)=U_0(z)$, where $U_0(z)$ is the unique
solution of the stationary problem
\begin{equation}\label{eq-phi}
\left\{\begin{array}{ll}
{U_0} '' +f(U_0)=0 \vspace{3pt}\\
U_0(-\infty)= \am,\quad U_0(0)=a,\quad U_0(+\infty)= \ap .
\end{array} \right.
\end{equation}
This solution represents the first approximation of the profile of
a transition layer around the interface observed in the stretched
coordinates. Note that the integral condition \eqref{int-f}
guarantees the existence of a solution of \eqref{eq-phi}. For
example, in the simple case where $f(u)=u(1-u^2)$, we have
$U_0(z)=\tanh (z /\sqrt 2)$. In the general case, the following
standard estimates hold:

\begin{lem}\label{est-phi}
There exist positive constants $C$ and $\lambda$ such that
$$
\begin{array}{ll}
0 <\ap-U_0(z)&\leq Ce^{-\lambda|z|} \quad \text{ for } z\geq 0
\vspace{3pt}\\
0 <U_0(z)-\am&\leq Ce^{-\lambda|z|} \quad \text{ for } z\leq 0.\\
\end{array}
$$
In addition, $U_0$ is a strictly increasing function and, for
$j=1, 2$,
\begin{equation}\label{diff-U0}
|D^jU_0(z)|\leq Ce^{-\lambda|z|} \quad \text{ for } z\in \R.
\end{equation}
\end{lem}

{\noindent \bf Proof.} We only give an outline. Rewriting the
equation in \eqref{eq-phi} as
\[
\dot{u}=v, \qquad \dot{v}=-f(u),
\]
we see that $(U_0(z),U_0'(z))$ is a heteroclinic orbit of the
above system connecting the equilibria $(\am,0)$ and $(\ap,0)$.
These equilibria are saddle points, with the linearized
eigenvalues $\{\lambda_-,\,-\lambda_-\}$ and
$\{\lambda_+,\,-\lambda_+\}$, respectively, where
\[
\lambda_-=\sqrt{-f'(\am)},\qquad \lambda_+=\sqrt{-f'(\ap)}.
\]
Consequently, we have
\begin{equation}\label{U0-asympt}
U_0(z)=\left\{
\begin{array}{ll}
\am+C_1\,e^{\,\lambda_- z}+o(e^{\,\lambda_- z})
\qquad&\hbox{as} \ z\to -\infty,\\
\ap+C_2\,e^{-\lambda_+ z}+o(e^{-\lambda_+ z}) \qquad&\hbox{as} \
z\to +\infty,
\end{array}\right.
\end{equation}
for some constants $C_1, C_2$.  The desired estimates now follow
by setting $\lambda=\min (\lambda_+,\,\lambda_-)$. \qed
\\

Next we collect the $\ep ^{-1}$ terms.  Recalling that $\n
U_{0z}=0$ and that $|\nabla\tilde{d}|=1$ near $\Gamma_t$, we get
\begin{equation}\label{eqU1}
\VVV + f'(\U)U_1 ={U_0}'(\widetilde d_{t}-\Delta \widetilde d)+
g(x,t,\U).
\end{equation}
This equation can be seen as a linearized problem for
\eqref{eq-phi} with an inhomogeneous term. As is well-known (see,
for instance, \cite{NMHS}), the solvability condition for the
above equation plays the key role in determining the equation of
interface motion. The following lemma is rather standard, but we
give an outline of the proof for the convenience of the reader.

\begin{lem}[Solvability condition]\label{Fredholm}
Let $A(z)$ be a bounded function on $-\infty<z<\infty$. Then the
problem
\begin{equation}\label{eq-psi}
\left\{\begin{array}{l}
\psi_{zz} + f'(U_0(z))\psi=A(z),\qquad z\in \R \vsp \\
\psi(0)=0,\quad \psi \in L^\infty (\R),
\end{array} \right.
\end{equation}
has a solution if and only if
\begin{equation}\label{cond-sol}
\int_\R A(z){U_0} '(z)dz=0.
\end{equation}
Moreover the solution, if it exists, is unique and satisfies
\begin{equation}\label{psi-bound}
|\psi(z)|\leq C\Vert A\Vert_{L^\infty} \qquad\hbox{for} \ z\in\R ,
\end{equation}
for some constant $C>0$.
\end{lem}

{\noindent \bf Proof.} Multiplying the equation by ${U_0}'$ and
integrating it by parts, we easily see that the condition
\eqref{cond-sol} is necessary. Conversely, suppose that this
condition is satisfied. Then, since ${U_0} '$ is a bounded
positive solution to the homogeneous equation $\psi_{zz} +
f'(U_0(z))\psi=0$, one can use the method of variation of
constants to find the above solution $\psi$ explicitly.  More
precisely,
\begin{equation}\label{psi}
\begin{array}{rl}
\psi(z)&\di =\varphi(z)\int_0^z\Big(\varphi^{-2}(\zeta)
\int_{-\infty}^{\,\zeta} A(\xi)\varphi(\xi)\,d\xi\Big)d\zeta
\vspace{5pt}\\
&\di =-\varphi(z)\int_0^z\Big(\varphi^{-2}(\zeta)
\int_\zeta^\infty A(\xi)\varphi(\xi)\,d\xi\Big)d\zeta,
\end{array}
\end{equation}
where $\varphi:={U_0}'$.  The estimate \eqref{psi-bound} now
follows from the above expression and \eqref{U0-asympt}. \qed
\\

{}From the above lemma, the solvability condition for \eqref{eqU1}
is given by
$$
\int_\R \Big[ {{U_0}'}^2(z)(\widetilde d_{t}-\Delta \widetilde
d)(x,t) + g(x,t,U_0(z)){U_0}'(z) \Big ]dz= 0,
$$
for all $(x,t) \in Q_T$. Hence we get
\begin{equation*}
\widetilde d_{t} - \Delta \widetilde d = -\frac{\int_\R
g(x,t,U_0(z)){U_0}'(z)\,dz}{\int_\R {{U_0}'}^2(z) \,dz},
\end{equation*}
which gives
\begin{equation*}
\widetilde d_{t}= \Delta \widetilde d -
\frac{G(x,t,\ap)-G(x,t,\am)}{\int_\R {{U_0}'}^2(z) \,dz}.
\end{equation*}
Moreover, multiplying equation \eqref{eq-phi} by ${U_0}'$ and
integrating it from $-\infty$ to $z$, we obtain
$$
\begin{array}{ll}
0&=\di{\int_{-\infty}^z}\big({U_0} ''{U_0}'+f(U_0){U_0}'\big)(s)ds \vsp \\
&=\di{\frac 12} {{U_0} '}^2(z)-W(U_0(z))+W(\am),
\end{array}
$$
where we have also used the fact that $U_0(-\infty)=\am$ and
${U_0} '(-\infty)=0.$ This implies that
$$
{U_0}'(z)=\sqrt 2 \big(W(U_0(z))-W(\am)\big)^{ 1/2},
$$
and therefore
$$
\begin{array}{ll}
\di{\int_\R} {{U_0} '}^2(z)dz &=\di{\int_\R} {U_0}'(z)\sqrt 2
\big(W(U_0(z))-W(\am)\big)^{ 1/2}dz\vsp \\
&=\sqrt 2\di{\int_{\am}^{\ap}} (W(s) - W(\am))^{1/2} ds.
\end{array}
$$
It then follows, in view of the definition of $c_0$ in \eqref{c0},
that
\begin{equation}\label{eq-d}
\widetilde d_{t} = \Delta \widetilde d  -
c_0(G(x,t,\ap)-G(x,t,\am)).
\end{equation}
We are now ready to derive the equation of interface motion. Since
$\n \widetilde d$ $(=\nabla_x\, \widetilde d(x,t))$ coincides with
the outward normal unit vector to the hypersurface $\Gamma _t$, we
have $\widetilde d_{t}(x,t)=-V_n$, where $V_n$ is the normal
velocity of the interface $\Gamma _t$. It is also known that the
mean curvature $\kappa$ of the interface is equal to $\Delta
\widetilde d/(N-1)$. Thus the equation of interface motion is
given by:
\begin{equation}\label{eq:2-6}
V_n=-(N-1)\kappa + c_0(G(x,t,\ap)-G(x,t,\am)) \quad \mbox{on}
\enskip \Gamma _t.
\end{equation}
Summarizing, under the assumption that the solution $\ue$ of
problem $\Pe$ satisfies
\begin{equation*}
\ue\to
\begin{cases}
\ap &\quad  \textrm{ in } Q_T^+ \\
\am &\quad  \textrm{ in } Q_T^-,
\end{cases}\qquad\hbox{as}\ \ \ep\to 0,
\end{equation*}
we have formally proved that the boundary $\Gamma _t$ between
$\Omega_t^-$ and $\Omega_t^+$ moves according to the law
\eqref{eq:2-6}.\\

To conclude this section, we give basic estimates for
$U_1(x,t,z)$, which we will need in Section \ref{s:motion} to
study the motion of interface.  Substituting \eqref{eq-d} into
\eqref{eqU1} gives
\begin{equation}\label{eqU1-a}
\left\{\begin{array}{ll}
U_{1zz}+f'(U_0(z))U_1=g(x,t,U_0(z))- \cha (x,t){U_0}'(z),\vsp \\
U_1(x,t,0)=0, \quad \quad \quad \quad \quad \quad \quad
U_1(x,t,\cdot) \in L^\infty(\R),
\end{array}\right.
\end{equation}
where
\begin{equation}\label{gam}
\cha (x,t)= c_0 (G(x,t,\ap)-G(x,t,\am)).
\end{equation}
Thus $U_1(x,t,z)$ is a solution of \eqref{eq-psi} with
\begin{equation}\label{A0}
A=A_0(x,t,z):=g(x,t,U_0(z))- \cha (x,t){U_0}'(z),
\end{equation}
where the variables $x,t$ are considered parameters.  The problem
\eqref{eqU1-a} has a unique solution by virtue of Lemma
\ref{Fredholm}. Moreover, since $A_0(x,t,z)$ remains bounded as
$(x,t,z)$ varies in $\ombar\times[0,T]\times\R$, the estimate
\eqref{psi-bound} implies
\begin{equation}\label{def-M}
|U_1(x,t,z)|\leq M\qquad\hbox{for}\ \ x\in\ombar,\;
t\in[0,T],\;z\in\R,
\end{equation}
for some constant $M>0$.  Similarly, since $\nabla U_1$ is a
solution of \eqref{eq-psi} with
\[
A= \nabla_x A_0(x,t,z) \;\Big(= \nabla_x \big(g(x,t,U_0(z))-
\cha(x,t){U_0}'(z)\big)\;\Big),
\]
and since $g$ is assumed to be $C^1$ in $x$, we obtain
\begin{equation}\label{def-M2}
|\nabla_x U_1(x,t,z)|\leq M\qquad\hbox{for}\ \ x\in\ombar,\;
t\in[0,T],\;z\in\R,
\end{equation}
for some constant $M>0$.

To obtain estimates as $z\to\pm\infty$, we first observe that
\eqref{U0-asympt} implies
\begin{equation}\label{A0-g}
A_0(x,t,z)-g(x,t,\apm)=O(e^{-\lambda |z|})\qquad\hbox{as}\ \ z\to
\pm \infty,
\end{equation}
uniformly in $x\in\ombar,\,t\in[0,T]$.  We then apply the
following general estimates:

\begin{lem}\label{psi-decay}
Let the assumptions of Lemma \ref{Fredholm} hold, and assume
further that $A(z)-A^{\pm}=O(e^{-\delta |z|})$ as $z\rightarrow
\pm \infty$ for some constants $A^+,\,A^-$ and $\delta >0$.  Then
there exists a constant $\lambda>0$ such that
\begin{equation}\label{U1-exp-bound}
\psi(z)-\frac{A^{\pm}}{f'(\apm)}=O(e^{- \lambda |z|}),\quad\
|\psi'(z)|+|\psi''(z)|=O(e^{- \lambda |z|}),
\end{equation}
as $z\to\pm\infty$.
\end{lem}

{\noindent \bf Proof.} We only state the outline. To derive the
former estimate, we need a slightly more elaborate version of
\eqref{U0-asympt}.  Since $f(u)$ is $C^2$, we have
$f(u)=(u-\apm)f'(\apm)+O\big((u-\apm)^2\big)$.  Consequently,
\begin{equation}\label{U0-asympt-a}
U_0(z)=\left\{
\begin{array}{ll}
\am+C_1\,e^{\,\lambda_- z}+O(e^{\,2\lambda_- z})
\qquad&\hbox{as} \ z\to -\infty,\\
\ap+C_2\,e^{-\lambda_+ z}+O(e^{-2\lambda_+ z}) \qquad&\hbox{as} \
z\to +\infty.
\end{array}\right.
\end{equation}
Using the expression \eqref{psi} along with the estimate
$A(z)-A^{\pm}=O(e^{-\delta |z|})$ and \eqref{U0-asympt-a}, we see
that
\[
\psi(z)=-\frac{A^\pm}{(\lambda_\pm)^2}
+O\big(\,|z|e^{-\lambda_\pm|z|}\,\big)+
O\big(\,e^{-\min(\delta,\lambda_\pm)|z|}\,\big) \quad\ \hbox{as} \
z\to \pm\infty.
\]
This implies the former estimate in \eqref{U1-exp-bound}, where
$\lambda$ can be any constant satisfying
$0<\lambda<\min(\lambda_-,\lambda_+,\delta)$.  Substituting this
into equation \eqref{eq-psi} gives the estimate for $\psi_{zz}$.
Finally, the estimate for $\psi_{z}$ follows by integrating
$\psi_{zz}$ from $\pm\infty$ to $z$.
\qed\\

{}From the above lemma and \eqref{A0-g} we obtain the estimate
\begin{equation}\label{est-psi}
|U_{1z}(x,t,z)|+|U_{1zz}(x,t,z)| \leq Ce^{-\lambda |z|},
\end{equation}
for $x\in\ombar,\,t\in [0,T],\,z\in\R$. Similarly, since
\eqref{diff-U0} implies
\[
(\n _xA_0)(x,t,z)-(\n _x g)(x,t,\apm)=O(e^{-\lambda
|z|})\qquad\hbox{as}\ \ z\to \pm \infty,
\]
we can apply Lemma \ref{psi-decay} to $\psi=\nabla_x U_1$, to
obtain
\begin{equation*}
|\n _x U_{1z}(x,t,z)|+|\n _x U_{1zz}(x,t,z)| \leq Ce^{-\lambda|z|}
\end{equation*}
for $x\in\ombar,\,t\in [0,T],\,z\in\R$.  As a consequence, there
is a constant, which we denote again by $M$, such that
\begin{equation}\label{def-M3}
|\n _x U_{1z}(x,t,z)|\leq M.
\end{equation}

Next we consider the boundary condition.  Note that
\eqref{g-neumann} implies
\begin{equation}\label{droite-neumann}
\frac{\partial}{\partial \nu} A_0= \frac{\partial}{\partial
\nu}\Big[g(x,t,U_0(z))-\cha(x,t){U_0}'(z)\Big]=0 \quad\text{ on }
\partial \Omega.
\end{equation}
Consequently, from the expression \eqref{psi}, or equivalently the
expression
\[
U_1(x,t,z) =U_0'(z)\int_0^z\Big(\big(U_0'(\zeta)\big)^{-2}
\int_{-\infty}^{\,\zeta} A_0(x,t,\xi)U_0'(\xi)\,d\xi\Big)d\zeta,
\]
we see that
\begin{equation}\label{U1-neumann}
\frac{\partial U_1}{\partial \nu} =0 \quad\text{ on }\dom.
\end{equation}

\section{Generation of interface: the case $\g\equiv 0$}
\label{s:generation-0}

This section deals with the generation of interface, namely the
rapid formation of internal layers that takes place in a
neighborhood of $\Gamma_0=\{x\in \om,\,u_0(x)=a\}$ within the time
span of order $\ep^2 |\ln\ep|$.  For the time being we focus on
the special case where $\g\equiv 0$.  We will discuss the general
case in Section \ref{s:generation-g}. In the sequel, $\eta _0$
will stand for the following quantity:
\[
\eta _0:= \min (a-\am,\ap -a).
\]
Our main result in this section is the following:

\begin{thm}\label{g-th-gen}
Let $\eta \in (0,\eta _0)$ be arbitrary and define $\mu$ as the
derivative of $f(u)$ at the unstable zero $u=a$, that is
\begin{equation}\label{g-def-mu}
\mu=f'(a).
\end{equation}
Then there exist positive constants $\ep_0$ and $M_0$ such that,
for all $\,\ep \in (0,\ep _0)$,
\begin{enumerate}
\item for all $x\in\om$,
\begin{equation}\label{g-part1}
\am-\eta \leq u^\ep(x,\mu ^{-1} \ep ^2 | \ln \ep |) \leq \ap+\eta,
\end{equation}
\item for all $x\in\om$ such that $|u_0(x)-a|\geq M_0 \ep$, we
have that
\begin{align}
&\text{if}\;~~u_0(x)\geq a+M_0\ep\;~~\text{then}\;~~u^\ep(x,\mu
^{-1} \ep ^2 | \ln \ep |)
\geq \ap-\eta,\label{g-part2}\\
&\text{if}\;~~u_0(x)\leq a-M_0\ep\;~~\text{then}\;~~u^\ep(x,\mu
^{-1} \ep ^2 | \ln \ep |)\leq \am+\eta \label{g-part3}.
\end{align}
\end{enumerate}
\end{thm}

The above theorem will be proved by constructing a suitable pair
of sub- and super-solutions. Note that we do not need condition
\eqref{int-f} in proving this theorem.

\subsection{The bistable ordinary differential equation}
\label{SS:bistable ODE}

Let us first consider the problem without diffusion:
\begin{equation*}\label{no-diffusion}
\bar{u}_t=\frac{1}{\ep^2}\,f(\bar{u}), \qquad \bar{u}(x,0)=u_0(x).
\end{equation*}
This solution is written in the form
\[
\bar{u}(x,t)=Y\Big(\frac{t}{\ep^2},\,u_0(x)\Big),
\]
where $Y(\tau,\xi)$ denotes the solution of the ordinary
differential equation
\begin{equation}\label{g-ode}
\left\{\begin{array}{ll} Y_\tau (\tau,\xi)&=f(Y(\tau,\xi)) \quad
\text { for } \tau >0 \vspace{3pt}\\
Y(0,\xi)&=\xi.
\end{array}\right.
\end{equation}
Here $\xi$ ranges over the interval $(-2C_0,2C_0)$, with $C_0$
being the constant defined in \eqref{int1}. We first study basic
properties of $Y$.

\begin{lem}\label{g-Y1}
We have $Y_\xi >0$, for all $\xi \notin \{\am,a,\ap\}$, $\tau >
0$. Furthermore,
$$
Y_\xi(\tau,\xi)=\frac{f (Y(\tau,\xi))}{f (\xi)}.
$$
\end{lem}

{\noindent \bf Proof.} First, differentiating equation
\eqref{g-ode} by $\xi$, we obtain
\begin{equation*}
\left\{\begin{array}{ll} Y_{\xi\tau}=Y_\xi f'(Y)\vspace{3pt}\\
Y_\xi(0,\xi)=1,
\end{array}\right.
\end{equation*}
which is integrated as follows:
\begin{equation}\label{g-Y2}
Y_\xi(\tau,\xi)=\exp \Big[\int_0^\tau f'(Y(s,\xi))ds\Big]>0.
\end{equation}
We then differentiate equation \eqref{g-ode} by $\tau$ and obtain
\begin{equation*}
\left\{\begin{array}{ll} Y_{\tau\tau}=Y_\tau f'(Y)\vspace{3pt}\\
Y_\tau(0,\xi)=f (\xi),
\end{array}\right.
\end{equation*}
which in turn implies
\[
\begin{array}{ll}Y_\tau(\tau,\xi)&=f (\xi) \exp
\Big[\di{\int_0^\tau} f'(Y(s,\xi))ds\Big]\vsp \\ &=f
(\xi)Y_\xi(\tau,\xi).\end{array}
\]
This last equality, in view of \eqref{g-ode}, completes the proof
of Lemma \ref{g-Y1}. \qed

\vskip 8pt For $\xi \notin \{\am,a,\ap\}$, we define a function
$A(\tau,\xi)$ by
\begin{equation}\label{g-A}
A(\tau,\xi)=\frac{f'(Y(\tau,\xi))-f'(\xi)}{f (\xi)}.
\end{equation}

\begin{lem}\label{g-Y5}
We have, for all $\xi \notin \{\am,a,\ap\}$, $\tau > 0$,
$$
A(\tau,\xi)=\int_0^\tau f''(Y(s,\xi))Y_\xi(s,\xi)ds.
$$
\end{lem}

{\noindent \bf Proof.} Differentiating by $\xi$ the equality of
Lemma \ref{g-Y1} leads to
\begin{equation}\label{g-Y4}
Y_{\xi\xi}=A(\tau,\xi)Y_\xi,
\end{equation}
whereas differentiating \eqref{g-Y2} by $\xi$ yields
$$
Y_{\xi\xi}=Y_\xi \int_0^\tau f''(Y(s,\xi))Y_\xi(s,\xi)ds.
$$
These two last results complete the proof of Lemma \ref{g-Y5}.
\qed

\vskip 8pt Next we need some estimates on $Y$ and its derivatives.
First, we estimate the speed of the evolution of $Y$ when the
initial value $\xi$ lies between $\am+\eta$ and $\ap-\eta$.

\vskip 8pt
\begin{lem}\label{g-est-derY-A-milieu}
Let $\eta \in (0,\eta _0)$ be arbitrary. Then there exist positive
constants $\tilde C _1=\tilde C _1(\eta)$, $\tilde C _2=\tilde C
_2(\eta)$ and $C_3=C_3(\eta)$ such that
\begin{enumerate}
\item if $\xi \in (a,\ap-\eta)$ then, for every $\tau>0$ such that
$Y(\tau,\xi)$ remains in the interval $(a,\ap-\eta)$, we have
\begin{equation}\label{g-est-Y3}
\tilde C _1e^{\mu\tau}\leq Y_\xi(\tau,\xi) \leq \tilde C _2
e^{\mu\tau},
\end{equation}
\begin{equation}\label{g-est-A-milieu}
|A(\tau,\xi)|\leq C_3(e^{\mu \tau}-1),
\end{equation}
where $\mu$ is the constant defined in \eqref{g-def-mu}; \item if
$\xi\in (\am+\eta,a)$ then, for every $\tau>0$ such that
$Y(\tau,\xi)$ remains in the interval $(\am+\eta,a)$, we have
\eqref{g-est-Y3} and \eqref{g-est-A-milieu}.
\end{enumerate}
\end{lem}

{\noindent \bf Proof.} We take $\xi \in (a,\ap-\eta)$ and suppose
that, for $s \in (0,\tau)$, $Y(s,\xi)$ remains in the interval
$(a,\ap-\eta)$. Integrating the equality
$$
\frac{Y_\tau(s,\xi)}{f (Y(s,\xi))}=1
$$
from $0$ to $\tau$ yields
\begin{equation}\label{g-tt}
\int_0^{\tau} \frac{Y_\tau (s,\xi)}{f (Y(s,\xi))}ds \quad =\tau.
\end{equation}
Hence by the change of variable $q=Y(s,\xi)$ we get
\begin{equation}\label{g-tau}
\int _\xi ^{Y(\tau,\xi)} \frac{dq}{f (q)}=\tau.
\end{equation}
Moreover, the equality of Lemma \ref{g-Y1} leads to
\begin{equation}\label{naka}
\begin{array}{lll}
\ln Y_\xi (\tau,\xi)&=
\di{\int _ \xi ^{Y(\tau,\xi)}} \frac{f'(q)}{f(q)}dq \vsp \\
&=\di{\int _ \xi ^{Y(\tau,\xi)}}\big
[\frac{f'(a)}{f(q)}+\frac{f'(q)-f'(a)}{f(q)}\big ]dq \vsp \\
&=\mu \tau+\di{\int _ \xi ^{Y(\tau,\xi)}}h(q)dq,
\end{array}
\end{equation}
where $h(q)=(f'(q)-\mu)/f(q)$. As $h(q)$ tends to $f''(a)/f'(a)$
when $q$ tends to $a$, $h$ is continuous on $[a,\ap-\eta]$. Hence
we can define
\[
H=H(\eta):=\Vert h \Vert _{L^\infty (a,\ap-\eta)}.
\]
Since $|Y(\tau,\xi)-\xi|$ takes its value in the interval $[0,\ap
-a-\eta]\subset[0,\ap-a]$, it follows from \eqref{naka} that
\[
\mu \tau -H(\ap-a) \leq \ln Y_\xi(\tau,\xi) \leq \mu
\tau+H(\ap-a),
\]
which, in turn, proves \eqref{g-est-Y3}. Next Lemma \ref{g-Y5} and
\eqref{g-est-Y3} yield
\[
\begin{array}{ll}|A(\tau,\xi)| &\leq \Vert f'' \Vert
_{L^\infty(\am,\ap)} \di{\int_0^\tau} \tilde C _2 e^{\mu s}ds \vsp \\
&\leq C_3(e^{\mu\tau}-1),
\end{array}
\]
which completes the proof of \eqref{g-est-A-milieu}. The case
where $\xi$ and $Y(\tau,\xi)$ are in $(\am+\eta,a)$ is similar and
omitted. \qed

\vskip 8pt
\begin{cor}\label{g-est-Y-milieu} Let $\eta \in (0,\eta
_0)$ be arbitrary. Then there exist positive constants
$C_1=C_1(\eta)$ and $C_2=C_2(\eta)$ such that
\begin{enumerate}
\item if $\xi\in (a,\ap-\eta)$ then, for every $\tau>0$ such that
$Y(\tau,\xi)$ remains in the interval $(a,\ap-\eta)$, we have
\begin{equation}\label{g-est-Y-1}
C_1e^{\mu \tau}(\xi-a)\leq Y(\tau,\xi)-a \leq C_2e^{\mu
\tau}(\xi-a);
\end{equation}
\item if $\xi \in (\am +\eta,a)$ then, for every $\tau>0$ such
that $Y(\tau,\xi)$ remains in the interval $(\am +\eta,a)$, we
have
\begin{equation}\label{g-est-Y-2}
C_2e^{\mu \tau}(\xi-a)\leq Y(\tau,\xi)-a \leq C_1e^{\mu
\tau}(\xi-a).
\end{equation}
\end{enumerate}
\end{cor}

{\noindent \bf Proof.} We can find $B_1=B_1(\eta)>0$ and
$B_2=B_2(\eta)>0$ such that, for all $q\in (a,\ap-\eta)$,
\begin{equation}\label{g-ineg-fq}
B_1(q-a)\leq f (q) \leq B_2(q-a).
\end{equation}
We use this inequality for $a<Y(\tau,\xi) < \ap-\eta$ to obtain
$$
B_1(Y(\tau,\xi)-a)\leq f (Y(\tau,\xi))\leq B_2(Y(\tau,\xi)-a).
$$
We also use this inequality for $a<\xi < \ap-\eta$ to obtain
$$
B_1(\xi-a) \leq f(\xi) \leq B_2(\xi-a).
$$
Next we use the equality $Y_\xi=f (Y)/f (\xi)$ of Lemma \ref{g-Y1}
to deduce that
$$
\frac{B_1}{B_2}(Y(\tau,\xi)-a)\leq (\xi-a)Y_\xi(\tau,\xi) \leq
\frac{B_2}{B_1}(Y(\tau,\xi)-a),
$$
which, in view of \eqref{g-est-Y3}, implies that
$$
\frac{B_1}{B_2}\tilde C_1 e^{\mu \tau}(\xi-a) \leq Y(\tau,\xi)-a
\leq \frac{B_2}{B_1}\tilde C_2 e^{\mu \tau}(\xi-a).
$$
This proves \eqref{g-est-Y-1}. The proof of \eqref{g-est-Y-2} is
similar and is omitted. \qed

\vskip 8pt We now present estimates in the case where the initial
value $\xi$ is smaller than $\am+\eta$ or larger than $\ap-\eta$.

\begin{lem}\label{g-est-bords}
Let $\eta \in (0,\eta _0)$ and $M>0$ be arbitrary.  Then there
exists a positive constant $C_4=C_4(\eta,M)$ such that
\begin{enumerate}
\item if $\xi \in [\ap-\eta,\ap+M]$, then, for all $\tau> 0$,
$Y(\tau,\xi)$ remains in the interval $[\ap-\eta,\ap+M]$ and
\begin{equation}\label{g-est-A-bords}
|A(\tau,\xi)|\leq C_4\tau \quad\hbox{for}\ \ \tau> 0 \,;
\end{equation}
\item if $\xi \in [\am-M,\am+\eta]$, then, for all $\tau> 0$,
$Y(\tau,\xi)$ remains in the interval $[\am-M,\am+\eta]$ and
\eqref{g-est-A-bords} holds.
\end{enumerate}
\end{lem}

{\noindent \bf Proof.} Since statement (i) and statement (ii) can
be treated in the same way, we will only prove the former. The
fact that $Y(\tau,\xi)$ remains in the interval $[\ap-\eta,\ap+M]$
directly follows from the bistable properties of $f$, or, more
precisely, from the sign conditions $f(\ap-\eta)>0$, $f(\ap+M)<0$.

To prove \eqref{g-est-A-bords}, suppose first that $\xi\in
[\ap,\ap+M]$. In view of \eqref{der-f}, $f'$ is strictly negative
in an interval of the form $[\ap,\ap +c]$ and $f$ is negative in
$[\ap,\infty)$.  We denote by $-m<0$ the maximum of $f$ on
$[\ap+c,M]$.  Then, as long as $Y(\tau,\xi)$ remains in the
interval $[\ap+c,M]$, the ordinary differential equation
\eqref{g-ode} implies
$$
Y_\tau \leq -  m .
$$
This means that, for any $\xi\in [\ap,\ap+M]$, we have
\[
Y(\tau,\xi) \in [\ap,\ap +c]\qquad\hbox{for}\ \ \tau
\geq\overline{\tau}:=\frac{M-c}{m}.
\]
In view of this, and considering that $f'(Y)<0$ for $Y\in
[\ap,\ap+c]$, we see from the expression \eqref{g-Y2} that
\[
\begin{array}{ll}
Y_\xi(\tau,\xi) &\di\vspace{6pt}
=\exp\Big[\int_0^{\overline{\tau}}f'(Y(s,\xi))ds\Big]\,
\exp\Big[\int_{\overline{\tau}}^\tau f'(Y(s,\xi))ds\Big]\\
&\di\vspace{6pt}
\leq \exp \Big[\int _0^{\overline{\tau}}f'(Y(s,\xi))ds\Big]\\
&\di \leq \exp\Big[\int _0 ^{\overline{\tau}}
\sup_{z\in[\am-M,\ap+M]} |f'(z)|ds\Big] =:\tilde{C}_4,
\end{array}
\]
for all $\tau\geq\overline{\tau}$. It is clear from the same
estimate \eqref{g-Y2} that $Y_\xi\leq\tilde{C}_4$ holds also for
$0\leq\tau\leq\overline{\tau}$.  We can then use Lemma \ref{g-Y5}
to deduce that
$$
|A(\tau,\xi)|\leq \tilde{C}_4 \int_0^\tau|f''(Y(s,\xi))|ds \leq
C_4\tau.
$$
The case $\xi\in [\ap-\eta,\ap]$ can be treated in the same way.
This completes the proof of the lemma. \qed

\vskip 8pt Now we choose the constant $M$ in the above lemma
sufficiently large so that $[-2C_0,2C_0]\subset [\am-M,\ap+M]$,
and fix $M$ hereafter. Then $C_4$ only depends on $\eta$. Using
the fact that $\tau=O(e^{\mu\tau}-1)$ for $\tau>0$, one can easily
deduce from \eqref{g-est-A-milieu} and \eqref{g-est-A-bords} the
following general estimate.

\begin{lem}\label{g-EST-A}
Let $\eta \in (0,\eta _0)$ be arbitrary and let $C_0$ be the
constant defined in \eqref{int1}. Then there exists a positive
constant $C_5=C_5(\eta)$ such that, for all $\tau>0$ and all $\xi
\in(-2C_0,2C_0)$,
$$
|A(\tau,\xi)|\leq C_5(e^{\mu \tau}-1).
$$
\end{lem}

\subsection{Construction of sub- and super-solutions}

We are now ready to construct the sub- and super-solutions for the
study of generation of interface.  For simplicity, we first
consider the case where
\begin{equation}\label{int2}
\frac{\partial u_0}{\partial \nu} =0 \quad\text{ on }\dom.
\end{equation}
In this case, our sub- and super-solutions are given by
\begin{equation}\label{w+-}
w_\ep^\pm(x,t)=
Y\Big(\frac{t}{\ep^2},\,u_0(x)\pm\ep^2C_6(\emut-1)\Big).
\end{equation}
In the general case where \eqref{int2} does not necessarily hold,
we have to slightly modify $w_\ep^\pm(x,t)$ near the boundary
$\dom$. This will be discussed later.

\begin{lem}\label{g-w}
Assume \eqref{int2}.  Then there exist positive constants $\ep_0$
and $C_6$ such that, for all $\, \ep \in (0,\ep _0)$,
$(w_\ep^-,w_\ep^+)$ is a pair of sub- and super-solutions for
problem $\Pe$, in the domain $\ombar\times [0,\mu ^{-1} \ep^2|\ln
\ep|]$, satisfying $w^-_\ep(x,0)=w^+(x,0)=u_0(x)$. Consequently
\begin{equation}\label{g-coincee1}
w_\ep^-(x,t) \leq u^\ep(x,t) \leq w_\ep^+(x,t) \quad\ \hbox{for} \
x\in\ombar,\;0\leq t \leq \mu ^{-1} \ep^2|\ln \ep|.
\end{equation}
\end{lem}

{\noindent \bf Proof.}  The assumption \eqref{int2} implies
$$
\frac{\partial w_\ep^\pm}{\partial \nu}=0 \quad\text{ on }\partial
\Omega \times (0,+\infty).
$$
Now we define an operator ${\cal L}_0$ by
\[
{\cal L}_0 u :=u_t -\Delta u -\edeux f(u),
\]
and prove that ${\cal L}_0 w_\ep^+\geq 0$ . Straightforward
computations yield
$$
{\cal L}_0 w_\ep^+ =\edeux Y_\tau + C_6\,\mu\, \emut Y_\xi -\Delta
u_0 Y _ \xi -|\n u_0|^2 Y_{\xi \xi}-\edeux f(Y),
$$
therefore, in view of the ordinary differential equation
\eqref{g-ode},
$$
{\cal L}_0w_\ep^+=\Big[C_6 \,\mu\, \emut-\Delta u_0-
{\frac{Y_{\xi\xi}}{Y_\xi}}\,|\nabla u_0|^2\Big]Y_\xi.
$$
We note that, in the range $0 \leq t \leq \mu ^{-1} \ep ^2|\ln
\ep|$, we have, for $\ep _0$ sufficiently small,
$$
0\leq \ep^2 C_6(\emut-1) \leq \ep ^2C_6(\ep^{-1}-1) \leq C_0,
$$
where $C_0$ is the constant defined in \eqref{int1}. Hence
$$
\xi:=u_0(x)\pm C_6(\emut-1)\in (-2C_0,2C_0),
$$
and it follows from the estimate of $A=Y_{\xi\xi}/Y_\xi$ in Lemma
\ref{g-EST-A}, with the choice $\tau:=t/\ep ^2$, that
$$
\begin{array}{ll}
{\cal L}_0 w_\ep^+ &\geq \Big[C_6\,\mu \emut-|\Delta
u_0|-C_5(\emut-1)|\nabla u_0|^2\Big]Y_\xi \vsp \\
&\geq \Big[(C_6\,\mu-C_5 |\nabla u_0|^2)\emut-|\Delta u_0| +C_5
|\nabla u_0|^2\Big]Y_\xi.
\end{array}
$$
Since $Y_\xi
>0$, this inequality implies that, for $C_6$ large
enough,
$$
{\cal L}_0w^+_\ep \geq \Big[C_6 \mu-C_5 {C_0}^2 -C_0\Big]Y_\xi\geq
0.
$$
Hence $w_\ep^+$ is a super-solution for problem $\Pe$. Similarly
$w_\ep^-$ is a sub-solution. Obviously
$w^-_\ep(x,0)=w^+(x,0)=u_0(x)$.  Lemma \ref{g-w} is proved. \qed\\

In the more general case where \eqref{int2} is not necessarily
valid, one can proceed as follows: in view of \eqref{dalltint} and
\eqref{initial-data} there exist positive constants $d_1$, $\rho$
such that $u_0(x) \geq a+\rho$ if $d(x,\partial \Omega) \leq d_1$.
Let $\chi$ be a smooth cut-off function defined on $[0,+\infty)$
such that $0\leq \chi \leq 1$, $\chi(0)=\chi '(0)=0$ and $\chi
(z)=1$ for $z\geq d_1$. Then we define
\[
\begin{array}{l}\di\vspace{3pt}
u_0 ^+(x)=\chi(d(x,\partial \Omega))\,u_0(x)+\big[1-\chi
(d(x,\partial\Omega))\big]\max_{x\in\overline \Omega}\, u_0(x),
\\ \di
u_0 ^-(x)=\chi(d(x,\partial \Omega))\,u_0(x)+\big[1-
\chi(d(x,\partial \Omega))\big](a+\rho).
\end{array}
\]
Clearly, $u_0 ^- \leq u_0 \leq u_0^+$, and both $u_0^+$ and
$u_0^+$ satisfy \eqref{int2}. Now we set
\begin{equation*}\label{w+-2}
\tilde{w}_\ep^\pm(x,t)=
Y\Big(\frac{t}{\ep^2},\,u^\pm_0(x)\pm\ep^2C_6(\emut-1)\Big).
\end{equation*}
Then the same argument as in Lemma \ref{g-w} shows that
$(\tilde{w}_\ep^-,\tilde{w}_\ep^+)$ is a pair of sub- and
super-solutions for problem $\Pe$.  Furthermore, since
$\tilde{w}_\ep^-(x,0)=u_0^-(x)\leq u_0(x)\leq u_0^+(x)
=\tilde{w}_\ep^+(x,0)$, the comparison principle asserts that
\begin{equation}\label{g-coincee}
\tilde{w}_\ep^-(x,t) \leq u^\ep(x,t) \leq \tilde{w}_\ep^+(x,t)
\quad\ \hbox{for }x\in\ombar,\;0\leq t \leq \mu ^{-1} \ep^2|\ln
\ep|.
\end{equation}

\subsection{Proof of Theorem \ref{g-th-gen}}

In order to prove Theorem \ref{g-th-gen} we first present basic
estimates of the function $Y$ after a time of order $\tau\sim |\ln
\ep|.$

\begin{lem}\label{after-time}
Let $\eta \in (0,\eta _0)$ be arbitrary; there exist positive
constants $\ep_0$ and $C_7$ such that, for all $\,\ep\in(0,\ep
_0)$,
\begin{enumerate}
\item for all $\xi\in (-2C_0,2C_0)$,
\begin{equation}\label{g-part11}
\am-\eta \leq Y(\mu ^{-1} | \ln \ep |,\xi) \leq \ap+\eta;
\end{equation}
\item for all $\xi\in (-2C_0,2C_0)$ such that $|\xi-a|\geq C_7
\ep$, we have that
\begin{align}
&\text{if}\;~~\xi\geq a+C_7 \ep\;~~\text{then}\;~~Y(\mu ^{-1}| \ln
\ep |,\xi)
\geq \ap-\eta,\label{g-part22}\\
&\text{if}\;~~\xi\leq a-C_7 \ep\;~~\text{then}\;~~Y(\mu ^{-1}| \ln
\ep |,\xi)\leq \am+\eta \label{g-part33}.
\end{align}
\end{enumerate}
\end{lem}

{\noindent \bf Proof.} We first prove \eqref{g-part22}. For $\xi
\geq a+C_7\ep$, as long as $Y(\tau,\xi)$ has not reached
$\ap-\eta$, we can use \eqref{g-est-Y-1} to deduce that
$$
\begin{array}{lll}
Y(\tau,\xi) & \geq a +C_1 e^{\mu \tau}
(\xi-a)\vsp \\
& \geq a + C_1 C_7 e^{\mu \tau} \ep \vsp \\
& \geq \ap-\eta,
\end{array}
$$
provided that $\tau$ satisfies
$$
\tau \geq \tau ^\ep =: \mu ^{-1}\ln \frac {\ap -a-\eta}{C_1C_7
\ep}.
$$
Choosing
$$
C_7=\frac{\max (a-\am,\ap-a)-\eta}{C_1},
$$
we see that $\mu ^{-1}|\ln \ep| \geq \tau ^\ep$, which completes
the proof of \eqref{g-part22}. Using \eqref{g-est-Y-2}, one easily
proves \eqref{g-part33}.

Next we prove \eqref{g-part11}. First, in view of the profile of
$f$, if we leave from a $\xi \in [\am-\eta,\ap+\eta]$ then
$Y(\tau,\xi)$ will remain in $[\am-\eta,\ap+\eta]$. Now suppose
that $\ap+\eta \leq \xi \leq 2C_0$. We check below that $Y(\mu
^{-1}|\ln \ep|,\xi)\leq \ap+\eta$. First, in view of
\eqref{der-f}, we can find $p>0$ such that
\begin{equation}\label{g-pente}
\begin{array}{ll}\text { if } \quad \ap \leq
u \leq 2C_0 & \text { then } \quad  f(u)
\leq p(\ap-u)\vspace{3pt}\\
\text { if } \quad -2C_0 \leq u \leq \am & \text { then } \quad
f(u) \geq -p(u-\am).
\end{array}
\end{equation}
We then use the ordinary differential equation to obtain, as long
as $\ap+\eta \leq Y \leq 2C_0$, the inequality $Y_\tau \leq
p(\ap-Y)$. It follows that
$$
\frac{Y_\tau}{Y-\ap}\leq -p.
$$
Integrating this inequality from $0$ to $\tau$ leads to
$$
\begin{array}{ll}Y(\tau,\xi) &\leq \ap+(\xi-\ap)e^{-p\tau}\vspace{3pt}\\
&\leq \ap+(2C_0-\ap)e^{-p\tau}.
\end{array}
$$
One easily checks that, for $\ep \in (0,\ep _0)$, with $\ep_0=\ep
_0(\eta)$ small enough, we have $Y(\mu ^{-1}| \ln \ep|,\xi)\leq
\ap+\eta$, which completes the proof of \eqref{g-part11}. \qed

\vskip 8pt We are now ready to prove Theorem \ref{g-th-gen}. By
setting $t=\mu ^{-1} \ep ^2|\ln \ep|$ in \eqref{g-coincee}, we
obtain
\begin{multline}\label{g-gr}
Y\big(\mu ^{-1}|\ln \ep|, u_0^-(x)-(C_6 \ep -C_6 \ep ^2 )\big)\\
\leq u^\ep(x,\mu ^{-1} \ep^2|\ln \ep|) \leq Y\big(\mu ^{-1}|\ln
\ep|, u_0^+(x)+C_6 \ep -C_6 \ep ^2\big).
\end{multline}
Furthermore, by the definition of $C_0$ in \eqref{int1}, we have,
for $\ep_0$ small enough,
$$
- 2 C_0 \leq u_0^\pm(x)\pm (C_6 \ep -C_6 \ep ^2) \leq 2C_0
\qquad\hbox{for}\ x\in\om.
$$
Thus the assertion \eqref{g-part1} of Theorem \ref{g-th-gen} is a
direct consequence of \eqref{g-part11} and \eqref{g-gr}.

Next we prove \eqref{g-part2}. We choose $M_0$ large enough so
that $M_0\ep- C_6 \ep+C_6 \ep ^2 \geq C_7\ep$.  Then, for any
$x\in \om$ such that $u_0^-(x)\geq a+M_0 \ep$, we have
$$
u_0^-(x)-(C_6 \ep -C_6 \ep ^2)\geq a+M_0\ep- C_6 \ep+C_6 \ep ^2
\geq a+C_7\ep.
$$
Combining this, \eqref{g-gr} and \eqref{g-part22}, we see that
\[
u^\ep(x,\mu^{-1} \ep ^2 | \ln \ep |)\geq \ap-\eta,
\]
for any $x\in\om$ with $u_0^-(x)\geq a+M_0 \ep$.  From the
definition of $u_0^-$ it is clear that $u_0^-(x)\geq a+M_0 \ep$ if
and only if $u_0(x)\geq a+M_0 \ep$, provided that $\ep$ is small
enough. This proves \eqref{g-part2}.  The inequality
\eqref{g-part3} can be shown the same way.  This completes the
proof of Theorem \ref{g-th-gen}.\qed

\subsection{Optimality of the generation time}

To conclude this section we show that the generation time
${t}^{\,\ep}:=\mu^{-1}\ep^2|\ln\ep|$ that appears in Theorem
\ref{g-th-gen} is optimal.  In other words, the interface will not
be fully developed until $t$ comes close to $t^{\,\ep}$.

\begin{prop}\label{pr:optimal-time}
Denote by $t^{\,\ep}_{min}$ the smallest time such that
\eqref{resultat} holds for all $t\in [\,t^{\,\ep}_{min},\,T\,].$
Then there exists a constant $b=b(C)$ such that
\[
{t}^{\,\ep}_{min}\geq \mu^{-1} \ep^2 (|\ln\ep| - b),
\]
for all $\,\ep \in (0,\ep_0).$
\end{prop}

{\noindent \bf Proof.} For simplicity, we deal with the case where
\eqref{int2} is valid. In that case, \eqref{g-coincee1} holds for
all small $\ep>0$. For each $b>0$, we put
\[
t ^\ep (b):=\mu ^{-1} \ep ^2 (|\ln \ep|-b),
\]
and evaluate $\ue(x,t^\ep(b))$ at a point $x \in \om _0 ^+$ where
$dist(x,\Gamma _0)=C\ep$.  Since $u_0=a$ on $\Gamma_t$ and since
$|\nabla u_0|\leq C_0$ by \eqref{int1}, we have
\begin{equation}\label{m-M}
u_0(x) \leq a +C_0 C \ep.
\end{equation}
It follows from this and \eqref{g-est-Y-1} that
\[
\begin{array}{lll}
w^+ _\ep (x, t^\ep(b))& =Y\Big ( \mu ^{-1} (|\ln \ep|-b),
u_0(x)+\ep C_6 e^{-b}-\ep ^2 C_6 \Big ) \vsp \\
& \leq a+C_2e^{|\ln \ep|-b} \big(u_0(x)+\ep C_6 e^{-b}-\ep^2
C_6-a ) \vsp \\
& \leq a+C_2 \ep^{-1} e^{-b} (C_0C \ep+\ep C_6 e^{-b}) \vsp \\
& = a+C_2 e^{-b} (C_0C + C_6 e^{-b}).
\end{array}
\]
Now we choose $b$ to be sufficiently large, so that
\[
a+C_2 e^{-b} (C_0C + C_6 e^{-b}) <\ap-\eta.
\]
Then the above estimate and \eqref{g-coincee1} yield
\[
\ue(x,t^\ep(b))\leq w^+ _\ep(x, t^\ep(b))<\ap-\eta.
\]
This implies that \eqref{resultat} does not hold at $t=t^\ep(b)$,
hence $t^\ep(b)<t^\ep_{min}$.  The lemma is proved. \qed

\section{Generation of interface in the general case}
\label{s:generation-g}

In this section we extend Theorem \ref{g-th-gen} to the case where
$\g\not\equiv 0$. The proof is more technical than the case
$\g\equiv 0$, but the underlying ideas are the same. Hence we will
basically follow the argument of Section \ref{s:generation-0},
simply pointing out the main differences.

\subsection{The perturbed ordinary differential equation}

We first consider a slightly perturbed nonlinearity:
$$
\f(u)=f(u)+\delta ,
$$
where $\de$ is any constant. For $|\delta|$ small enough, this
function is still bistable. More precisely, $\f$ has the following
properties, whose proof is omitted:

\begin{lem}
Let $\de_0$ be small enough.  Then for any $\de\in(-\de_0,\de_0)$,
\begin{enumerate}
\item $\f$ has exactly three zeros, namely $\hm < \h < \hp$, and
there exists a positive constant $C$ such that
\begin{equation}\label{h}
|\hm -\am|+|\h -a|+|\hp-\ap|\leq C|\de|.
\end{equation}
\item We have
\begin{equation}\label{signeF}
\begin{array}{ll}
\f>0 \quad \text {in}\quad  (-\infty,\hm)\cup(\h,\hp),
\vspace{3pt}\\
\f<0 \quad \text{in}\quad(\hm,\h)\cup(\hp,+\infty).
\end{array}
\end{equation}
\item There exists a positive constant, denoted again by $C$, such
that
\begin{equation}\label{mu}
| \mm -\mu| \leq C|\de|,
\end{equation}
where
\[
\mm:= \f '(\h)=f'(\h).
\]
\end{enumerate}
\end{lem}

Now for each $\de\in(-\de_0,\de_0)$, we define $Y(\tau,\xi;\de)$
as the solution of the following ordinary differential equation:
\begin{equation}\label{ode}
\left\{\begin{array}{ll} Y_\tau (\tau,\xi;\de)&=\f
(Y(\tau,\xi;\de)) \quad \text { for }\tau >0\vspace{3pt}\\
Y(0,\xi;\de)&=\xi ,
\end{array}\right.
\end{equation}
where $\xi$ varies in $(-2C_0,2C_0)$, with $C_0$ being the
constant defined in \eqref{int1}.

To prove Theorem \ref{g-th-gen}, we will construct a pair of sub-
and super-solutions for $\Pe$ by simply replacing the function
$Y(\tau,\xi)$ in \eqref{w+-} by $Y(\tau,\xi;\de)$, with an
appropriate choice of $\de$.  For this strategy to work, we have
to check that the basic properties of $Y(\tau,\xi)$ in Subsection
\ref{SS:bistable ODE} carry over to $Y(\tau,\xi;\de)$.

First, it is clear that all the differential and integral
identities in Subsection \ref{SS:bistable ODE} that follow
directly from \eqref{g-ode} are still valid for \eqref{ode}. In
particular, Lemmas \ref{g-Y1} and \ref{g-Y5} remain to hold if we
replace $Y(\tau,\xi)$ by $Y(\tau,\xi;\de)$, $f$ by $\f$ and
$A(\tau,\xi)$ by $A(\tau,\xi;\de)$, where
\begin{equation*}\label{g-Ad}
A(\tau,\xi,\de)=\frac{\f'(Y(\tau,\xi;\de))-\f'(\xi)}{\f(\xi)}.
\end{equation*}
Next let us show that the basic estimates which we have
established in Subsection \ref{SS:bistable ODE} are also valid for
$Y(\tau,\xi;\de)$. The following lemma, which is an analogue of
Lemma \ref{g-est-derY-A-milieu}, is fundamental.

\begin{lem}\label{est-derY-A-milieu}
Let $\eta \in (0,\eta _0)$ be arbitrary. Then there exist positive
constants $\de_0=\de_0(\eta)$, $\widetilde C _1=\widetilde C
_1(\eta)$, $\tilde C _2=\tilde C _2(\eta)$ and $C_3=C_3(\eta)$
such that, for any $\de \in [-\de _0,\de_0]$,
\begin{enumerate}
\item if $\xi \in (\h,\,\ap-\eta)$ then, for every $\tau>0$ such
that $Y(\tau,\xi;\de)$ remains in the interval $(\h,\,\ap-\eta)$,
we have
\begin{equation}\label{est-Y3}
\tilde C _1e^{\m\tau}\leq Y_\xi(\tau,\xi;\de) \leq \tilde C _2
e^{\m\tau},
\end{equation}
\begin{equation}\label{est-A-milieu}
|A(\tau,\xi;\de)|\leq C_3(e^{\m\tau}-1);
\end{equation}
\item the same estimates as above hold if the interval
$(\h,\,\ap-\eta)$ is replaced by $(\am+\eta,\,\h)$ .
\end{enumerate}
\end{lem}

{\noindent \bf Proof.} In view of \eqref{h}, we can choose a small
constant $\de_0=\de_0(\eta)>0$ such that $(\h,\ap-\eta)\subset
(\h,\ap (\de))$, for every $\de\in [-\de _0,\de _0]$.  Therefore
$\f(q)$ does not change sign in the interval $(\h,\ap-\eta)$.
Thus, in order to prove the lemma, we just have to write again the
proof of Lemma \ref{g-est-derY-A-milieu}, simply replacing
$Y(\tau,\xi)$ by $Y(\tau,\xi;\de)$. We do not repeat the entire
proof here. Instead, let us explain why $\tilde C_1$, $\tilde C_2$
and $C_3$ can be chosen independent of $\de$.  In view of the
proof of Lemma \ref{g-est-derY-A-milieu}, it is sufficient to
estimate, for $q\in(\;\!a(\de),\,\ap-\eta]$, the modulus of the
quantity
\[
h_\de (q):=\frac{f'(q)-f'(a(\de))}{f_\de(q)}
\]
by a constant depending on $\eta$, but not on
$\de\in[-\de_0,\de_0]$. Since
\[
h_\de(q)\to \frac{\f''(a(\de))}{\f'(a(\de))}
=\frac{f''(a(\de))}{f'(a(\de))} \quad\ \hbox{as} \ \ q \to a(\de),
\]
we see that the function $(q,\de)\mapsto h_\de(q)$ is continuous
in the compact region $\{\,|\de|\leq \de_0,\;a(\de)\leq q\leq
\ap-\eta\,\}$. It follows that $|h_\de(q)|$ is bounded as
$(q,\de)$ varies in this region. This completes the proof of Lemma
\ref{est-derY-A-milieu}. \qed

\begin{cor}\label{est-Y-milieu}
Let $\eta \in (0,\eta _0)$ be arbitrary. Then there exist positive
constants $\de_0=\de_0(\eta)$, $C_1=C_1(\eta)$ and $C_2=C_2(\eta)$
such that, for any $\de \in [-\de _0,\de _0]$,
\begin{enumerate}
\item if $\xi\in (\h,\ap-\eta)$ then, for every $\tau>0$ such that
$Y(\tau,\xi;\de)$ remains in the interval $(\h,\ap-\eta)$, we have
\begin{equation}\label{est-Y-1}
C_1e^{\mm \tau}(\xi-\h)\leq Y(\tau,\xi;\de)-\h \leq C_2e^{\mm
\tau}(\xi-\h);
\end{equation}
\item if $\xi \in (\am +\eta,\h)$ then, for every $\tau>0$ such
that $Y(\tau,\xi;\de)$ remains in the interval $(\am +\eta,\h)$,
we have
\begin{equation}\label{est-Y-2}
C_2e^{\mm \tau}(\xi-\h)\leq Y(\tau,\xi;\de)-\h \leq C_1e^{\mm
\tau}(\xi-\h).
\end{equation}
\end{enumerate}
\end{cor}

{\noindent \bf Proof.} We can simply follow the proof of Corollary
\ref{g-est-Y-milieu}. In order to prove that $C_1$ and $C_2$ are
independent of $\de$, all we have to do is to find constants
$B_1=B_1(\eta)>0$ and $B_2=B_2(\eta)>0$ such that, for all $\de\in
[-\de_0,\de_0]$ and all $q\in (a(\de),\ap-\eta)$,
\begin{equation}\label{naka1}
B_1(q-a(\de))\leq f_ \de (q) \leq B_2(q-a(\de)).
\end{equation}
This can be easily done, since $(q,\de)\mapsto\f(q)/(q-a(\de))$ is
a positive continuous function on the compact region
$\{\,|\de|\leq \de_0,\;a(\de)\leq q\leq \ap-\eta\,\}$. \qed

\vskip 8pt Now, it is no trouble to establish an analogue of
Lemmas \ref{g-est-bords} and \ref{g-EST-A} with constants
independent of $\de$. We claim, without proof, that:

\begin{lem}\label{est-bords}
Let $\eta \in (0,\eta _0)$ and $M>0$ be arbitrary. Then there
exist positive constants $\de _0=\de _0(\eta,M)$ and
$C_4=C_4(\eta,M)$ such that, for any $\de \in [-\de _0,\de _0]$,
\begin{enumerate}
\item if $\xi \in [\ap-\eta,\ap+M]$, then, for all $\tau > 0$,
$Y(\tau,\xi;\de)$ remains in the interval $[\ap-\eta,\ap+M]$ and
\begin{equation}\label{est-A-bords}
|A(\tau,\xi;\de)|\leq C_4\tau \quad\hbox{for}\ \ \tau>0 \,;
\end{equation}
\item if $\xi \in [\am-M,\am+\eta]$, then, for all $\tau >0$,
$Y(\tau,\xi;\de)$ remains in the interval $[\am-M,\am+\eta]$ and
\eqref{est-A-bords} holds.
\end{enumerate}
\end{lem}

\begin{lem}\label{EST-A}
Let $\eta \in (0,\eta _0)$ be arbitrary and let $C_0$ be the
constant defined in \eqref{int1}. Then there exist positive
constants $\de_0=\de_0(\eta)$, $C_5=C_5(\eta)$ such that, for all
$\de \in [-\de _0,\de _0]$, for all $\tau>0$ and all $\xi \in
(-2C_0,2C_0)$,
$$
|A(\tau,\xi;\de)|\leq C_5(e^{\m\tau}-1).
$$
\end{lem}

\subsection{Construction of sub- and super-solutions}

We now construct a pair of sub- and super-solutions by modifying
the definition \eqref{w+-}. We set
\[
w_\ep^\pm(x,t)=Y\Big(\frac{t}{\ep^2},u_0(x)\pm\ep^2r(\pm \ep
\mathcal G,\frac{t}{\ep^2});\pm \ep \mathcal G\Big)
\]
where the function $r(\de,\tau)$ is given by
\[
r(\de,\tau)=C_6(e^{\m\tau}-1),
\]
and the constant $\mathcal G$ is chosen such that, for all small
$\ep>0$,
\[
|\g(x,t,u)| \leq \mathcal G \qquad\hbox{for} \  \ (x,t,u)\in
\overline{\Omega}\times [0,T] \times\R,
\]
which, in view of \eqref{g-est3}, is clearly possible.

\begin{lem}\label{w}
There exist positive constants  $\ep_0$ and $C_6$ such that for
all $\,\ep \in (0,\ep _0)$, $(w_\ep^-,w_\ep^+)$ is a pair of sub-
and super-solutions for problem $\Pe$, in the domain $\ombar\times
[0,\mu ^{-1} \ep^2|\ln \ep|]$, satisfying
$w^-_\ep(x,0)=w^+(x,0)=u_0(x)$.
\end{lem}

{\noindent \bf Proof.} First, the same cut-off argument as in
Subsection 3.2 enables us to assume \eqref{int2} for simplicity.
Hence $w_\ep ^\pm$ satisfy the Neumann boundary conditions. We
define an operator ${\cal L}$ by
\[
{\cal L} u:=u_t-\Delta u-\ep^{-2}(f(u)-\g(x,t,u)),
\]
and prove below that ${\cal L} w_\ep ^+ \geq 0$ by slightly
modifying the argument which we have used to prove ${\cal
L}_0w_\ep^+\geq 0$ in Section \ref{s:generation-0}.  A
straightforward calculation yields
$$
{\cal L} w_\ep^+=\frac {1}{\ep^2}\Big[Y_\tau-f(Y)+\ep \g
(x,t,Y)\Big]+Y_\xi\Big[C_6 \mu(\ep \mathcal G)e^{\mu(\ep\mathcal
G)\frac{t}{\ep ^2}} -\Delta
u_0-\displaystyle{\frac{Y_{\xi\xi}}{Y_\xi}}|\nabla u_0|^2\Big].
$$
If $\ep _0$ is sufficiently small, we note that $\pm \ep \mathcal
G \in (-\de _0,\de _0)$ and that, in the range $0 \leq t \leq \mu
^{-1} \ep ^2|\ln \ep|$,
$$
|\ep^2 C_6(e^{\mu(\pm \ep \mathcal G)t/\ep^2}-1)| \leq \ep
^2C_6(\ep^{-\mu(\pm \ep \mathcal G)/\mu}-1) \leq C_0,
$$
which implies that
$$
u_0(x)\pm \ep^2r(\pm \ep \mathcal G,\frac{t}{\ep^2}) \in
(-2C_0,2C_0).
$$
These observations allow us to use the results of the previous
subsection with the choices $\tau:=t/\ep^2$,
$\xi:=u_0(x)+\ep^2r(\ep \mathcal G,t/\ep ^2)$ and $\de:=\ep
\mathcal G$. In particular, the ordinary differential equation
\eqref{ode} yields $Y_\tau=f(Y)+\ep \mathcal G$, which implies
that
$$
{\cal L} w_\ep^+=\frac {1}{\ep}\Big[\mathcal
G+\g(x,t,Y)\Big]+Y_\xi\Big[C_6 \mu(\ep \mathcal G)\emutt-\Delta
u_0-\displaystyle{\frac{Y_{\xi\xi}}{Y_\xi}}|\nabla u_0|^2\Big].
$$
By the choice of $\mathcal G$ the first term of the right-hand
side member is positive. Using the estimate of
$A=Y_{\xi\xi}/Y_\xi$ in Lemma \ref{EST-A}, we obtain, for a
constant $C_5$ that is independent of $\ep$,
\[
\begin{array}{ll}{\cal L} w_\ep^+
&\geq Y_\xi\Big[C_6\mu(\ep \mathcal G)\emutt-|\Delta
u_0|-C_5(\emutt -1)|\nabla u_0|^2\Big] \vsp \\
&\geq Y_\xi\Big[(C_6\mu(\ep \mathcal G)-C_5 |\nabla
u_0|^2)\emutt-|\Delta u_0| +C_5 |\nabla
u_0|^2\Big]. \vsp \\
\end{array}
\]
In view of \eqref{mu}, this inequality implies that, for $\ep \in
(0,\ep _0)$, with $\ep _0$ small enough, and for $C_6$ large
enough,
\[
{\cal L} w_\ep^+ \geq \Big[C_6 \frac 12\mu-C_5{C_0}^2-C_0\Big]
\geq 0.
\]
This completes the proof of the lemma. \qed

\vskip 8pt Hence, as in Section \ref{s:generation-0}, the
comparison principle can be applied to deduce
\begin{equation}\label{g-coincee2}
w_\ep^-(x,t) \leq u^\ep(x,t) \leq w_\ep^+(x,t) \quad\ \hbox{for} \
x\in\ombar,\;0\leq t\leq \mu^{-1}\ep ^2|\ln \ep|.
\end{equation}

\subsection{Proof of Theorem \ref{g-th-gen} for the general case}

As in Subsection 3.3, we first present a key estimate of the
function $Y$ after a time interval of order $\tau\sim |\ln \ep|$.
Roughly speaking, a perturbation $\de$ of order $\ep$ does not
affect the result of Lemma \ref{after-time}.

\begin{lem}
Let $\eta \in (0,\eta _0)$ be arbitrary. Then there exist positive
constants $\ep_0$ and $C_7$ such that, for all $\ep \in
(0,\ep_0)$,
\begin{enumerate}
\item for all $\xi\in (-2C_0,2C_0)$,
\begin{equation}\label{part11}
\am-\eta \leq Y(\mu ^{-1} | \ln \ep |,\xi;\pm \ep \mathcal G) \leq
\ap+\eta,
\end{equation}
\item for all $\xi\in (-2C_0,2C_0)$ such that $|\xi-a|\geq C_7
\ep$, we have that
\begin{align}
&\text{if}\;~~\xi\geq a+C_7 \ep\;~~\text{then}\;~~Y(\mu ^{-1}| \ln
\ep |,\xi;\pm \ep \mathcal G)
\geq \ap-\eta,\label{part22}\vspace{3pt}\\
&\text{if}\;~~\xi\leq a-C_7 \ep\;~~\text{then}\;~~Y(\mu ^{-1}| \ln
\ep |,\xi;\pm \ep \mathcal G)\leq \am+\eta \label{part33}.
\end{align}
\end{enumerate}
\end{lem}

{\noindent \bf Proof.} In the sequel, by $\ep$ we always mean $\ep
\in (0,\ep _0)$, with $\ep_0=\ep_0(\eta)$ small enough. In view of
\eqref{h}, we have, for $C_7$ large enough, $a+C_7\ep \geq
a(\pm\ep \mathcal G)+\frac 1 2 C_7 \ep$. Hence for $\xi \geq
a+C_7\ep$, as long as $Y(\tau,\xi;\pm \ep \mathcal G)$ has not
reached $\ap-\eta$, we can use \eqref{est-Y-1} to deduce, as in
Section \ref{s:generation-0}, that \eqref{part22} is valid
provided that
\begin{equation*}
\tau \geq  \frac{1}{\mu (\pm \ep \mathcal G)}\ln \frac{m_0-\eta+C
\mathcal G \ep}{\frac 1 2 C_1 C_7 \ep}=:\mu ^{-1}(\ep)|\ln \ep|,
\end{equation*}
where $m_0=\max (a-\am,\ap-a)$. To complete the proof of
\eqref{part22} we must choose $C_7$ so that $\mu ^{-1}|\ln
\ep|-\mu ^{-1}(\ep)|\ln \ep| \geq 0$. A simple computation shows
that
\[
\mu ^{-1}|\ln \ep|-\mu ^{-1}(\ep)|\ln \ep| =\frac {\mu(\pm \ep
\mathcal G)-\mu}{\mu(\pm \ep \mathcal G)\mu}|\ln \ep|-\frac
{1}{\mu(\pm \ep \mathcal G)}\ln \frac{m_0-\eta+C \mathcal G\ep}
{\frac 1 2 C_1 C_7}.
\]
The first term, thanks to \eqref{mu}, is of order $\ep|\ln \ep|$.
Hence, for $C_7$ large enough, the upper quantity can be made
positive for all $\ep$. The proof of \eqref{part33} is similar and
omitted.

Next we prove \eqref{part11}. First, we can assume that the stable
zeros of $f_{\pm \ep\mathcal G}$, $\am (\pm \ep \mathcal G)$ and
$\ap(\pm \ep \mathcal G)$, are in $[\am-\eta,\ap+\eta]$. Hence, in
view of the profile of $f_{\pm \ep\mathcal G}$, if we leave from a
$\xi \in [\am-\eta,\ap+\eta]$ then $Y(\tau,\xi;\pm \ep \mathcal
G)$ will remain in $[\am-\eta,\ap+\eta]$. Now suppose that
$\ap+\eta \leq \xi \leq 2C_0$. We check below that $Y(\mu ^{-1} |
\ln \ep|,\xi;\pm \ep \mathcal G)\leq \ap+\eta$. As in Section
\ref{s:generation-0}, as long as $\ap+\eta \leq Y \leq2C_0$,
\eqref{g-pente} leads to the inequality $Y_\tau \leq p(\ap-Y)+\ep
\mathcal G$. It follows that
\[
\frac{Y_\tau}{Y-\ap}\leq -p+\ep\frac {\mathcal G}{\eta},
\]
which implies, by integration from $0$ to $\tau$, that
\[
Y(\tau,\xi;\pm \ep \mathcal G) \leq \ap+(2C_0-\ap)e^{(-p + \ep
\frac {\mathcal G}{ \eta})\tau}.
\]
One easily checks that, for $\ep$, we have $Y(\mu ^{-1}| \ln
\ep|,\xi;\pm \ep \mathcal G)\leq \ap+\eta$, which completes the
proof of \eqref{part11}. \qed

\vskip 8pt We are now ready to prove Theorem \ref{g-th-gen} in the
general case. By setting $t=\mu ^{-1} \ep ^2|\ln \ep|$ in
\eqref{g-coincee2}, we get
\begin{multline}\label{gr}
Y\Big(\mu ^{-1}|\ln \ep|, u_0(x)-\ep^2 r(-\ep \mathcal G, \mu
^{-1}|\ln \ep|);-\ep \mathcal
G\Big)\\
\leq u^\ep(x,\mu ^{-1} \ep^2|\ln \ep|) \leq Y\Big(\mu ^{-1}|\ln
\ep|, u_0(x)+\ep^2 r(\ep \mathcal G, \mu ^{-1}|\ln \ep|);+\ep
\mathcal G \Big).
\end{multline}
The point will be that, in view of \eqref{mu},
\begin{equation}\label{point}
\lim _{\ep \rightarrow 0} \frac{\mu-\mu(\pm\ep \mathcal
G)}{\mu}\ln \ep=0.
\end{equation}
It follows that
$$
\ep ^2r(\pm \ep \mathcal G,\mu ^{-1} |\ln \ep|)= C_6\ep(\ep
^{(\mu-\mu(\pm\ep \mathcal G))/\mu}-\ep) \in (\frac 1 2 C_6\ep,
\frac 3 2 C_6 \ep).
$$
Hence, as in Section \ref{s:generation-0}, the result
\eqref{g-part1} of Theorem \ref{g-th-gen} is a direct consequence
of \eqref{part11} and \eqref{gr}.

Next we prove \eqref{g-part2}. We take $x\in \om$ such that
$u_0(x)\geq a+M_0 \ep$; then
$$
\begin{array}{ll}u_0(x)-\ep^2r(-\ep \mathcal G,
\mu ^{-1}(\ep)|\ln \ep|)
&\geq a+M_0\ep-\frac 3 2 C_6 \ep\vspace{3pt}\\
&\geq a+C_7 \ep,
\end{array}
$$
if we choose $M_0$ large enough. Using \eqref{gr} and
\eqref{part22} we obtain \eqref{g-part2} which completes the proof
of Theorem \ref{g-th-gen}.\qed

\section{Motion of interface}\label{s:motion}

In Sections \ref{s:generation-0} and \ref{s:generation-g}, we have
proved that the solution $u^\ep$ develops a clear transition layer
within a very short time.  The aim of the present section is to
show that, once such a clear transition layer is formed, it
persists for the rest of time and that its law of motion is well
approximated by the interface equation $\Pz$.

Let us formulate the above assertion more clearly. By taking the
first two terms of the formal asymptotic expansion \eqref{inner},
we get a formal approximation of a solution up to order $\ep\,$:
\begin{equation}\label{ue-tilde}
\ue(x,t)\,\approx\,\tilde{u}^\ep(x,t):= \U\Big(\frac{\widetilde
d(x,t)}{\ep}\Big) +\ep U_1\Big(x,t,\frac{\widetilde
d(x,t)}{\ep}\Big).
\end{equation}
Here $U_0,\,U_1$ are as defined in \eqref{eq-phi} and
\eqref{eqU1-a}. The right-hand side has a clear transition layer
which lies exactly on $\Gamma_t$. Our goal is to show that this
function is a good approximation of a real solution; more
precisely:
\begin{quote}
{\it If $u^\ep$ becomes close to $\tilde{u}^\ep$ at some $t=t_0$,
then it stays close to $\tilde{u}^\ep$ for the rest of time.
Consequently, $\Gamma^\ep_t$ evolves roughly like $\Gamma_t$.}
\end{quote}

In order to prove this assertion, we will construct a pair of sub-
and super-solutions $u_\ep^-$ and $u_\ep^+$ for problem $\Pe$ by
slightly modifying the above function $\tilde{u}^\ep$. It then
follows that, if the solution $\ue$ satisfies
\[
u_\ep^-(x,t_0)\leq \ue(x,t_0)\leq  u_\ep^+(x,t_0),
\]
for some $t_0\geq 0$, then
\[
u_\ep^-(x,t)\leq \ue(x,t)\leq  u_\ep^+(x,t),
\]
for $t_0\leq t\leq T$, which implies that the solution $\ue$ stays
close to $\tilde{u}^\ep$.

The rest of this section is devoted to the construction of these
sub- and super-solutions. We begin with some preparations.

\subsection{A modified signed distance function}

For our later analysis, it is convenient to introduce a ``cut-off
signed distance function" $d$, which is defined as follows. First,
choose $d_0>0$ small enough so that the signed distance function
$\widetilde d$ defined in \eqref{eq:dist} is smooth in the
following tubular neighborhood of $\Gamma$:
\[
 \{(x,t) \in \overline{Q_T},\;|\widetilde{d}(x,t)|<3d_0\},
\]
and that
\begin{equation}\label{front}
 dist(\Gamma_t,\partial \Omega)\geq 3d_0 \quad \textrm{ for all }
 t\in[0,T].
\end{equation}
Next let $\zeta(s)$ be a smooth increasing function on $\R$ such
that
\[
 \zeta(s)= \left\{\begin{array}{ll}
 s &\textrm{ if }\ |s| \leq d_0\vspace{4pt}\\
 -2d_0 &\textrm{ if } \ s \leq -2d_0\vspace{4pt}\\
 2d_0 &\textrm{ if } \ s \geq 2d_0.
 \end{array}\right.
\]
We then define the cut-off signed distance function $d$ by
\begin{equation}
d(x,t)=\zeta\big(\tilde{d}(x,t)\big).
\end{equation}
Note that $|\nabla d|=1$ in the region $\{(x,t) \in
\overline{Q_T},\,|\widetilde{d}(x,t)|<d_0\}$ and that, in view of
\eqref{front}, $\n d=0$ in a neighborhood of $\partial \Omega$.
Note also that the equation of motion $\Pz$, which is equivalent
to \eqref{eq-d}, is now written as
\begin{equation}\label{interface}
 d_t=\Delta d  - \cha (x,t) \quad\
 \textrm{on}\ \; \Gamma_t,
\end{equation}
where $\cha (x,t)$ is the function defined in \eqref{gam}.

\subsection{Construction of sub- and super-solutions}

As we stated earlier, we now construct sub- and super-solutions by
modifying the function $\tilde{u}^\ep$ in \eqref{ue-tilde}.
Concerning the second term $U_1$, which is defined in
\eqref{eqU1-a}, the terms $\Delta U_1$ and $U_{1t}$ do not make
sense as we only assume that $g(\cdot,\cdot,u) \in
C^{1+\vartheta,\frac{1+\vartheta}{2}}$. In order to cope with this
lack of smoothness, we replace $U_1$ by a smooth function $\Vep$,
which is defined by
\begin{equation}\label{eqU1-ep}
\left\{\begin{array}{ll}
\VVVep+f'(U_0(z))\Vep=g^\ep(x,t,U_0(z))-\chaep(x,t){U_0}'(z),\vsp\\
\Vep(x,t,0)=0, \qquad\quad \Vep(x,t,\cdot) \in L^\infty(\R),
\end{array}\right.
\end{equation}
where
\begin{equation}\label{gam-ep}
\chaep (x,t)= c_0 (G^\ep(x,t,\ap)-G^\ep(x,t,\am)),
\end{equation}
with $G^\ep(x,t,s)=\int _a ^s g^\ep(x,t,r)dr$. Thus $U_1
^\ep(x,t,z)$ is a solution of \eqref{eq-psi} with
\begin{equation}\label{A0-ep}
A=A_0^\ep(x,t,z):=g^\ep(x,t,U_0(z))- \chaep (x,t){U_0}'(z),
\end{equation}
where the variables $x,t,\ep$ are considered parameters. Using
\eqref{g-est3} and the same arguments as in the end of Section
\ref{s:formal}, we obtain estimates analogous to \eqref{def-M} and
\eqref{def-M2}, with a constant $M$ independent of $\ep$:
\begin{equation}\label{def-M-M2-bis}
|\Vep (x,t,z)|\leq M, \quad  |\n _x \Vep (x,t,z)|\leq M.
\end{equation}
Moreover, $g^\ep$ being $C^2$ in $x$ and $C^1$ in $t$, $\Delta_x
\Vep$ and $U_{1t}^\ep$ are solutions of \eqref{eq-psi} with
$A=\Delta_x A_0^\ep$ and $A=A_{0t}^\ep$, respectively. Thus, in
view of \eqref{g-est1}, we obtain
\begin{equation}\label{delta-U1}
|\Delta_x \Vep (x,t,z)|\leq C/\ep, \quad |U_{1t}^\ep(x,t,z)|\leq
C/\ep,
\end{equation}
with some constant $C$ independent of $\ep$. Similarly,
\eqref{g-est3} and Lemma \ref{psi-decay} yield estimates analogous
to \eqref{est-psi} and \eqref{def-M3} for $\Vep$, with $C$ and $M$
independent of $\ep$:
\begin{equation}\label{est-psi-bis}
|\VVep (x,t,z)|+|\VVVep (x,t,z)|\leq Ce^{-\lambda |z|},
\end{equation}
\begin{equation}\label{def-M3-bis}
|\n _x \VVep (x,t,z)|\leq M.
\end{equation}
In the rest of this section, $C$ and $M$ will stand for the
constants that appear in inequalities
\eqref{def-M-M2-bis}--\eqref{def-M3-bis}. Note also that
\eqref{g-neumann} implies the Neumann boundary conditions
\eqref{U1-neumann} for $\Vep$.

We look for a pair of sub- and super-solutions $u_\ep^{\pm}$ for
$\Pe$ of the form
\begin{equation}\label{sub}
u_\ep^{\pm}(x,t)=U_0\Big(\frac{d(x,t) \pm \ep p(t)}{\ep}\Big)+\ep
\Vep\Big(x,t,\frac{d(x,t) \pm \ep p(t)}{\ep}\Big)\pm q(t),
\end{equation}
where
\[
\begin{array}{lll}
p(t)=-\EB+e^{Lt}+ K\vsp, \\
q(t)=\sigma \big( \beta \EB+\ep^2 Le^{Lt}\big).
\end{array}
\]
Note that $q=\sigma\ep^2\,p_t$.  It is clear from the definition
of $u_\ep^\pm$ that
\begin{equation}\label{sub-lim}
\lim_{\ep\rightarrow 0} u_\ep^\pm(x,t)= \left\{
\begin{array}{ll}
\ap &\textrm { for all } (x,t) \in Q_T^+ \vspace{4pt}\\
\am &\textrm { for all } (x,t) \in Q_T^-.\\
\end{array}\right.
\end{equation}

The main result of this section is the following:

\begin{lem}\label{fix}
Choose $\beta,\,\sigma>0$ appropriately. Then for any $K>1$, there
exist constants $\ep_0,\,L>0$ such that, for any $\ep\in(0,\ep
_0)$, the functions $(u_\ep^-,u_\ep^+)$ are a pair of sub- and
super-solutions for $\Pe$ in the domain $\ombar\times [0,T]$.
\end{lem}

\subsection{Proof of lemma \ref{fix}}\label{ss:E1-E7}

By virtue of \eqref{U1-neumann} and the fact that $\n d=0$ near
$\partial \Omega $, we have
\[
\di{\frac{\partial u_\ep^\pm}{\partial \nu}}=0 \qquad \textrm{on}
\ \ \partial \Omega \times [0,T].
\]
What we have to show is
$$
{\cal L}u_\ep^+:=(u_\ep^+)_t-\Delta u_\ep^+-\frac{1}{\ep
^2}(f(u_\ep^+) -\ep \g(x,t,u_\ep^+))\geq 0,
$$
and that ${\cal L}u_\ep^-\leq 0$.  We will prove only the former
inequality for $u_\ep^+$, since the latter follows by the same
argument.

\subsubsection {Computation of ${\cal L}u_\ep^+$}
Straightforward computations yield
$$
\begin{array}{lll}
(u_\ep^+)_t= {U_0}'(\displaystyle{\frac{d_t}{\ep}}+p_t) +\ep
U_{1t}^\ep+ \VVep(d_t+\ep p_t) +q_t\vsp \\
\nabla u_\ep^+ = {U_0}' \displaystyle{\frac{\nabla d}{\ep}} +
\ep\nabla \Vep + \VVep \nabla
d\vsp \\
\Delta u_\ep^+= {U_0}''\displaystyle{\frac{|\nabla d|^2}{\ep ^2}}
+ {U_0}'\displaystyle{\frac{\Delta d}{\ep}} + \ep\Delta \Vep+2
\nabla  \VVep\cdot \nabla d + \VVVep \displaystyle{\frac{|\nabla
d|^2}{\ep}} + \VVep \Delta  d,
\end{array}
$$
where the function $U_0$, as well as its derivatives, are
evaluated at $z=\big(d (x,t)+\ep p(t)\big)/ \ep $, whereas the
function $\Vep$, as well as its derivatives, are evaluated at
$\Big(x,t,\big(d (x,t)+\ep p(t)\big)/ \ep \Big)$.  Note that
$\nabla$ and $\Delta$ stand for $\nabla _x$ and $\Delta _x$,
respectively.  We also have
$$
\begin{array}{l}f(u_\ep^+)=f(U_0)+(\ep \Vep + q)f '(U_0) +
\di{\frac 12} (\ep \Vep +q)^2f ''(\theta) \vspace{4pt}\\
g(x,t,u_\ep^+)=g(x,t,U_0)+(\ep \Vep+q)g_u(x,t,\omega),
\end{array}
$$
where $\theta(x,t)$ and $\omega(x,t)$ are some functions
satisfying $U_0<\theta<u_\ep^+,\;U_0<\omega<u_\ep^+$. Writing
$\g=g+\g-g$ and combining the above expressions with
\eqref{eq-phi} and \eqref{eqU1-ep}, we obtain
\[
{\cal L}u_\ep^+=E_1+\cdots+E_7,
\]
where:
\vsp \\
$\qquad\quad  E_1=- \edeux q\,\Big(f'(U_0)
+\frac 12 q f ''(\theta)\Big)+{U_0}'p_t+q_t$\vsp \\
$\qquad\quad  E_2=\displaystyle{ \Big(\frac{{U_0}''}{\ep^2}
+ \frac {\VVVep}{\ep}\Big)}(1-|\n d|^2)$\vsp \\
$\qquad\quad  E_3=\displaystyle{ \Big(\frac {{U_0}'}{\ep}+
\VVep\Big)}(d_t-\Delta d +\cha)$\vsp \\
$\qquad\quad  E_4=\ep \VVep\, p_t+\di{\frac{1}{\ep}}
q\,\big(\,g_u(x,t,\omega)-\Vep f''(\theta)\,\big)$\vsp\\
$\di\qquad\quad  E_5=-\cha\, \VVep-\frac 12 \,(\Vep)^2
f''(\theta)+
\Vep g_u(x,t,\omega)-2\,\n  \VVep \cdot\n d $\vsp \\
$\qquad\quad  E_6=\ep U_{1t}^\ep-\ep \Delta  \Vep$\vsp\\
$\qquad\quad E_7=\di{\frac{1}{\ep}}(\g -g)(x,t,u_\ep^+)
-\di{\frac{1}{\ep}}(\g-g)(x,t,U_0)+\di{\frac{1}{\ep}}
(\chaep-\cha)(x,t){U_0}'\,.$\\

Before starting to estimate each of the above terms, let us
present some useful inequalities.  First, by assumption
\eqref{der-f}, there exist positive constants $b,\,m$ such that
\begin{equation}\label{bords}
f'(U_0(z))\leq -m \qquad \hbox{if} \quad U_0(z)\in
[\am,\,\am+b]\cup[\ap-b,\,\ap].
\end{equation}
On the other hand, since the region $\{z\in\R\,|\,U_0(z)\in
[\am+b,\,\ap-b] \,\}$ is compact and since ${U_0}'>0$ on $\R$,
there exists a constant $a_1>0$ such that
\begin{equation}\label{milieu}
{U_0}'(z) \geq a_1 \qquad\hbox{if} \quad U_0(z)\in
[\am+b,\,\ap-b].
\end{equation}
We set
\begin{equation}\label{beta}
\beta= \frac{m}{4}\,,
\end{equation}
and choose $\sigma$ that satisfies
\begin{equation}\label{sigma}
0< \sigma \leq \min\,(\sigma_0,\sigma_1,\sigma_2),
\end{equation}
where
\[
\sigma_0:=\frac{a_1}{\di m+F_1},\quad
\sigma_1:=\frac{1}{\beta+1},\quad \sigma
_2:=\frac{4\beta}{F_2(\beta+1)},
\]
\begin{equation*}\label{F2}
F_1:=\Vert f'\Vert _{L^\infty(\am,\ap)}, \qquad F_2:=\Vert
f''\Vert _{L ^\infty (\am-2,\ap+2)}.
\end{equation*}
Combining \eqref{bords} and \eqref{milieu}, and considering that
$\sigma \leq \sigma _0$, we obtain
\begin{equation}\label{U0-f}
U_0'(z)-\sigma f'(U_0(z))\geq \sigma m \qquad \hbox{for} \ \
-\infty<z<\infty.
\end{equation}

Now let $K>1$ be arbitrary. In what follows we will show that
${\cal L} u _\ep ^+ \geq 0$ provided that the constants $\ep_0$
and $L$ are appropriately chosen. We recall that $\am <U_0<\ap$.
We go on under the following assumption
\begin{equation}\label{ep0M}
\ep_0M\leq 1, \qquad \ep _0^2 Le^{LT} \leq 1\, .
\end{equation}
Then, given any $\ep\in(0,\ep_0)$, we have $\ep |\Vep(x,t,z)|\leq
1$ and, since $\sigma \leq \sigma _1$, $0\leq q(t)\leq 1$, so that
\begin{equation}\label{uep-pm}
\am-2\leq u_\ep^\pm(x,t) \leq \ap+2\, .
\end{equation}

\subsubsection {The term $E_1$}

Direct computation gives
$$
E_1=\frac{\beta}{\ep^2}\,\EB(I-\sigma\beta)+Le^{Lt}(I+\ep^2\sigma
L),
$$
where
$$
I=U_0'-\sigma f '(U_0)-\frac {\sigma^2}2 f
''(\theta)(\beta\EB+\ep^2 Le^{Lt}).
$$
In virtue of \eqref{U0-f} and \eqref{uep-pm}, we have
\[
I\geq \sigma m-\frac {\sigma^2}{2} F_2(\beta+\ep^2 Le^{LT}).
\]
Combining this, \eqref{ep0M} and the inequality $\sigma \leq
\sigma _2$, we obtain $ I \geq 2\sigma\beta$. Consequently, we
have
$$
E_1\geq \frac{\sigma\beta^2}{\ep^2}\EB + 2\sigma\beta L e^{Lt}.
$$

\subsubsection {The term $E_2$}

First, in the region where $|d|\leq d_0$, we have $|\n d|=1$,
hence $E_2=0$. Next we consider the region where $|d|\geq d_0.$ We
deduce from Lemma \ref{est-phi} and from \eqref{est-psi-bis} that
:
$$
|E_2|\leq C(\frac{1}{\ep^2}+\frac{1}{\ep})e^{-\lambda|d+\ep p|/
\ep}\leq \frac{2C}{\ep^2}e^{-\lambda(d_0 / \ep-|p|)}.
$$
We remark that $0<K-1 \leq p \leq e^{LT} +K$. Consequently, if we
assume
\begin{equation}\label{ga}
e^{LT}+K \leq \frac{d_0}{2\ep_0},
\end{equation}
then $\displaystyle{\frac{d_0}{\ep}}-|p|\geq
\displaystyle{\frac{d_0}{2\ep}}$, so that
$$
|E_2|\leq \frac{2C}{\ep^2}e^{-\lambda d_0 / (2\ep)} \leq C_2
:=\frac{32C}{(e\lambda d_0)^2}.
$$

\subsubsection {The term $E_3$}

By \eqref{interface} and \eqref{gam}, we have
$$
(d_t-\Delta d +\cha)(x,t)=0 \qquad \textrm{on} \quad \Gamma_t=\{x
\in \om,\; d(x,t)=0\}.
$$
Since $\gamma$ is of class $C^{1+\vartheta,
\frac{1+\vartheta}{2}}$ by virtue of \eqref{g-est3}, we see that
the interface $\Gamma_t$ is of class $C^{3+\vartheta,
\frac{3+\vartheta}{2}}$. Therefore both $\Delta d$ and $d_t$ are
Lipschitz continuous near $\Gamma_t$.  It follows that there
exists a constant $N>0$ such that:
$$
|(d_t-\Delta d +\cha )(x,t)|\leq N|d(x,t)| \quad \textrm{ for all
}  (x,t) \in Q_T.
$$
Applying  Lemma \ref{est-phi} and the estimate \eqref{est-psi-bis}
we deduce that
$$
\begin{array}{lll}
|E_3|&\leq 2NC\displaystyle{\frac
{|d|}{\ep}}e^{-\lambda| d/ \ep +p|}\vsp \\
&\leq 2NC \max_{\xi \in \R }|\xi|e^{-
\lambda|\xi +p|}\vsp \\
&\leq 2NC\max (|p|,\di{\frac 1 \lambda}).
\end{array}
$$
Thus, recalling that $|p|\leq e^{Lt}+K$, we obtain
$$
|E_3|\leq C_3(e^{Lt}+K)+{C_3}', $$ where $C_3:=2NC$ and
${C_3}':=2NC/\lambda$.

\subsubsection {The term $E_4$}

In view of \eqref{g-est2} and \eqref{est-psi-bis}, both $g_u$ and
$|\VVep|$ are bounded by some constant $C$. Hence, substituting
the expression for $p_t$ and $q$, we obtain
\[
|E_4| \leq C_4\big ( \frac 1 \ep \beta \EB+\ep L e^{Lt}\big ),
\]
where $C_4:=C+\sigma(C+MF_2)$.

\subsubsection {The term $E_5$}

In view of \eqref{gam}, the term $|\cha|$ is bounded by
$c_0(\ap-\am)C$ on $\ombar \times [0,T]$. Using \eqref{g-est2} and
\eqref{def-M3-bis}, we easily obtain $|E_5|\leq C_5$, where $C_5$
depends only on $C$, $M$, $F_2$.

\subsubsection {The term $E_6$}
We use \eqref{delta-U1} to deduce that $|E_6|\leq 2C=:C_6$.

\subsubsection {Finally the term $E_7$}
We recall that $|g^\ep-g| \leq C\ep$ so that $|\chaep -\cha | \leq
c_0(\ap-\am)C \ep$. It then follows that
\[
 |E_7|\le2C+Cc_0(\ap-\am)=:C_7.
\]

\subsubsection {Completion of the proof}

Collecting all these estimates gives
\begin{equation}\label{september}
 {\cal L} u_\ep^+\geq (\frac{\sigma \beta ^2}{\ep^2}-\frac{C_4\beta}{\ep})
 e^{-\beta t/\ep ^2}+ (2\sigma \beta L-C_3-\ep C_4 L)e^{Lt}-C_8,
\end{equation}
where $ C_8:=C_2+KC_3+{C_3}'+C_5+C_6+C_7$. Now we set
\[
 L:=\frac 1 T\ln \frac {d_0}{4\ep _0},
\]
which, for $\ep_ 0$ small enough, validates assumptions
\eqref{ep0M} and \eqref{ga}. For $\ep_0$ small enough, the first
term of the right-hand side of \eqref{september} is positive,
hence
\[
 {\cal L} u_\ep^+\geq \big[\sigma \beta L-C_3]e^{Lt}-C_8 \geq  \frac 1 2 \sigma \beta L -C_8 \geq 0.
\]
The proof of Lemma \ref{fix} is now complete, with the choice of
the constants $\beta, \sigma$ as in \eqref{beta},
\eqref{sigma}.\qed

\section{Proof of the main results }\label{s:proof}

\subsection{ Proof of Theorem \ref{width}}\label{ss:proof-w}

Let $\eta \in (0,\eta _0)$ be arbitrary. Choose $\beta$ and
$\sigma$ that satisfy \eqref{beta}, \eqref{sigma} and
\begin{equation}\label{eta}
\sigma \beta \leq \frac \eta 3.
\end{equation}
By Theorem \ref{g-th-gen}, there exist positive constants $\ep_0$
and $M_0$ such that \eqref{g-part1}, \eqref{g-part2} and
\eqref{g-part3} hold with the constant $\eta$ replaced by $ \sigma
\beta /2$. Since $\n u_0 \cdot n \neq 0$ everywhere on $\Gamma
_0=\{x\in\om, \; u_0(x)=a\}$ and since $\Gamma _0$ is a compact
hypersurface, we can find a positive constant $M_1$ such that
\begin{equation}\label{corres}
\begin{array}{ll}\text { if } \quad d_0 (x) \geq \ M_1 \ep
&\text { then } \quad u_0(x) \geq a +M_0 \ep\vspace{3pt}\\
\text { if } \quad d_0 (x) \leq -M_1 \ep & \text { then } \quad
u_0(x) \leq a -M_0 \ep.
\end{array}
\end{equation}
Here $d_0(x):=\tilde d(x,0)$ denotes the signed distance function
associated with the hypersurface $\Gamma_0$. Now we define
functions $H^+(x), H^-(x)$ by
\[
\begin{array}{l}
H^+(x)=\left\{
\begin{array}{ll}
\ap+\sigma\beta/2\quad\ &\hbox{if}\ \ d_0(x)\geq -M_1\ep\\
\am+\sigma\beta/2\quad\ &\hbox{if}\ \ d_0(x)<  -M_1\ep,
\end{array}\right.
\vsp\\
H^-(x)=\left\{
\begin{array}{ll}
\ap-\sigma\beta/2\quad\ &\hbox{if}\ \ d_0(x)\geq \;M_1\ep\\
\am-\sigma\beta/2\quad\ &\hbox{if}\ \ d_0(x)<  \;M_1\ep.
\end{array}\right.
\end{array}
\]
Then from the above observation we see that
\begin{equation}\label{H-u}
H^-(x) \,\leq\, u^\ep(x,\mu^{-1} \ep^2|\ln \ep|) \,\leq\,
H^+(x)\qquad \hbox{for}\ \ x\in\Omega.
\end{equation}

Next we fix a sufficiently large constant $K>0$ such that
\begin{equation}\label{K}
U_0(-M_1+K) \geq \ap-\frac {\sigma \beta}{3} \quad \text { and }
\quad U_0(M_1-K) \leq \am+\frac {\sigma \beta}{3}.
\end{equation}
For this $K$, we choose $\ep _0$ and $L$  as in Lemma \ref{fix}.
We claim that
\begin{equation}\label{uep-H}
u_\ep^-(x,0)\leq H^-(x),\quad\ H^+(x)\leq u_\ep^+(x,0) \qquad
\hbox{for} \ \ x\in\Omega.
\end{equation}
We only prove the former inequality, as the proof of the latter is
virtually the same. Then it amounts to showing that
\begin{equation}\label{c3}
u_\ep ^- (x,0)=U_0\big(\frac {d_0(x)}{\ep}-K\big)+\ep
\Vep\big(x,0, \frac{d_0(x)}{\ep}-K\big)-\sigma (\beta+\ep ^2 L)
\;\leq\; H^-(x).
\end{equation}
By \eqref{def-M-M2-bis} we have $|\Vep|\leq M$. Therefore, by
choosing $\ep_0$ small enough so that $\ep_0 M \leq
\sigma\beta/6$, we see that
\begin{align*}
u_\ep ^- (x,0) &\leq\;
U_0\big(\frac {d_0(x)}{\ep}-K\big)+\ep M-\sigma(\beta+\ep ^2 L)\\
& \leq\;U_0\big(\frac {d_0(x)}{\ep}-K\big)-\frac{5}{6}\sigma\beta.
\end{align*}
In the range where $d_0(x) < M_1 \ep$, the second inequality in
\eqref{K} and the fact that $U_0$ is an increasing function imply
\[
U_0\big(\frac {d_0(x)}{\ep}-K\big)-\frac{5}{6}\sigma\beta \;\leq\;
\am-\frac {\sigma \beta}{2}\;=\;H^-(x).
\]
On the other hand, in the range where $d_0(x) \geq M_1 \ep$, we
have
\[
U_0\big(\frac {d_0(x)}{\ep}-K\big)-\frac{5}{6}\sigma\beta \;\leq\;
\ap-\frac{5}{6}\sigma\beta \;\leq\;H^-(x).
\]
This proves \eqref{c3}, hence \eqref{uep-H} is established.

Combining \eqref{H-u} and \eqref{uep-H}, we obtain
$$
u_\ep^-(x,0)\leq u^\ep(x,\mu ^{-1} \ep ^2|\ln \ep|) \leq
u_\ep^+(x,0).
$$
Since $u_\ep^-$ and $u_\ep^+$ are sub- and super-solutions of
$\Pe$ thanks to Lemma \ref{fix}, the comparison principle yields
\begin{equation}\label{ok}
u_\ep^-(x,t) \leq u^\ep (x,t+t^\ep) \leq u_\ep^+(x,t) \quad \text
{ for } 0 \leq t \leq T-t^\ep,
\end{equation}
where $t^\ep=\mu ^{-1} \ep ^2|\ln \ep|$. Note that, in view of
\eqref{sub-lim}, this is enough to prove Corollary \ref{total}.
Now let $C$ be a positive constant such that
\begin{equation}\label{C}
U_0(C-e^{LT}-K) \geq \ap-\frac \eta 2 \quad \text { and } \quad
U_0(-C+e^{LT}+K) \leq \am+\frac \eta 2.
\end{equation}
One then easily checks, using \eqref{ok} and \eqref{eta}, that,
for $\ep _0$ small enough, for $0\leq t \leq T-t^\ep$, we have
\begin{equation}\label{correspon}
\begin{array}{ll}\text { if } \quad d(x,t) \geq \ C \ep &
\text { then } \quad
u^\ep (x,t+t^\ep) \geq \ap -\eta\vspace{3pt}\\
\text { if } \quad d(x,t) \leq -C \ep & \text { then } \quad u^\ep
(x,t+t^\ep) \leq \am +\eta,
\end{array}
\end{equation}
and
$$
u^\ep (x,t+t^\ep) \in [\am-\eta,\ap+\eta],
$$
which completes the proof of Theorem \ref{width}.\qed

\subsection{Proof of Theorem \ref{error}}

In the case where $ \mu ^{-1} \ep ^2 |\ln \ep|\leq t \leq T$, the
assertion of the theorem is a direct consequence of Theorem
\ref{width}.  Thus, all we have to consider is the case where
$0\leq t\leq \mu ^{-1} \ep ^2 |\ln \ep|$.  We first need the
following lemma concerning $Y$, the solution of the perturbed
ordinary differential equation \eqref{ode}.

\begin{lem}\label{est-Y-debut}
There exists a constant $C_8 >0$ such that
\begin{equation}\label{debut}
\begin{array}{ll}\text { if } \quad \xi \geq a+C_8 \ep &
\text{then} \quad Y(\tau,\xi;\pm \ep \mathcal G)>a\quad \hbox{for}
\quad 0\leq \tau \leq \mu ^{-1}|\ln \ep|,
\vspace{3pt}\\
\text { if } \quad \xi \leq a-C_8 \ep & \text {then} \quad
Y(\tau,\xi;\pm \ep \mathcal G)<a \quad \hbox{for} \quad 0\leq \tau
\leq \mu ^{-1}|\ln \ep|.
\end{array}
\end{equation}
\end{lem}

{\noindent \bf Proof.} We only prove the first inequality. In view
of estimates \eqref{est-Y-1} and \eqref{h}, we obtain, for $\xi
\geq a+C_8 \ep$,
\[
\begin{array}{llll}Y(\tau,\xi;\pm \ep \mathcal
G)&\geq a(\pm \ep \mathcal G)+C_1 e^{\mu(\pm \ep \mathcal
G)\tau}(a+C_8 \ep-a(\pm \ep \mathcal G)) \vsp \\ &\geq a-C\mathcal
G \ep+C_1(-C\mathcal G \ep+C_8 \ep)\vsp \\ &\geq
a+\ep(C_1C_8-C\mathcal G(C_1 +1))\vsp \\ &>a,
\end{array}
\]
if we choose $C_8$ large enough. \qed

\vskip 8pt Now we turn to the proof of Theorem \ref{error}. We
first claim that there exists a positive constant $M_2$ such that
for all $t \in [0,\mu ^{-1} \ep ^2 |\ln \ep|]$,
\begin{equation}\label{claim1}
\Gamma _t ^\ep \subset \mathcal N _{M_2\ep} (\Gamma _0).
\end{equation}
To see this, we choose $M_0$ large enough, so that $M_0 \geq C_8 +
2 C_6$ holds in addition to \eqref{g-part1}, \eqref{g-part2} and
\eqref{g-part3}. We then choose $M_2 > M_1$, where $M_1$ is as
defined in \eqref{corres}. In view of this last condition, we see
that if $\ep _0$ is small enough and if $d_0(x)\geq M_2 \ep$, then
for $0 \leq t \leq \mu^{-1} \ep ^2 | \ln \ep|$,
$$
\begin{array}{llll}u_0(x)-\ep ^2 r(-\ep \mathcal G,\di{\frac
{t}{\ep ^2}})&\geq a+M_0 \ep-\ep^2 C_6\big[e^{\mu(-\ep \mathcal
G)|\ln \ep|/ \mu}-1\big]\vsp \\ &\geq a+\ep\big[M_0-C_6 \ep
^{(\mu-\mu(\pm\ep \mathcal G))/\mu}+\ep
C_6\big]\vsp \\
&\geq a+\ep(M_0 -2 C_6) \hspace{32pt} \big(\ \leftarrow \text {
thanks to \eqref{point} }\big)\vsp
\\
&\geq a+C_8 \ep.
\end{array}
$$
This inequality and Lemma \ref{est-Y-debut} imply
$w_\ep^-(x,t)>a$, where $w_\ep^-$ is the sub-solution defined in
\eqref{w+-}.  Consequently, by \eqref{g-coincee1},
\[
u^\ep(x,t)>a\qquad \hbox{if} \ \ d_0(x)\geq M_2 \ep.
\]
In the case where $d_0(x)\leq -M_2\ep$, similar arguments lead to
$u^\ep (x,t)<a$. This completes the proof of \eqref{claim1}. Note
that we have proved that, for all $0 \leq t \leq \mu ^{-1} \ep ^2
|\ln \ep|$,
\begin{equation}\label{rost1}
\begin{array}{ll}
u^\ep(x,t) >a \quad \text {if} \quad x \in \Omega _0 ^+ \setminus
\mathcal N _{M_2\ep}(\Gamma _0),\vspace{3pt}\\
u^\ep(x,t) <a \quad \text {if} \quad x \in \Omega _0 ^- \setminus
\mathcal N _{M_2\ep}(\Gamma _0).
\end{array}
\end{equation}
Since $\Gamma _t$ depends on $t$ smoothly, there is a constant
$\tilde C >0$ such that, for all $t \in [0,\mu ^{-1} \ep ^2 |\ln
\ep|]$,
\begin{equation}\label{claim2}
\Gamma _0 \subset \mathcal N _{\tilde C\ep ^2|\ln \ep|} (\Gamma
_t),
\end{equation}
and
\begin{equation}\label{rost2}
\begin{array}{ll}
\Omega _t ^+ \setminus \mathcal N _{\tilde C \ep}(\Gamma _t)
\subset
\Omega _0 ^+ \setminus \mathcal N _{M_2\ep}(\Gamma _0),\vspace{3pt}\\
\Omega _t ^- \setminus \mathcal N _{\tilde C \ep}(\Gamma _t)
\subset \Omega _0 ^- \setminus \mathcal N _{M_2\ep}(\Gamma _0).
\end{array}
\end{equation}
As a consequence of \eqref{claim1} and \eqref{claim2} we get
$$
\Gamma _t ^\ep \subset \mathcal N _{M_2 \ep+\tilde C \ep ^2|\ln
\ep|} (\Gamma _t) \subset \mathcal N _{C \ep} (\Gamma _t),
$$
which completes the proof of Theorem
\ref{error}.\qed\\

{\bf \noindent Proof of Corollary \ref{total-2}.} In view of
Theorem \ref{error} and the definition of the Hausdorff distance,
to prove this corollary we only need to show that
\begin{equation}\label{obj}
\Gamma _t \subset \mathcal N _{C' \ep} (\Gamma _t ^\ep)\quad \text
{ for } \quad 0 \leq t \leq T,
\end{equation}
for some constant $C'>0$. To that purpose let $C'$ be a constant
satisfying $C'>\max (\tilde C,C)$, where $C$ is as in Theorem
\ref{width} and $\tilde C$ as in \eqref{rost2}. Choose $t \in
[0,T]$ and $x_0 \in \Gamma _t$ arbitrarily and, $n$ being the
Euclidian normal vector exterior to $\Gamma _t$ at point $x_0$,
define a pair of points:
$$
x^+:=x_0+C'\ep n \quad \text{and}\quad  x^-:=x_0-C' \ep n.
$$
Since $C'>C$ and since the curvature of $\Gamma _t$ is uniformly
bounded as $t$ varies over $[0,T]$, we see that
$$
x^+ \in \Omega ^+ _t \setminus \mathcal N _{C\ep}(\Gamma _t) \quad
\text {and} \quad  x^- \in \Omega ^- _t \setminus \mathcal N
_{C\ep}(\Gamma _t),
$$
if $\ep$ is sufficiently small. Therefore, if $t\in [\mu ^{-1} \ep
^2|\ln \ep|,T]$, then, by Theorem \ref{width}, we have
\begin{equation}\label{rostand}
u^\ep(x^-,t)<a<u^\ep(x^+,t).
\end{equation}
On the other hand, if $t\in[0,\mu ^{-1} \ep ^2|\ln \ep|]$, then
from \eqref{rost1}, \eqref{rost2} and the fact that $C'>\tilde C$,
we again obtain \eqref{rostand}. Thus \eqref{rostand} holds for
all $t\in[0,T]$. Now, by the mean value theorem, we see that for
each $t\in[0,T]$ there exists a point $x_1$ on the line segment
$[x^-,x^+]$ such that $u^\ep(x_1,t)=a$. This implies $x_1 \in
\Gamma _t ^\ep$. Furthermore we have $d(x_0,x_1)\leq C' \ep$,
since $x_1$ lies on the line segment $[x^-,x^+]$. This proves
\eqref{obj}.\qed

\section{Application to reaction-diffusion systems}
\label{s:RD}

In this section we discuss the singular limit of the
reaction-diffusion system $\RD$ and prove Theorems \ref{width-N},
\ref{error-N} and their corollaries.  Our strategy is to regard
the first equation of $\RD$ as a perturbed Allen-Cahn equation and
apply what we have already proved for this equation.

\subsection{Preliminaries: global existence}\label{ss:global}

Before studying the singular limit of $\RD$, we first show that
the solution of this system exists globally for $t\geq 0$,
provided that $\ep$ is sufficiently small.  Recall that the system
$\RD$ is written in the form
\[
\begin{cases}
\,u_t=\Delta u+\edeux\,
\big(f(u)+\ep\,f_1(u,v)+O(\ep^2)\big),\vspace{4pt}\\
\,v_t=D \Delta v +h(u,v),
\end{cases}\hspace{30pt}
\]
where $h(u,v)$ satisfies the hypothesis ({\bf H}).  The standard
parabolic theory guarantees the existence of local solutions for
$\RD$.  In order to prove that the solution exists globally for
$t\geq 0$, it suffices to show that the solution remains uniformly
bounded.  This will be done by using the well-known method of
invariant rectangles.

Given arbitrary $u_0,\,v_0\in C(\ombar)$, we choose a constant
$L>0$ such that
\begin{equation}\label{L}
f(-L)>0>f(L),\qquad -L\leq u_0(x)\leq L\ \ \ \hbox{for} \ \
x\in\ombar.
\end{equation}
Such a constant $L$ exists since $f(u)>0$ for $u<\am$ and $f(u)<0$
for $u>\ap$.  By hypothesis ({\bf H}), we can choose a constant
$M_1$ satisfying
\[
M_1\geq \Vert v_0\Vert_{L^\infty(\om)},
\]
along with the condition \eqref{h(u,v)}, namely
\begin{equation}\label{side-v}
 h(u,-M_1)\geq 0\geq h(u,M_1) \qquad\hbox{for}\ \ |u|\leq L.
\end{equation}

Now we consider the rectangle
\[
 {\cal R}:=\{\,(u,v)\in\R^2\,\big|\, |u|\leq L,\,|v|\leq M_1\,\}.
\]
It follows from \eqref{L} that, for all sufficiently small
$\ep>0$,
\begin{equation}\label{side-u}
 f^\ep(-L,v)>0>f^\ep(L,v)\qquad\hbox{for}\ \ |v|\leq M_1.
\end{equation}
The inequalities \eqref{side-v} and \eqref{side-u} imply that the
rectangle ${\cal R}$ is a positively invariant region for the
system of ordinary differential equations
\[
\begin{cases}
\,u_t=\di{\frac{1}{\ep ^2}}f^\ep(u,v),\vspace{5pt}\\
\,v_t=h(u,v),
\end{cases}
\]
since the vector field $(\ep^{-2}f^\ep(u,v),h(u,v))$ points
inwards everywhere on the boundary of ${\cal R}$.  The maximum
principle then implies that ${\cal R}$ is also positively
invariant for the system $\RD$. Consequently, since
$(u_0(x),v_0(x))\in{\cal R}$ for $x\in\ombar$, we have
\[
(u(x,t),v(x,t))\in {\cal R}\qquad\hbox{for}\ \ x\in\ombar,\;t\geq
0,
\]
so long as the solution is defined. This uniform bound then
implies that the solution exists globally for $t\geq 0$.

In the case of equations for which only nonnegative solutions are
to be considered (see Remark \ref{rm:positive}), we can argue just
similarly, by replacing ${\cal R}$ by the rectangle ${\cal
R}_+:=\{(u,v)\,|\,0\leq u\leq L,\,0\leq v\leq M_1\}$. Summarizing,
we have proved the following proposition:

\begin{prop}\label{pr:global}
Let $(u_0,v_0)\in C(\ombar)\times C(\ombar)$.  In the case where
the conditions of Remark \ref{rm:positive} apply, assume further
that $u_0, v_0\geq 0$.  Then there exists $\ep_0>0$ such that for
any $\ep\in(0,\ep_0)$, the solution of $\RD$ exists globally for
$t\geq 0$ and is uniformly bounded.
\end{prop}

\vskip 5pt
\begin{rem}
For the details of the method of invariant rectangles, we refer
the reader to the book \cite[Chapter 14, Corollary 14.8]{S}. See
also \cite{CCS} and \cite{CH}.  It should be noted that \cite{KMY}
makes a much earlier study of invariant rectangles for a
finite-difference scheme for reaction diffusion systems.
\end{rem}

\subsection{Proof of the main results}
\label{ss:proof-results-RD}

Now we turn to the reaction-diffusion system $\RD$ and explain our
strategy for proving Theorems \ref{width-N}, \ref{error-N} and
their corollaries.

In what follows, we fix the initial data $(u_0,\,v_0)$ and denote
by $(\ue,\,\ve)$ the solution of the system $(RD^\ep)$.  The
solution of the associated moving boundary problem $(RD ^0)$ will
be denoted by $(\Gamma, \tilde v)$.

Given a function $v(x,t)$ on $\ombar\times[0,+\infty)$, we set
\begin{equation}\label{ge-g}
\begin{array}{rl}
g^\ep[v](x,t,u)&:=-f_1(u,v(x,t))-\ep\,\fee(u,v(x,t)),
\vspace{5pt}\\
g[v](x,t,u)&:=-f_1(u,v(x,t)),
\end{array}
\end{equation}
where $f_1,\,\fee$ are as in \eqref{f(u,v)}. The first equation of
$\RD$ is then written in the form
\begin{equation}\label{RD-1}
u_t=\Delta u + \edeux \big(f(u) - \ep\,g^\ep[v](x,t,u)\big),
\end{equation}
so that $\ue (x,t)$ is the solution of $\Pe$ with the choice of
the perturbation term $g^\ep (x,t,u)=g^\ep[\ve](x,t,u)$. On the
other hand, the equation of surface motion in the limit problem
$\RDz$ is written in the form
\begin{equation}\label{RDz-1}
 \di \, V_{n}=-(N-1)\kappa + c_0
 \int^{\ap}_{\am} g[\tilde{v}](x,t,r)\,dr
 \quad \text { on } \Gamma_t ,
\end{equation}
so that $\Gamma$ is the solution of $(P^0)$ with the choice
$g(x,t,u)=g[\tilde v](x,t,u)$.

Thus Theorems \ref{width-N}, \ref{error-N} and their corollaries
will follow from what we have shown for the single equation $\Pe$.
In order for Theorems \ref{width} and \ref{error} for $\Pe$ to be
applicable to the present reaction-diffusion system $(RD ^\ep)$,
all we have to do is to verify the conditions \eqref{g-est1} to
\eqref{g-est4}. More precisely, we have to show that, for all
small $\ep>0$,
\[
 |\Delta_{x} \g[\ve](x,t,u)| \leq C \ep^{-1} \quad \text{ and }
 \quad | \partial_t\,\g[\ve](x,t,u)| \leq C \ep^{-1},
\]
\[
 |{\partial_u}\,\g [\ve](x,t,u)| \leq C,
\]
\[
 \Vert \g[\ve](\cdot,\cdot,u)
 \Vert_{C^{1+\vartheta,\frac{1+\vartheta}{2}}(\ombar\times[0,T])}
 \leq C,
\]
\[
 \big|\g [\ve](x,t,u)-g[\tilde{v}](x,t,u)\big| \leq C \ep.
 \]
Since $g[v],\,\g[v]$ are defined by \eqref{ge-g} and since $f_1,\,
\fee$ are smooth, it suffices to prove the following estimates for
some $C>0$ and for all small $\ep>0$:
\begin{equation}\label{v-est1}
 |\Delta_{x} \ve(x,t)| \leq C \ep^{-1} \quad \text{ and }
 \quad | \partial_t \ve(x,t)| \leq C \ep^{-1},
\end{equation}
\begin{equation}\label{v-est3}
 \Vert \ve
 \Vert_{C^{1+\vartheta,\frac{1+\vartheta}{2}}(\ombar\times[0,T])}
 \leq C,
\end{equation}
\begin{equation}\label{v-est4}
 |\ve(x,t)-\tilde{v}(x,t)| \leq C \ep.
\end{equation}

The estimates \eqref{v-est1} and \eqref{v-est3} are elementary,
but \eqref{v-est4} requires far more elaborate analysis. In this
subsection we prove \eqref{v-est1}, \eqref{v-est3} and give an
outline of the proof of \eqref{v-est4}.  A full proof of
\eqref{v-est4} will be given later.

\begin{proof}[{\bf Proof of \eqref{v-est3} and \eqref{v-est1}}]

Since $\ve$ satisfies
\begin{equation}\label{ve}
\ve_t=D\Delta \ve+h(\ue,\ve)\quad\ \textrm{in}\ \
\Omega\times(0,T],
\end{equation}
along with the Neumann boundary conditions, it can be expressed as
\begin{equation}\label{ve-G}
\ve(x,t)\;=\;I_1+I_2,
\end{equation}
where \vspace{-10pt}
\[
\begin{split}
I_1&:=\int_{\om}G(x,y,t)v_0(y)dy,\\
I_2& :=\int_0^t\int_{\om}G(x,y,t-s)\;\!
h(\ue(y,s),\ve(y,s))\,dy\;\!ds,\\
\end{split}
\]
with $G(x,y,t)$ being the fundamental solution for equation
$v_t=D\Delta v$ under the Neumann boundary conditions. Since
$h(\ue,\ve)$ is uniformly bounded, standard estimates of
$G(x,y,t)$ imply \eqref{v-est3} for any $\vartheta\in (0,1)$.

In the mean while, the same rescaling argument as in Remark
\ref{rm:thickness} yields
\begin{equation}\label{resc-arg}
\Vert\ue\Vert_{C^{\vartheta,\frac{\vartheta}{2}}(\ombar\times[0,T])}
\leq C\ep^{-\vartheta}.
\end{equation}
Indeed, since $\n _y u$, $u_ \tau$ are bounded, where $y=x/\ep$,
$\tau=t/\ep^2$, we have $\n_x u=O(1/\ep)$, $u_t=O(1/\ep^2)$.
Consequently we have
\[
\begin{array}{rl}
\di{\frac{|u(x,t)-u(x',t')|}{|x-x'|^{\vartheta}+|t-t'|^{\vartheta/2}}}
&\leq\di{\frac{|u(x,t)-u(x',t)|}{|x-x'|^{\vartheta}}}+
\di{\frac{|u(x',t)-u(x',t')|}{|t-t'|^{\vartheta/2}}}
\vspace{5pt}\\
&\leq |u(x,t)-u(x',t)|^{1-\vartheta}\di{\frac{|u(x,t)-u(x',t)
|^{\vartheta}}{|x-x'|^{\vartheta}}}\vspace{5pt}\\
& \ +|u(x',t)-u(x',t')|^{1-\vartheta/2}
\di{\frac{|u(x',t)-u(x',t')|^{\vartheta/2}}{|t-t'|^{\vartheta/2}}}
\vspace{5pt}\\
&\leq(2\Vert u\Vert _{L^\infty})^{1-\vartheta}\Vert \n_xu
\Vert_{L^\infty}^{\vartheta}+(2\Vert u\Vert_{L^\infty})
^{1-\vartheta/2}\Vert u_t\Vert_{L^\infty}^{\vartheta/2}
\vspace{6pt}\\
&\leq C \ep^{-\vartheta}.
\end{array}
\]
Combining \eqref{resc-arg} and \eqref{v-est3}, we see that $\Vert
h(\ue,\ve)\Vert_{C^{\vartheta,\frac{\vartheta}{2}}(\ombar
\times[0,T])}\leq C\ep^{-\vartheta}$, hence, by the Schauder
estimate,
\[
\Vert I_2
 \Vert_{C^{2+\vartheta,1+\frac{\vartheta}{2}}(\ombar\times[0,T])}
 \leq C\ep^{-\vartheta}.
\]
Here the constant $C$ may depend on the choice of
$\vartheta\in(0,1)$. On the other hand, $I_1$ is bounded in
$C^{2,1}(\ombar\times[0,T])$ since $v_0\in C^2(\ombar)$. Combining
these, we obtain $|\Delta_{x} \ve(x,t)| = O(\ep^{-\vartheta})$,
hence $O(\ep^{-1})$. Substituting this into \eqref{ve} yields the
second inequality in \eqref{v-est1}.
\end{proof}

%%%%%%%%%%%%%%%%
\noindent{\bf Outline of the proof of \eqref{v-est4}.} We decouple
the system $(RD^\ep)$ as follows. As mentioned earlier, $\ue
(x,t)$ is the solution of $\Pe$ with the choice of the
perturbation term $g^\ep(x,t,u)=g^\ep [v^\ep](x,t,u)$, that is,
\[
 (\star)\quad\begin{cases}
 u_t=\Delta u+\edeux (\F)\\
 \dudn = 0\vspace{3pt}\\
 u(x,0)=u_0(x).
 \end{cases}
\]
Once the solution $\ue$ is determined, $\ve$ is the solution of
the problem
\begin{equation*}
(\star\star)\quad \begin{cases}
\,v_t=D \Delta v + h(u,v)\vspace{4pt}\\
\,\di\frac{\partial v}{\partial \nu}= 0\vspace{4pt}\\
\,v(x,0)=v_0 (x),
\end{cases}
\end{equation*}
with the choice $u=\ue$. This means that $\ve (x,t)$ is a fixed
point of the following map $\Phi ^\ep:=\Phi _2 \circ \Phi _1
^\ep$\,:
\begin{equation*}
\Phi^\ep:v \; \xrightarrow {\;\Phi _1 ^\ep \; \text{via}\;
(\star)\;} \ \bar u^\ep \ \xrightarrow {\;\Phi _2 \; \text{via}\;
(\star\star)\;}\ \bar v^\ep,
\end{equation*}
where $\Phi _1 ^\ep$ maps a function $v(x,t)$ to the solution
$\bar u ^\ep (x,t)$ of $(\star)$ for the choice
$g^\ep(x,t,u)=g^\ep [v](x,t,u)$, and $\Phi _2$ maps a function
$\bar u (x,t)$ to the solution $\bar v (x,t)$ of $(\star\star)$
for the choice $u=\bar u$.

On the other hand, as for the limit problem $(RD ^0)$, the
solution $\tilde v (x,t)$ can be regarded as a fixed point of the
map $\Phi ^0:=\Phi _2 \circ \Phi _1 ^0$\,:
\begin{equation*}
\Phi ^0:v \; \xrightarrow {\;\Phi _1 ^0\; \text{via}\;
(\star\star\star)\;}\ \hat u \ \xrightarrow {\;\Phi _2\;
\text{via}\;(\star\star)\;}\ \hat v,
\end{equation*}
where $\Phi _1 ^0$ maps a function $v(x,t)$ to the step function
\begin{equation*}
 \hat u(x,t)=\begin{cases}
 \, \ap &\text{in}\ \; \om^+ (\Gamma _t[v])\\
 \, \am &\text{in}\ \; \om^- (\Gamma _t[v]),
 \end{cases}
\end{equation*}
where $\Gamma _t[v]$ is the solution of the equation of surface
motion
\begin{equation*}
(\star\star\star)\quad \begin{cases}
 \di \, V_{n}=-(N-1)\kappa -c_0 \int _{\am} ^{\ap} f_1 (r,v(x,t))dr
 \quad \text { on } \Gamma_t \vspace{3pt}\\
 \, \Gamma_t\big|_{t=0}=\Gamma_0,
\end{cases}
\end{equation*}
and $\om ^- (\Gamma)$ denotes the region enclosed by the
hypersurface $\Gamma$ and $\om ^+ (\Gamma)$ the region between
$\dom$ and $\Gamma$.\vsp

In what follows we set
\[
Q_t:=\om\times (0,t)\quad\ \hbox{for}\ \ 0<t\leq T.
\]
Given $\de _0 >0$, we define $\tmaxi=\tmaxi(\de _0)>0$ by
\begin{equation}\label{tmax}
\tmaxi=\max\{t\in[0,T], \Vert v^\ep -\tilde v\Vert_
{L^\infty(Q_t)} \leq \de _0\}.
\end{equation}
The key estimates for proving \eqref{v-est4} are the following:

\begin{clm}\label{claim1-phi}
There exist constants $\de _0 >0$ and $C>0$ such that, for any
$t\in(0,\tmaxi]$, we have
\begin{equation}\label{v1-v2}
\Vert \Phi^0 (v^\ep)-\Phi^0(\tilde v) \Vert_{L^\infty(Q_t)} \leq
C\int_0^t \frac{1}{\sqrt{t-s}}\Vert v^\ep-\tilde v
\Vert_{L^\infty(Q_s)}ds.
\end{equation}
\end{clm}

\begin{clm}\label{claim2-v}
There exists a constant $A>0$ such that, for any $v$ satisfying
the estimates \eqref{v-est1}, \eqref{v-est3} and the Neumann
boundary conditions, and for any $t\in(0,T]$, we have
\begin{equation}\label{Phie-Phi0}
\Vert \Phi^\ep(v)-\Phi^0(v)\Vert _{L^\infty(Q_t)}\leq A\:\!\ep.
\end{equation}
\end{clm}
\vskip 8 pt The proof of these claims will be given later. For the
moment, let us simply mention that Claim \ref{claim2-v} can be
shown by the following two-step argument:  first, our results on
the single equation $\Pe$ yields
\[
\Phi_1^\ep(v)-\Phi_1^0(v)=O(\ep),
\]
in the sense that the transition layer of $\Phi_1^\ep(v)$ and that
of $\Phi_1^0(v)$ are within an $O(\ep)$ distance; this observation
and an estimate of the heat kernel yield \eqref{Phie-Phi0}.  To
prove Claim \ref{claim1-phi}, we also use a similar estimate of
the heat kernel, see Lemma \ref{lm:int-layer} and Subsection
\ref{ss:proof-claim1-phi} for details.

Combining these estimates, we obtain, for any $t\in(0,\tmaxi]$,
\[
\begin{array}{ll}
\Vert \ve - \tilde v\Vert _{L^\infty(Q_t)} &=\Vert \Phi^\ep(\ve) -
\Phi^0(\tilde v)\Vert _{L^\infty(Q_t)}
\vspace{4pt}\\
&\leq \Vert \Phi^\ep(\ve)- \Phi^0(\ve)\Vert _{L^\infty(Q_t)}
 +\Vert \Phi^0(\ve)-\Phi^0(\tilde v)\Vert _{L^\infty(Q_t)}
 \vspace{4pt}\\
 &\leq\di A\:\!\ep+C\int_0^t \frac{1}{\sqrt{t-s}}\Vert \ve -
 \tilde v \Vert_{L^\infty(Q_s)}ds.
 \end{array}
\]
As we will see later in Lemma \ref{lem:k-kbar}, this implies
\begin{equation}\label{bound-kbar}
\Vert \ve - \tilde v\Vert _{L^\infty(Q_t)} \leq A\:\!\ep \bar k(t)
\quad \ \ (0\leq t\leq \tmaxi),
\end{equation}
where $\bar k$ is the function determined by the integral
equality:
\begin{equation}\label{kbar}
\bar k(t)=1+C\int_0^t \frac{\bar k(s)}{\sqrt{t-s}}\,ds.
\end{equation}
Since $\bar k(t)$ is bounded on any finite interval $[0,T]$, we
obtain $\Vert \ve - \tilde v\Vert _{L^\infty(Q_t)}=O(\ep)$, for
$0\leq t \leq \tmaxi$. This implies, first of all, that $\tmaxi=T$
if $\ep$ is small enough, hence it proves \eqref{v-est4}, for
$0\leq t\leq T$. \qed

\vskip 8pt The rest of this section gives a detailed account of
the proof of \eqref{v-est4}. We begin with some notations to
clarify the statements of the above claims.

\subsection{Some notations}\label{ss:notations}

Given any function $\barg (x,t,u)$ satisfying the conditions
\eqref{g-est1} and \eqref{g-est3}, we can define a classical
solution of the interface equation $\Pz$ on some time interval
$0\leq t<\tblow(\barg)$. We denote this solution by
$\Gamma_t\gbragbar$ in order to clarify its dependence on $\barg$.
More precisely, $\Gamma_t\gbragbar$ is a solution of the problem
\[
 (P ^0 _ {\barg}) \quad\begin{cases}\di\vspace{3pt}
 \, V_{n}=-(N-1)\kappa + c_0 \int_{\am}^{\ap}\barg(x,t,r)dr
 \quad \text { on } \Gamma_t\\
 \, \Gamma_t\big|_{t=0}=\Gamma _0.
\end{cases}
\]
Also, we denote by $u^\ep \gbragbar(x,t)$ the solution of the
problem
\[
 (P ^\ep _{\barg}) \quad\begin{cases}
 u_t=\Delta u+\edeux (f(u)-\ep \barg(x,t,u))
 &\text{in }\om \times (0,+\infty)\vspace{3pt}\\
 \dudn = 0 &\text{on }\partial \om \times (0,+\infty)\vspace{3pt}\\
 u(x,0)=u_0(x) &\text{in }\om.
 \end{cases}
\]
Once the interface $\Gamma_t \gbragbar$ is given, we denote by
$\Omega^-_t\gbragbar,\,\Omega^+_t\gbragbar$ the region enclosed by
$\Gamma_t\gbragbar$ and the one enclosed between $\partial\Omega$
and $\Gamma_t\gbragbar$, respectively. As in \eqref{u}, we define
the step function $\tilde{u}\gbragbar(x,t)$ by
\begin{equation}\label{u-g}
 \tilde u\gbragbar(x,t)=\begin{cases}
 \, \ap &\text{in}\ \; \om^+_t\gbragbar\\
 \, \am &\text{in}\ \; \om^-_t\gbragbar
 \end{cases} \quad\text{for } t\in[0,\tblow(\barg)).
\end{equation}

Next, given any function $u(x,t)$ on $\ombar\times[0,T]$, we
denote by $V[u](x,t)$ the solution of the problem
\begin{equation}\label{eqV}
\begin{cases}
\,V_t=D \Delta V + h(u(x,t),V)
\quad &\textrm{in}\ \ \Omega\times(0,T] \vspace{4pt}\\
\,\di\frac{\partial V}{\partial \nu}= 0 \quad&
\textrm{on}\ \ \partial \Omega \times (0,T]\vspace{4pt}\\
\,V(x,0)=v_0 (x) \quad&\textrm{in}\ \ \Omega.
\end{cases}
\end{equation}

In view of \eqref{ge-g} and the above notations, the solution
$(\ue,\ve)$ of $\RD$ is expressed as
\[
\ue=\ue[\;\!\g[\ve]],\quad \ve=V[\;\!\ue].
\]
On the other hand the solution $(\Gamma,\tilde v)$ of $(RD ^0)$ is
expressed as
\[
\Gamma _t=\Gamma _t[\;\!g[\tilde v]],\quad \tilde v=V[\;\!\tilde
u],
\]
the step function $\tilde{u}$ in $(RD ^0)$ being given by
\[
\tilde{u}=\tilde{u}[\;\!g[\tilde{v}]].
\]
Finally, the maps $\Phi ^\ep _1$, $\Phi ^0 _1$ and $\Phi _2$ are
now written as
$$
\Phi ^\ep _1: v \rightarrow \ue[\g[v]], \;\;\; \Phi ^0 _1: v
\rightarrow \tilde u [g[v]],\;\;\; \Phi _2: u \rightarrow V[u].
$$

\subsection{Interface motion under perturbation}
\label{ss:interface}

In this subsection we show that the interface $\Gamma_t\gbragbar$
depends continuously on the pressure term induced by $\barg$. To
this end, we first fix constants $C_*,\,T_*>0$, $\vartheta\in
(0,1)$, and denote by ${\cal Y}$ the set of functions
$\barg(x,t,u)$ on $\ombar\times [0,\,T_*]\times\R$ satisfying
\begin{equation}\label{g-est3-2}
 \sup_{u\in\R}\,\Vert \barg(\cdot,\cdot,u)
 \Vert_{C^{1+\vartheta,\frac{1+\vartheta}{2}}(\ombar\times[0,T_*])} \leq C_*.
\end{equation}

\begin{prop}\label{pr:gamma-g}
Let $\barg\in {\mathcal Y}$. Let $T \in(0,\tblow(\barg))$. Then
there exist positive constants $\delta,\,K,\,M$ such that, for any
$\tildeg \in{\cal Y}$ satisfying
\begin{equation}
\Vert \tildeg-\barg\Vert_{L^\infty(\om\times(0,T)\times\R)}
 \leq\delta,
\end{equation}
it holds that $\tblow(\tildeg)>T$, where we recall that
$\tblow(\tildeg)$ is the maximum time of existence of a classical
solution of Problem $(P^{\;\!0}_{\tildeg})$. Furthermore, for each
$t\in [0,T]$,
\begin{equation}\label{gamma-lip}
 d_{\mathcal H}
 (\,\Gamma_t[\tildeg], \Gamma_t\gbragbar\,)\leq
 K(e^{Mt}-1)\,\Vert
 \tildeg-\barg\Vert_{L^\infty(\om\times(0,t)\times\R)}\,.
\end{equation}
\end{prop}

{\noindent \bf Proof.} First, the assertion that $\tblow(\tilde
g)>T$, for $\tildeg$ sufficiently close to $\barg$, follows from
the standard local existence theory for quasi-linear parabolic
equations. In fact, by using appropriate parametrization, one can
express $\Gamma _ t[\barg]$ and $\Gamma _t[\tildeg]$, as graphes
over ${\cal M}$, where ${\cal M}$ is a $N-1$ dimensional manifold
without boundary, and transfer the motion equations $\Pzggbar$ and
$(P^{\;\!0}_{\tildeg})$, into quasi-linear parabolic equations on
the manifold ${\cal M}$, at least locally in time. For more
details we refer to \cite{CR}. Since $\tildeg$ and $\barg$ satisfy
\eqref{g-est3-2}, and since the embedding
$$
C^{1+\vartheta,\frac{1+\vartheta}{2}}\hookrightarrow
C^{1+\vartheta',\frac{1+\vartheta'}{2}}
$$
is compact if $0<\vartheta'<\vartheta$, the assumption $\Vert
\tildeg -\barg\Vert_{L^\infty} \leq \de$ implies
\[
\Vert \tildeg(\cdot,\cdot,u)-\barg(\cdot,\cdot,u)
\Vert_{C^{1+\vartheta',\frac{1+\vartheta'}{2}}} \leq c(\de),
\]
where $c(\de)$ is a constant satisfying $c(\de)\to 0$, as $\de \to
0$. Consequently, the coefficients appearing in $\Pzggbar$ and
$(P^{\;\!0}_{\tildeg})$ satisfy
\[
 \big\Vert \int^{\ap}_{\am} \tildeg(\cdot,\cdot,r)dr -
 \int^{\ap}_{\am}\barg(\cdot,\cdot,r) dr \,
 \big\Vert_{C^{1+\vartheta',\frac{1+\vartheta'}{2}}} \leq (\ap-\am)c(\de).
\]
Hence, the two solutions $\Gamma _ t[\barg]$ and $\Gamma
_t[\tildeg]$ stay close to each other, at least locally in time,
and, by repeating this argument, one can prove that
$\tblow(\tildeg)>T$, for $\de$ sufficiently small.

Next we prove the estimate \eqref{gamma-lip}. This will be done by
using the maximum principle. Let us introduce some notation. For
each $\barg\in {\cal Y}$, we denote by $d(x,t;\barg)$ the signed
distance function associated with the interface
$\Gamma_t\gbragbar$. By $\bar {\Gamma}_t\preceq\tilde{\Gamma}_t$
we mean that $\bar{\Gamma}_t$ lies inside of $\tilde{\Gamma}_t$.
Clearly we have
\begin{equation}\label{gamma-d}
\Gamma_t\gbragbar \preceq\Gamma_t\gbragtilde\;
\Longleftrightarrow\; d(x,t;\barg)\geq
d(x,t;\tildeg)\quad\hbox{for}\ x\in\ombar.
\end{equation}
Now we choose $t_0\in[0,T]$ arbitrarily and put
\[
\eta_0:=\Vert \tildeg-\barg\Vert_{L^\infty(\om\times(0,t_0)\times
\R)}.
\]
Then
\[
 \barg(x,t,u)-\eta_0 \leq \tildeg(x,t,u)\leq \barg(x,t,u) +\eta_0
 \qquad\hbox{for}\ \ 0\leq t\leq t_0.
\]
The comparison principle then yields
\[
 \Gamma_t[\;\!\barg-\eta_0]\preceq\Gamma_t
 \gbragtilde \preceq\Gamma_t[\;\!\barg+\eta_0]
 \quad \hbox{for}\ \ 0\leq t\leq t_0.
\]
Thus, in order to prove \eqref{gamma-lip}, it suffices to show
that there exist constants $K,M>0$ such that, for all small
$\eta_0>0$,
\begin{equation}\label{gamma-lip-2}
\left\{\begin{array}{l}
 d_{\mathcal H}(\Gamma_t[\;\!\barg-\eta_0],\,
 \Gamma_t\gbragbar)\,\leq\, K\eta_0\;\!(e^{Mt}-1)\vsp\\
 d_{\mathcal H}(\Gamma_t[\;\!\barg+\eta_0],\,
 \Gamma_t\gbragbar)
 \,\leq\, K\eta_0\;\!(e^{Mt}-1),
\end{array}\right.
\end{equation}
for $0\leq t\leq t_0$. We will only show the latter inequality for
$\Gamma_t[\;\!\barg+\eta_0]$ since the former can be shown in the
same manner.

Recall that $d(x,t;\barg)$ satisfies the equation
\eqref{interface}, namely
\begin{equation}\label{interface-1}
 d_t=\Delta d  -c_0 \int^{\ap}_{\am}\barg(x,t,r) dr \quad\
 \textrm{on}\ \; \Gamma_t\gbragbar.
 \end{equation}
Choose a constant $d_0>0$ such that $d(x,t;\barg)$ is smooth
--- say, $C^3$ in $x$ and $C^{3/2}$ in $t$--- in the neighborhood
${\cal N}_{d_0}(\Gamma_t\gbragbar),\;0\leq t\leq T$. By the
smoothness of $d(x,t;\barg)$ and equality \eqref{interface-1},
there exists a constant $N>0$ such that
\[
 \big| d_t-\Delta d  +c_0 \int^{\ap}_{\am}\barg(x,t,r) dr \big|
 \leq N|d| \quad\ \textrm{in}\ \; {\cal N}_{d_0/2}(\Gamma_t\gbragbar).
\]
Now we put
\[
\begin{split}
\artificial (x,t):= d(x,t;\barg)-K\eta_0\;\!(e^{2Nt}-1),\\
\tilde{\Gamma}_t:=\{x\in\om\,|\,\artificial(x,t)=0\},
\end{split}
\]
where the constant $K$ is to be determined later. If
\[
\eta_0\leq \eta_0^*:=\frac{e^{-2NT}\,d_0}{4}K^{-1},
\]
then $\tilde{\Gamma}_t$ lies within the neighborhood ${\cal
N}_{d_0/2}(\Gamma_t\gbragbar)$. Observe that
\[
\begin{split}
(\artificial) _t-\Delta\artificial &=d_t-2NK\eta_0\;\!
e^{2Nt}-\Delta d\\
&\leq -c_0\int^{\ap}_{\am}\barg(x,t,r) dr
+N|d|-2NK\eta_0\;\!e^{2Nt}.
\end{split}
\]
Since $d=K\eta_0\;\!(e^{2Nt}-1)$ on $\tilde{\Gamma}_t$, we obtain
$$
(\artificial)_t-\Delta\artificial \leq
-c_0\di\int^{\ap}_{\am}\barg(x,t,r)dr
 -NK\eta_0 \quad\ \textrm{on}\
 \;\tilde{\Gamma}_t.
$$
Now we set
\[
K=(\ap-\am)c_0 N^{-1}.
\]
Then it follows from the above inequality that
\[
(\artificial) _t\leq
\Delta\artificial-c_0\int^{\ap}_{\am}\barg(x,t,r)dr
 -(\ap-\am)\:\!c_0\:\! \eta_0
 \quad\ \textrm{on}\ \; \tilde{\Gamma}_t.
\]
This inequality and the fact that $\artificial
(x,0)=d(x,0\,;\barg)$ imply that $\tilde{\Gamma}_t$ satisfies
\[
 \begin{cases}\di
 \, V_{n}\geq -(N-1)\kappa + c_0 \int^{\ap}_{\am}
 \big(\barg(x,t,r)+\eta_0\big)\,dr
 \quad \text{on}\ \;\tilde{\Gamma}_t, \vspace{3pt}\\
 \, \tilde{\Gamma}_t\big|_{t=0}= \Gamma _0.
\end{cases}
\]
On the other hand, $\Gamma_t[\barg+\eta_0]$ satisfies
\[
 \begin{cases}\di
 \, V_{n}= -(N-1)\kappa + c_0 \int^{\ap}_{\am}
 \big(\barg(x,t,r)+\eta_0\big)\,dr
 \quad \text{on}\ \;\Gamma_t[\barg+\eta_0], \vspace{3pt}\\
 \, \Gamma_t[\barg+\eta_0]\big|_{t=0}= \gzero.
\end{cases}
\]
By the comparison principle, we obtain
\[
\Gamma_t\gbragbar\preceq \Gamma_t[\barg+\eta_0]\preceq
\tilde{\Gamma}_t \qquad\hbox{for}\ \; 0\leq t\leq t_0.
\]
Consequently,
\[
d_{\mathcal H}(\Gamma_t[\barg+\eta_0],\,\Gamma_t\gbragbar) \leq
d_{\mathcal H}(\tilde{\Gamma}_t,\,\Gamma_t\gbragbar) \leq
K\eta_0\;\!(e^{2Nt}-1),
\]
for $0\leq t\leq t_0$.  The lemma is proved. \qed

\subsection{Proof of Claim \ref{claim2-v}}\label{preuveclaim2}

For a function $v(x,t)$ satisfying the estimates \eqref{v-est1},
\eqref{v-est3} and the Neumann boundary conditions, we compare
below $\Phi ^\ep (v)=\Phi _2 \circ \Phi ^\ep _1 (v)$ and $\Phi ^0
(v)=\Phi _2 \circ \Phi ^0
_1(v)$.\\

\noindent{\bf Action of $\Phi ^\ep _1$ and $\Phi ^0 _1$.} Let us
compare $\Phi ^\ep _1 (v)=\ue[\g[v]]$ with the step function $\Phi
^0 _1(v)=\tilde{u}[g[v]]$. By the definitions in \eqref{ge-g} we
have $\g[v]=g[v]+O(\ep)$, and all the conditions in
\eqref{g-est1}--\eqref{g-neumann} are satisfied. It follows that
our results for the single equation apply and, in particular,
\[
 u_\ep^-(x,t) \leq \ue[\g[v]](x,t+t^\ep)\leq
 u_\ep^+(x,t) \quad \text { for } 0 \leq t \leq
 T-t^\ep,
\]
where $u_\ep^{\pm}$ are as in \eqref{sub}, $d$ being the signed
distance function associated with the interface $\Gamma _t
[g[v]]$. Since the term $\EB$ in $q(t)$ --- that appears in
\eqref{sub}--- quickly becomes small,
\[
 \begin{array}{l}\di
 \Big|\ue[\g[v]]-\tilde{u}[g[v]]\Big|(x,t)
 \leq \ap-U_0\Big(\frac{d(x,t)-\ep p(t)}{\ep}\Big)
 +O(\ep)\quad\hbox{in}\  \om^+_t[\;\!g[v]]
 \vspace{4pt}\\ \di
  \Big|\ue[\g[v]]-\tilde{u}[g[v]]\Big|(x,t)
 \leq U_0\Big(\frac{d(x,t)+\ep p(t)}{\ep}\Big)-\am
 +O(\ep) \quad\hbox{in}\  \om^-_t[\;\!g[v]]
 \end{array}
\]
for $\mu_1\ep^2|\ln \ep| \leq t\leq T$, provided that we choose
the constant $\mu_1$ large enough. Consequently, by Lemma
\ref{est-phi}, there exist constants $B,\,C>0$ such that
\begin{equation}\label{difference}
\big|\Phi ^\ep _1(v)-\Phi ^0 _1 (v)\big|(x,t)\leq B\exp
\Big(-\lambda \frac{|d(x,t)|}{\ep}\Big)+C\ep,
\end{equation}
for $(x,t)\in\ombar\times[\mu_1\ep^2|\ln
\ep|,\;T]$.\\

\noindent{\bf Action of $\Phi _2$.} Next we compare $\Phi ^\ep
(v)=V[\Phi ^\ep _1(v)]$ and $\Phi ^0 (v)=V[\Phi ^0 _1 (v)]$. Set
$w:=\Phi ^\ep  (v)-\Phi ^0 (v)$. By subtracting the equations for
$V[\Phi ^\ep _1(v)]$ and $V[\Phi ^0 _1 (v)]$, we obtain
\[ w_t=D\Delta w +
\Big(h(\phiunep,\phiep)-h(\phiunzero,\phizero)\Big).
\]
Since $|h(\phiunep,\phiep)-h(\phiunzero,\phizero)|\leq
C|w|+C|\phiunep-\phiunzero|$ for some constant $C>0$, the function
$\tilde{w}:=e^{-Ct}w$ satisfies
\[
\tilde{w}_t \leq D\Delta \tilde{w} + Ce^{-Ct}
|\phiunep(x,t)-\phiunzero(x,t)| + C(|\tilde{w}|-\tilde{w}),
\]
hence\vspace{-6pt}
\begin{equation}\label{tilde-w-sub}
\tilde{w}_t \leq D\Delta \tilde{w} +
C|\phiunep(x,t)-\phiunzero(x,t)| + C(|\tilde{w}|-\tilde{w}).
\end{equation}
Now let $W(x,t)$ be the solution of the equation
\[
W_t= D\Delta W + C|\phiunep(x,t)-\phiunzero(x,t)| + C(|W|-W),
\]
with initial data $W(x,0)=0$.  Then since \eqref{tilde-w-sub}
implies that $\tilde{w}$ is a sub-solution of the above equation,
and since $\tilde{w}(x,0)=0$, we have
\begin{equation}\label{tilde-w-W}
\tilde{w}(x,t)\leq W(x,t) \qquad \hbox{for}\ \ x\in\ombar, \ t\geq
0.
\end{equation}
Moreover, since $W\geq 0$, the above equation for $W$ can be
reduced to
\[
W_t= D\Delta W + C|\phiunep(x,t)-\phiunzero(x,t)|.
\]
In view of this, we see that
$$
W(x,t)=C\int_0^t\int_\om G(x,y,t-s)
 |\phiunep(y,s)-\phiunzero(y,s)|\,dy\;\!ds,
$$
$G(x,y,t)$ being the fundamental solution that appears in
\eqref{ve-G}. This and \eqref{tilde-w-W} yield
\begin{equation}\label{w-G0}
 |w(x,t)|\leq C e^{Ct}\int_0^t\int_\om G(x,y,t-s)
 |\phiunep(y,s)-\phiunzero(y,s)|\,dy\;\!ds.
\end{equation}
Combining this and \eqref{difference}, we obtain
\begin{equation}\label{w-G}
 |w(x,t)|\leq
 BCe^{Ct}\int_0^t\int_\om G(x,y,t-s)
 \exp\Big(-\lambda \frac{|d(y,s)|}{\ep}\Big)dy\;\!ds  +O(\ep).
\end{equation}
In order to estimate the above integral, we need the following
lemma:

\begin{lem}\label{lm:int-layer}
Let $\Gamma$ be a smooth closed hypersurface in $\Omega$ and
denote by $d(x)$ the signed distance function associated with
$\Gamma$. Then there exist constants $C, r_0>0$ such that for any
function $\eta(r)\geq 0$ on $\R$, it holds that
\begin{equation}\label{int-layer1}
 \int_{|d|\leq r_0} G(x,y,t)\;\!\eta(d(y))\,dy \leq
 \frac{C}{\sqrt{t}\;}\int_{-r_0}^{r_0}\eta(r)\,dr
 \quad\ \hbox{for}\ \ 0< t\leq T.
\end{equation}
\end{lem}

The proof of this lemma will be given in the next subsection. As
is easily seen from its proof, the above estimate remains to hold
if $\Gamma$ depends on $t$ smoothly; in other words, the constant
$C$ can be chosen uniformly as $\Gamma$ varies. Applying the above
estimate to $\Gamma_t[\;\!g[v]]$, $0<t\leq T$, we obtain
\[
 \begin{split}
 &\int_\om G(x,y,t-s)
 \exp\Big(-\lambda \frac{|d(y,s)|}{\ep}\Big)dy\\
 &\hspace{50pt}=\int_{|d| < r_0}+\int_{|d| \geq r_0}
  G(x,y,t-s)\exp\Big(-\lambda \frac{|d(y,s)|}{\ep}\Big)dy\\
 &\hspace{50pt}= O\big(\frac{\ep}{\sqrt{t-s}}\,\big)+
 O\big(e^{-\lambda r_0/\ep}\big)\\
 &\hspace{50pt}=O\big(\frac{\ep}{\sqrt{t-s}}\,\big).
 \end{split}
\]
It follows from this and \eqref{w-G} that
\begin{equation}\label{est-w}
 |w(x,t)|=O\Big(\ep \int_0^t\frac{1}{\sqrt{t-s}}\,ds\Big)+O(\ep)=O(\ep),
\end{equation}
which completes the proof of Claim \ref{claim2-v}.\qed

\subsection{Proof of Lemma \ref{lm:int-layer}}

We first show that
\begin{equation}\label{surface-int}
 \int_{\Gamma}G(x,y,t)\,dS_y \leq \frac{C}{\sqrt{t\,}\,} \qquad
 \hbox{for}\ \ x\in\Omega,\ 0< t\leq T.
\end{equation}
It suffices to prove this estimate on a small interval $[0,t_0]$,
since the estimate for the remaining interval $[t_0,T]$ will
follow by simply choosing a large constant $C$ (since $G$ is
bounded for $t$ large). Hereafter, we choose $t_0$ sufficiently
small. Then, for $0<t\leq t_0$, $G(x,y,t)$ is well approximated by
the fundamental solution on the entire space $\R^N$:
\[
G_0(x,y,t):=\frac{1}{(4\pi D t)^{N/2}}\,
\exp\Big(-\frac{|x-y|^2}{4Dt}\Big).
\]
In particular, there exists a constant $C>0$ such that
\[
0<G(x,y,t)\leq C \,G_0(x,y,t) \qquad\hbox{for}\ \ x,y\in\ombar, \
0<t\leq t_0,
\]
(see, for example, \cite[Section IV.2]{Ey}).  Thus it suffices to
prove \eqref{surface-int} for $G_0$ instead of $G$.

Given $x\in\Omega$, let $x_0$ be the point on $\Gamma$ that is
closest to $x$, and let $n(x_0)$ be the outward normal to $\Gamma$
at $x_0$.  Then $x-x_0=d(x)n(x_0)$.  Define
\[
\widetilde{Y}:=\{\,y\in \R^N\;,\;y\cdot n(x_0)=0\,\}, \qquad
Y_0=span\langle n(x_0)\rangle,
\]
where $\cdot$ denotes the Euclidean inner product in $\R^N$ and
$span\langle w\rangle$ the line spanned by the vector $w$. This
gives an orthogonal decomposition $\R^N=\widetilde Y\oplus Y_0$,
and $x_0+\widetilde Y$ is the tangent hyperplane of $\Gamma$ at
$x_0$. Since $\Gamma$ is smooth, it is expressed locally as the
graph of a map defined on a subset of $\widetilde Y$.  More
precisely, there exist a smooth map $h:\widetilde Y\to Y_0$ and a
constant $\delta>0$ such that $h(0)=0,\;\nabla h(0)=0$, and that
\begin{equation}\label{S-Gamma}
\begin{array}{c}
S:=\{\,x_0+\tilde y+h(\tilde y)\;,\;\tilde y\in\widetilde Y,
|\tilde y|<\delta\,\}\subset \Gamma,\vsp\\
dist(x_0,\Gamma\setminus S)\geq \delta.\qquad\qquad
\end{array}
\end{equation}
Now we decompose the integral \eqref{surface-int} for $G_0$ as
\[
 \int_{\Gamma}G_0(x,y,t)\,dS_y =
 \frac{1}{(4\pi D t)^{N/2}}\,\left(\int_S + \int_{\Gamma\setminus S}
 \exp\Big(-\frac{|x-y|^2}{4Dt}\Big)dS_y\right).
\]
Since $|x-y|\geq |d(x)|$ for every $y\in\Gamma$ and since
\[
|x-y|\geq
\big|\,|x-x_0|-|y-x_0|\,\big|=\big|\,|d(x)|-|y-x_0|\,\big|,
\]
we have
\[
|x-y|\geq \frac{|d(x)|+\big|\,|d(x)|-|y-x_0|\,\big|}{2}\geq
\frac{|y-x_0|}{2}.
\]
This and \eqref{S-Gamma} yield
\[
|x-y|\geq \frac{\delta}{2}\qquad\hbox{for}\ \ y\in\Gamma\setminus
S.
\]
Consequently,
\begin{equation}\label{int-GS}
\int_{\Gamma\setminus S}\exp\Big(-\frac{|x-y|^2}{4Dt}\Big)dS_y
\leq e^{-\delta^2/16Dt}\,|\Gamma|,
\end{equation}
where $|\Gamma|$ denotes the total area of $\Gamma$.

On the other hand, for each $y\in S$, we can express $y-x_0$ as
\[
y-x_0=\tilde y + h(\tilde y)\qquad (\tilde y\in\widetilde
Y,\;h(\tilde y)\in Y_0),
\]
and $\widetilde Y$ can be identified with $\R^{N-1}$.  Thus
\[
\begin{split}
\int_{S}\exp\Big(-\frac{|x-y|^2}{4Dt}\Big)dS_y\hspace{200pt}\\
=\int_{|\tilde y|<\delta} \exp\Big(-\frac{|x-x_0-\tilde y-h(\tilde
y)|^2}{4Dt}\Big) \sqrt{1+|\nabla h(\tilde y)|^2}\,d\tilde y.
\end{split}
\]
Since $\nabla h(0)=0$, there exists a constant $C_1>0$ such that
\begin{equation}\label{nabla-h}
|\nabla h(\tilde y)|\leq C_1|\tilde y|\qquad \hbox{for}\ \ |
\tilde y|<\delta.
\end{equation}
Note also that the orthogonality $(x-x_0-h(\tilde y))\,\bot\,
\tilde y$ implies
\[
|x-x_0-\tilde y-h(\tilde y)|^2= |x-x_0-h(\tilde y)|^2+|\tilde
y|^2\geq |\tilde y|^2.
\]
Combining these, we obtain
\[
\begin{array}{rl}\displaystyle
\int_{S}\exp\Big(-\frac{|x-y|^2}{4Dt}\Big)dS_y&
\displaystyle\vsp\leq \int_{|\tilde
y|<\delta}\exp\Big(-\frac{|\tilde y|^2}{4Dt}\Big)
\sqrt{1+C_1^2|\tilde y|^2}\,d\tilde y\\
&\displaystyle= t^{(N-1)/2}
\int_{|z|<\sqrt{t}^{-1}\delta}\,e^{-|z|^2/4D}
\sqrt{1+t\,C_1^2|z|^2}\,dz,
\end{array}
\]
where $z:=\tilde y/\sqrt{t}$.  Observe that, as $t\to 0$,
\[
\int_{|z|<\sqrt{t}^{-1}\delta}\,e^{-|z|^2/4D}
\sqrt{1+t\,C_1^2|z|^2}\,dz \:\to\:
\int_{\R^{N-1}}\,e^{-|z|^2/4D}\,dz = (4D\pi)^{(N-1)/2}.
\]
Consequently,
\[
 \frac{1}{(4\pi D t)^{N/2}}\,\int_S
 \exp\Big(-\frac{|x-y|^2}{4Dt}\Big)dS_y \leq
 \frac{1}{\sqrt{4\pi Dt}}+ o\Big(\frac{1}{\sqrt t}\;\Big).
\]

Combining the estimate above and \eqref{int-GS}, we obtain
\[
 \int_{\Gamma}G_0(x,y,t)\,dS_y =
 O\Big(\frac{1}{\sqrt{t}}\Big)+O\Big(\frac{1}{(\sqrt{t\,}\,)^N}
 e^{-\delta^2/16Dt}\Big)=O\Big(\frac{1}{\sqrt{t}\,}\Big).
\]
Since $\Gamma$ is a smooth compact hypersurface, its curvature is
bounded.  Therefore, the constants $\delta$ and $C_1$ that appear
in \eqref{int-GS}, \eqref{nabla-h} can be chosen independent of
the choice of $x_0\in\Gamma$.  Hence the above $O(1/\sqrt{t})$
estimate is uniform with respect to the choice of $x\in\om$.  This
proves the estimate \eqref{surface-int}.

Now, choose a sufficiently small constant $r_0>0$ such that the
signed distance function $d(x)$ is smooth in the region
$\{\,d(x)<2r_0\,\}$.  For each $r\in[-r_0,r_0]$, we define a
hypersurface $\Gamma(r)$ by
\[
\Gamma(r):=\{\,x\in\om,\;d(x)=r\,\}.
\]
Then the curvatures of $\Gamma(r)$ are uniformly bounded as $r$
varies, which implies that there exists some constant $C>0$ such
that
\[
\int_{\Gamma(r)}G(x,y,t)\,dS_y \leq \frac{C}{\sqrt{t\,}\,}
\qquad\hbox{for}\ \ 0< t\leq T, \ \;r\in[-r_0,r_0].
\]
The estimate \eqref{int-layer1} now follows by integrating in $r$.
\qed

\subsection{Proof of Claim
\ref{claim1-phi}}\label{ss:proof-claim1-phi}

We compare below $\Phi^0 (v^\ep)=\Phi _2 \circ \Phi ^0 _1 (v^\ep)$
and $\Phi ^0 (\tilde v)=\Phi _2
\circ \Phi ^0 _1(\tilde v)$.\\

\noindent{\bf Action of $\Phi ^0 _1$.} Let us compare the two step
functions $\Phi ^0 _1 (v^\ep)=\tilde u[g[v^\ep]]$ and $\Phi ^0 _1
(\tilde v)=\tilde u[g[\tilde v]]$. We want to apply Proposition
\ref{pr:gamma-g}, with $g[\tilde v]$ and $g[v^\ep]$, playing the
role of $\barg$ and $\tildeg$, respectively. (Hence, the role of
$\tblow(\barg)$ is played by $\tblow(g[\tilde v])$, which
corresponds to $\tblow$ in Lemma \ref{existence}). First, we
choose $C_*>0$ large enough so that both $g[v^\ep](x,t,u)$ and
$g[\tilde v](x,t,u)$ satisfy \eqref{g-est3-2}. For
$T\in(0,\tblow)$, we choose $\de$, $K$ and $M$ as in Proposition
\ref{pr:gamma-g}. Next, we define $K_1=\max_{(u,v)\in{\cal
R}}|\partial_v f_1(u,v)|$, with ${\cal R}$ being the rectangle
defined in Subsection \ref{ss:global}, and $\de _0=\de /K_1$.  We
observe that, using the definition of $\tmaxi$ in \eqref{tmax},
\[
\Vert g[v^\ep]-g[\tilde v]\Vert_{L^\infty(\om\times
(0,\tmaxi)\times\R)}\leq K_1\Vert v^\ep-\tilde
v\Vert_{L^\infty(Q_{\tmaxi})}\leq K_1\de _0=\de.
\]
By \eqref{gamma-lip}, it follows that, for any $t\in[0,\tmaxi]$,
\[
 d_{\mathcal H}
 (\,\Gamma_t[\;\!g[v^\ep]],\;\!\Gamma_t[g[\tilde v]]\,)\leq K(e^{Mt}-1)\,
 \Vert g[v^\ep]-g[\tilde v]\Vert_{L^\infty}.
\]
Combining these, we obtain
\begin{equation}\label{gamma-v}
 d_{\mathcal H}
 (\,\Gamma_t[\;\!g[v^\ep]],\;\!\Gamma_t[g[\tilde v]]\,)\leq
 KK_1(e^{Mt}-1)\,\Vert v^\ep-\tilde v\Vert_{L^\infty(Q_t)}\,.
\end{equation}

\noindent{\bf Action of $\Phi _2$.} Next we compare the two
functions $\Phi ^0 (v^\ep)=V[\Phi ^0 _1(v^\ep)]=V[\tilde u
[g[v^\ep]]]$ and $\Phi ^0 (\tilde v)=V[\Phi ^0 _1(\tilde
v)]=V[\tilde u [g[\tilde v]]]$. Since
\[
|\tilde{u}[g[v^\ep]]-\tilde{u}[g[\tilde v]]|\leq \ap -\am,
\]
and since the two step functions differ only in the region
enclosed between the two surfaces $\Gamma_t[\;\!g[v^\ep]]$ and
$\Gamma_t[\;\!g[\tilde v]]$, the estimates \eqref{w-G0} and
\eqref{int-layer1} imply that there exists a constant $B_1>0$ such
that
\[
 \Vert V[\tilde{u}[g[v^\ep]]]-V[\tilde u [g[\tilde v]]]\;\!
 \Vert_{L^\infty(Q_t)}\leq B_1\int_0^t\frac{d_{\mathcal H}
 (\,\Gamma_s[\;\!g[v^\ep]],\;\!\Gamma_s[g[\tilde v])}{\sqrt{t-s}\;}\,ds.
\]
Combining this and \eqref{gamma-v}, we obtain, for any
$t\in[0,\tmaxi]$,
\begin{equation}\label{Phi-Lip}
 \Vert \;\!\Phi ^0(v^\ep)-\Phi^0 (\tilde v)\;\!\Vert_{L^\infty(Q_t)}
 \leq C \int_0^t \frac{1}{\sqrt{t-s}}\Vert v^\ep-\tilde v
\Vert_{L^\infty(Q_s)}ds,
 \end{equation}
with $C=B_1 K K_1 (e^{MT}-1)$. The proof of Claim \ref{claim1-phi}
is complete.\qed

\subsection{Estimate of $\,\bar k(t)$}

In this subsection we justify the estimate \eqref{bound-kbar}. Let
$\bar k(t)$ be the function satisfying \eqref{kbar}, namely,
\begin{equation*}
 \bar k(t)=1+C\int_0^t \frac{\bar k(s)}{\sqrt{t-s}}\,ds
 \quad\ \ \hbox{for}\ \ t\geq 0.
\end{equation*}
We will show below that $\bar k$ is given by
\begin{equation}\label{kbar2}
 \bar k(t)=e^{C^2\pi t}\Big(1+C\int_0^t \frac{e^{-C^2\pi s}}
 {\sqrt{s}}\,ds\Big).
\end{equation}
The following lemma justifies \eqref{bound-kbar}:

\begin{lem}\label{lem:k-kbar}
Let $k(t)$ be a continuous function satisfying
\[
0\leq k(t)\leq A+C\int_0^t \frac{k(s)}{\sqrt{t-s}}\,ds
\qquad(0<t\leq T),
\]
for some constant $A>0$ and $T>0$.  Then
\begin{equation}\label{k-kbar}
0\leq k(t)\leq A \bar k(t)\quad\ \ \hbox{for}\ \ 0<t\leq T.
\end{equation}
\end{lem}

{\noindent \bf Proof.} Define
\[
\bar k_\ep(t):=(1+\ep)A \bar k(t).
\]
Then this function satisfies
\begin{equation}\label{k-ep}
 \bar k_\ep(t)=(1+\ep)A+C\int_0^t \frac{\bar k_\ep(s)}{\sqrt{t-s}}
 \,ds \qquad(0<t\leq T).
\end{equation}
In particular, we have $\bar k_\ep(0)=(1+\ep)A>A\geq k(0)$. Let us
show that
\begin{equation}\label{k-kep}
  k(t)< \bar k_\ep(t)\quad\ \hbox{for}\ \  0<t\leq T.
\end{equation}
Suppose that \eqref{k-kep} does not hold.  Then there exists
$t_0\in(0,T]$ such that
\begin{equation}\label{k(t0)}
  k(t)<\bar k_\ep(t)\quad\hbox{for}\ \ 0\leq t<t_0,\qquad
  k(t_0)=\bar k_\ep(t_0).
\end{equation}
Combining the first part of \eqref{k(t0)} and \eqref{k-ep}, we get
\[
\bar k_\ep(t_0)= (1+\ep)A+C\int_0^{t_0} \frac{\bar
k_\ep(s)}{\sqrt{t_0-s}}
> A+C\int_0^{t_0} \frac{k(s)}{\sqrt{t_0-s}}
\geq k(t_0),
\]
but this contradicts the second part of \eqref{k(t0)},
establishing \eqref{k-kep}.  Letting $\ep\to 0$, we obtain
\eqref{k-kbar}.\qed

\begin{cor}\label{cor:kbar}
Let $k(t)$ be a continuous function satisfying
\[
0\leq k(t)= A+C\int_0^t \frac{k(s)}{\sqrt{t-s}}\,ds \qquad(0<t\leq
T),
\]
for some constant $A>0$ and $T>0$. Then $k(t)\equiv A\bar k(t)$.
In particular, the function $\bar k(t)$ is uniquely determined by
the integral identity \eqref{kbar}.
\end{cor}

{\noindent \bf Proof.} Define $\hat k(t):=A^{-1}k(t)$.  Then $\hat
k(t)$ satisfies
\[
\hat k(t)=1+C\int_0^t \frac{\hat k(s)}{\sqrt{t-s}}\,ds
\qquad(0<t\leq T).
\]
By Lemma \ref{lem:k-kbar} we have $\hat k(t)\leq \bar k(t)$.
Exchanging the role of $\hat k$ and $\bar k$, we obtain the
opposite inequality, hence $\hat k(t)\equiv \bar k(t)$. Thus
$k(t)=A\bar k(t)$.\qed\\

Now let us prove \eqref{kbar2}. Integration by parts gives
\[
\bar k(t)=1+2C\sqrt{t}+2C\int_0^t\sqrt{t-s}\;{\bar k}'(s)\:\!ds.
\]
Hence
\[
\bar k'(t)=\frac{C}{\sqrt{t}}+C\int_0^t\frac{\bar k'(s)}
{\sqrt{t-s}}\,ds.
\]
Thus the function $m(t):=\bar k'(t)-C/\sqrt{t}$ satisfies
\[
m(t)=C\int_0^t \frac{m(s)+C(\sqrt{s}\,)^{-1}}{\sqrt{t-s}}\,ds
=C^2\int_0^t \frac{1}{\sqrt{s(t-s)}}\,ds+ C\int_0^t
\frac{m(s)}{\sqrt{t-s}}\,ds.
\]
Since the first integral on the right-hand side is equal to $\pi$,
we obtain
\[
m(t)=C^2\pi+C\int_0^t \frac{m(s)}{\sqrt{t-s}}\,ds.
\]
It follows from Corollary \ref{cor:kbar} that $m(t)=C^2\pi\bar
k(t)$, hence
\[
\bar k'(t)=C^2\pi\bar k(t)+\frac{C}{\sqrt{t}}.
\]
This and the equality $\bar k(0)=1$ yield \eqref{kbar2}.\qed

\end{document}